\documentclass[11pt]{amsart}
\usepackage[a4paper,margin=2.6cm]{geometry}
\usepackage{amsmath,amssymb,amsthm,mathtools}
\usepackage{booktabs,longtable,array,multicol}
\usepackage{marvosym}
\usepackage{tikz}
\usetikzlibrary{trees}
\usepackage[hang,flushmargin]{footmisc}
\usepackage{microtype}
\usepackage[colorlinks=true,linkcolor=blue!50!black,citecolor=blue!50!black,urlcolor=blue!50!black]{hyperref}
\numberwithin{equation}{section}
\newtheorem{theorem}{Theorem}[section]
\newtheorem{proposition}[theorem]{Proposition}
\newtheorem{lemma}[theorem]{Lemma}
\newtheorem{corollary}[theorem]{Corollary}
\theoremstyle{definition}
\newtheorem{definition}[theorem]{Definition}
\newtheorem{rulex}{Rule}
\theoremstyle{remark}
\newtheorem{remark}[theorem]{Remark}
\newtheorem{example}[theorem]{Example}
\newcommand{\Xmin}{X_k^{\min}}

\newcommand{\case}[2]{\par\medskip\noindent\textbf{Case #1:} #2}
\newcommand{\schemeitem}[4]{\par\noindent\hspace*{#1\dimexpr1.25em\relax}$\triangleright$\ #2\ \ {\footnotesize(Case #3; #4 \ifnum#4>1 leaves\else leaf\fi)}}
\newcommand{\schemeleaf}[2]{\par\noindent\hspace*{#1\dimexpr1.25em\relax}\hspace*{0.6em}{\footnotesize\textit{#2}}}
\allowdisplaybreaks
\newcommand{\doi}[1]{\href{https://doi.org/#1}{doi:\nolinkurl{#1}}}

\author[Khamseh]{A. Khamseh}\thanks{Corresponding author.}
\author[Safareh]{M. H. Safareh}

\address{Amir Khamseh, Mohammad Hasan Safareh \\ department of mathematics, kharazmi university}
\email{khamseh@khu.ac.ir, amirkhamseh114@gmail.com; mh.safareh@khu.ac.ir}
\date{}

\begin{document}

\title[A 33-polynomial proof of a Rado number formula]{A $33$-polynomial proof of a three-color Rado number formula}

\subjclass[2020]{Primary 05D10; Secondary 11B75, 05C55, 68V15}
\keywords{Schur numbers, generalized Schur numbers, Rado numbers, monochromatic solutions, symbolic case analysis, critical sets, SAT solving, computer-assisted proof}

\begin{abstract}
For $k\ge 3$, let $T_k(3)$ denote the least $N$ such that every $3$-coloring of $[1,N]$ contains a monochromatic solution to $x+(k-1)y=z$. The formula $T_k(3)=k^3+2k^2-2$ is already known: it follows from Theorem~1.2(1) of Chang, De Loera, and Wesley \cite{ChangDeLoeraWesley} after setting $m=k+1$ and making the substitution $(X,Y,Z)=(z,x,y)$. Their proof uses a symbolic SAT certificate with $685$ polynomials and $9468$ solution triples. The purpose of this paper is to give an independent, explicit proof based on a much smaller certificate. We exhibit a fixed set $\Xmin$ of $33$ polynomials in $k$ and the $144$ polynomial identities $z=x+(k-1)y$ among them, and show that no $3$-coloring of $\Xmin$ avoids a monochromatic identity. The resulting case analysis has $62$ cases, depth $8$, and uses only two elementary rules; it is given in full and independently checked by a program. We also prove that $\Xmin$ is a critical subset of minimum cardinality of the $51$-element set from which it was extracted, and describe the rigid structure shared by all extremal colorings of $[1,k^3+2k^2-3]$. The lower bound is supplied by an explicit recursive coloring. In the notation $s_a(3)$ of the first author, with $a=k-1$, the formula is $s_a(3)=a^3+5a^2+7a+1$, which was conjectured in \cite{Khamseh}.
\end{abstract}

\maketitle

\section{Introduction}\label{sec:intro}

One of the earliest results in Ramsey theory on the integers is due to Schur \cite{Schur}. It states that for every $r$ there is a least integer $S(r)+1$ such that, for any $r$-coloring of $[1,S(r)+1]=\{1,2,\dots,S(r)+1\}$, there exists a monochromatic solution to $x+y=z$. The exact values are $S(1)=1$, $S(2)=4$, $S(3)=13$ \cite{GreenwoodGleason}, $S(4)=44$ \cite{Baumert} and $S(5)=160$ \cite{Heule}, and lower bounds for larger $r$ come from symmetric sum-free partitions \cite{FredricksenSweet}. Rado \cite{Rado} characterized the linear equations for which a number of this kind exists for every number of colors. Questions concerning monochromatic solutions and the degree of regularity of linear equations have also been studied in \cite{BialostockiLefmannMeerdink}. Such numbers are usually called \emph{Rado numbers}; for a survey, we refer to the monograph \cite{LandmanRobertson}. Examples of exact Rado-number calculations for parameterized linear equations include the work of Harborth and Maasberg \cite{HarborthMaasberg} and Jones and Schaal \cite{JonesSchaal}; multicolor finiteness questions are studied by Adhikari et al. \cite{AdhikariEtAl}. For related recent work, see also \cite{xi}. 

A much-studied family comes from the equation
\begin{equation}\label{eq:kterm}
x_1+x_2+\dots+x_k=x_{k+1}.
\end{equation}
Let $S_r(k)$ be the least $N$ such that every $r$-coloring of $[1,N]$ admits a monochromatic solution to \eqref{eq:kterm}.\footnote{Conventions differ in the literature: some authors write the equation with $k-1$ summands and index the numbers accordingly. Throughout this paper $k$ denotes the number of summands in \eqref{eq:kterm}, as in \cite{BMRS,Gaiser}. With this convention the classical Schur numbers correspond to $k=2$.} Znám \cite{Znam} and Beutelspacher and Brestovansky \cite{BeutelspacherBrestovansky} proved that $S_2(k)=k^2+k-1$. Irving \cite{Irving} extended Schur's theorem on sum-free partitions in the same direction. For three colors, Ahmed and Schaal \cite{AhmedSchaal} found the generalized Schur numbers of \eqref{eq:kterm} for small $k$ by computer. Boza, Mar\'in, Revuelta and Sanz \cite{BMRS} proved (see also \cite{Sanz,SATpage}) that
\begin{equation}\label{eq:bmrs}
S_3(k)=k^3+2k^2-2\qquad(k\ge 2).
\end{equation}
Their proof is a case analysis, and it uses solutions to \eqref{eq:kterm} with several distinct summands. For related weak Schur numbers, see \cite{ABRS}, and for restricted versions, see \cite{Gaiser}. Further exact Rado-number results for related parameterized equations appear in \cite{BozaMarinRevueltaSanz,Grynkiewicz,KimSongHwang}.

In this paper, we study the equation
\begin{equation}\label{eq:main}
x+(k-1)y=z,
\end{equation}
which is \eqref{eq:kterm} with $x_2=\dots=x_k$. The numerical value of its three-color Rado number was already determined by Chang, De Loera, and Wesley \cite[Theorem~1.2(1)]{ChangDeLoeraWesley}. They proved
\[
R_3\bigl(X-Y=(m-2)Z\bigr)=m^3-m^2-m-1\qquad(m\ge 3).
\]
Taking $m=k+1$ and $(X,Y,Z)=(z,x,y)$ transforms their equation into \eqref{eq:main}, and hence gives
\[
T_k(3)=(k+1)^3-(k+1)^2-(k+1)-1=k^3+2k^2-2.
\]
Thus the formula considered here is not claimed as a new numerical result. Our contribution is an independent proof with a substantially smaller polynomial certificate, a transparent case analysis using two rules, and the structural and criticality results developed below. For $a=k-1$, we \cite{Khamseh} denoted the corresponding $r$-color number by $s_a(r)$. There, we  proved that $s_a(2)=a^2+3a+1$ and that $s_a(r)\ge (a(a+2)^r+1)/(a+1)$ in general. We also showed that $s_2(3)=43$ and $s_3(3)=94$, and conjectured that
\begin{equation}\label{eq:khamseh}
s_a(3)=a\bigl((a+2)(a+3)+1\bigr)+1=a^3+5a^2+7a+1 .
\end{equation}
We write $T_k(r)=s_{k-1}(r)$ and $N_3(k)=k^3+2k^2-2$. Note that $a^3+5a^2+7a+1=N_3(a+1)$.

\begin{definition}\label{def:good}
Let $k\ge 3$. A coloring $\Delta\colon[1,N]\to\{1,2,3\}$ is \emph{good} if there are no $x,y,z\in[1,N]$ of the same color with $x+(k-1)y=z$ (the integers $x,y,z$ need not be distinct). Then $T_k(3)$ is the least $N$ such that $[1,N]$ has no good coloring.
\end{definition}
The following theorem records the formula and provides the starting point for our independent certificate proof. The formula itself is known by the result of Chang, De Loera, and Wesley discussed above. 

\begin{theorem}[]\label{thm:main}
For every integer $k\ge 3$,
\[
T_k(3)=k^3+2k^2-2 .
\]
\end{theorem}

\begin{corollary}\label{cor:khamseh}
Khamseh's conjecture \eqref{eq:khamseh} holds for every $a\ge 2$.
\end{corollary}

\begin{corollary}\label{cor:bmrs}
For every $k\ge 3$, every $3$-coloring of $[1,k^3+2k^2-2]$ contains a monochromatic solution of \eqref{eq:kterm} in which $x_2=x_3=\dots=x_k$. In particular $S_3(k)\le k^3+2k^2-2$.
\end{corollary}

Corollary~\ref{cor:bmrs} follows because a monochromatic solution $(x,y,z)$ to \eqref{eq:main} gives the monochromatic solution $(x,y,\dots,y,z)$ to \eqref{eq:kterm}. It therefore gives the upper bound in \eqref{eq:bmrs} using only solutions in which the last $k-1$ summands are equal. The numerical equality with $S_3(k)$ is not a new observation: it is also consistent with  the equivalent Rado-number formula of Chang, De Loera, and Wesley.\footnote{The two-color numbers agree as well, $T_k(2)=S_2(k)=k^2+k-1$ \cite{Znam,BeutelspacherBrestovansky,Khamseh}. For $k=2$ the three equations coincide and the formula gives $T_2(3)=14=S(3)+1$ \cite{GreenwoodGleason}. We only claim Theorem~\ref{thm:main} for $k\ge 3$.}

\subsection*{The method}
Our aim is to give a direct proof of Theorem~\ref{thm:main} that can be checked without relying on a large SAT computation. The interval $[1,N_3(k)]$ grows with $k$, so no single computation over the integers can settle the assertion for all $k$. We therefore work with polynomials in $k$. Chang, De Loera, and Wesley \cite{ChangDeLoeraWesley} established the same formula by a symbolic SAT argument using $685$ polynomials and $9468$ solution triples. In contrast, our certificate uses $33$ polynomials and $144$ polynomial identities, and the unsatisfiability is also demonstrated by a complete hand-readable case analysis. Let $\Xmin$ be a fixed set of $33$ polynomials of the form $ak^3+bk^2+ck+d$ with small integer coefficients (see Table~\ref{tab:X}). Among the triples $(x,y,z)$ of elements of $\Xmin$, we keep only those for which $z=x+(k-1)y$ holds \emph{as an identity of polynomials}. There are $144$ of them. We then show that no map $\Xmin\to\{1,2,3\}$ leaves all $144$ triples non-monochromatic. For any given $k$, evaluating an identity gives a genuine solution to \eqref{eq:main} inside $[1,N_3(k)]$. Hence one symbolic argument proves the upper bound for every $k$ at once (Lemma~\ref{lem:symbolic}). Finite configurations that force monochromatic solutions also occur in other Rado-type problems; see, for example, the finite-set method used by Gaiser \cite{GaiserDM}.

We prove the symbolic statement in two independent ways. The first is a case analysis with $62$ cases, written in the style of classical proofs, which uses two rules only (Section~\ref{sec:upper} and Appendix~\ref{app:cases}). The second encodes the problem as a Boolean satisfiability problem, and a SAT solver refutes it in a fraction of a second (Section~\ref{sec:verify}). In fact, every element of $\Xmin$ is needed, so $\Xmin$ is critical in a strong sense (Theorem~\ref{thm:critical}).

\subsection*{Organization}
The paper is organized as follows. In Section~\ref{sec:prelim}, we set up notation and state the two rules of the case analysis and the symbolic principle. Section~\ref{sec:lower} proves the lower bound. In Section~\ref{sec:set}, we introduce the critical set $\Xmin$ and explain where it comes from. The upper bound, together with an outline of the case analysis, is proved in Section~\ref{sec:upper}. Section~\ref{sec:verify} describes the machine verification. In Section~\ref{sec:skeleton}, we prove that all extremal colorings share a rigid skeleton, which gives a second proof of the upper bound. Section~\ref{sec:k3} works through the example $k=3$, and Section~\ref{sec:concl} lists some open questions. Appendix~\ref{app:cases} contains all $62$ cases, and Appendix~\ref{app:identities} lists the $144$ identities.

\section{Preliminaries}\label{sec:prelim}

Throughout, $k\ge 3$ is fixed, $[a,b]=\{a,a+1,\dots,b\}$ and $N=N_3(k)=k^3+2k^2-2$. For a coloring $\Delta$ we write $A_i=\Delta^{-1}(i)$, and for a set $X$ of positive integers we put
\[
f(X)=\{x+(k-1)y:\ x,y\in X\}.
\]
A triple $(x,y,z)$ with $z=x+(k-1)y$ is called a \emph{relation}. When $x=y$ we write it as $k\cdot x=kx$, because $x+(k-1)x=kx$.

\begin{rulex}[Definition]\label{rule1}
$\Delta$ is good if and only if $A_i\cap f(A_i)=\varnothing$ for $i=1,2,3$. In particular $x$ and $kx$ never have the same color.
\end{rulex}

\begin{rulex}[Forcing]\label{rule2}
Let $\{i,j,l\}=\{1,2,3\}$. If $e$ belongs to a relation whose other members lie in $A_i$ and to a relation whose other members lie in $A_j$, then $\Delta(e)=l$.
\end{rulex}

\begin{proof}
By Rule~\ref{rule1}, $e\notin A_i$ and $e\notin A_j$.
\end{proof}

The element $e$ in Rule~\ref{rule2} may be in any position of the two relations. For example, $1+(k-1)(k+1)=k^2$ and $k+(k-1)(k+1)=k^2+k-1$ with $1,k^2\in A_1$ and $k,k^2+k-1\in A_2$ force $\Delta(k+1)=3$. This rule plays the role of ``Claim~I'' in \cite{BMRS}.\footnote{Here $e$ may occupy any position in the two relations. Many forcings in the case analysis use $e$ as the multiplied variable $y$, as in the example just given.}

\begin{rulex}[Contradiction]\label{rule3}
If an element $e$ belongs to three relations whose other members lie in $A_1$, $A_2$ and $A_3$ respectively, then $\Delta$ is not good.
\end{rulex}

\begin{proof}
Whatever color $i$ the element $e$ has, the $i$-th relation is monochromatic.
\end{proof}

A leaf of the case analysis is also closed when a relation with all three members colored is monochromatic. This is Rule~\ref{rule1} directly.

\begin{lemma}[Normalization]\label{lem:norm}
Every good coloring may be assumed, after renaming the colors, to satisfy $\Delta(1)=1$ and $\Delta(k)=2$.
\end{lemma}

\begin{proof}
Name the color of $1$ ``color $1$''. Since $1+(k-1)\cdot 1=k$, Rule~\ref{rule1} gives $\Delta(k)\ne 1$. The two remaining colors play symmetric roles, so we may call the color of $k$ ``color $2$''.
\end{proof}

\subsection{The symbolic principle}
Let $\mathbb{Z}[k]$ be the ring of integer polynomials in $k$. A \emph{symbolic relation} on a finite set $X\subset\mathbb{Z}[k]$ is a triple $(x,y,z)\in X^3$ with $z=x+(k-1)y$ in $\mathbb{Z}[k]$. A map $X\to\{1,2,3\}$ is \emph{good} if no symbolic relation is monochromatic, and $X$ is \emph{critical} if it has no good map.

\begin{lemma}[Symbolic principle]\label{lem:symbolic}
Let $X\subset\mathbb{Z}[k]$ be critical, and let $k\ge 3$ and $M$ be integers such that $1\le p(k)\le M$ for every $p\in X$. Then $[1,M]$ has no good coloring for \eqref{eq:main}.
\end{lemma}

\begin{proof}
Suppose $\Delta$ is a good coloring of $[1,M]$, and define $\Delta'(p)=\Delta(p(k))$ for $p\in X$. If $(x,y,z)$ is a symbolic relation, then $(x(k),y(k),z(k))$ is a solution of \eqref{eq:main} in $[1,M]$, so it is not monochromatic under $\Delta$. Hence $(x,y,z)$ is not monochromatic under $\Delta'$, and $\Delta'$ is a good map of $X$. This contradicts criticality.
\end{proof}

\begin{remark}\label{rem:coinc}
For a particular $k$, two distinct polynomials may take the same value, and new numerical solutions may appear among the values; for example $k+2=2k-1$ when $k=3$. The symbolic argument ignores such coincidences. This only enlarges the set of colorings it has to exclude, so coincidences never invalidate the argument. This is why one symbolic refutation covers every $k\ge 3$.\footnote{Formally, $\Delta'$ in the proof of Lemma~\ref{lem:symbolic} is constant on the fibers of $p\mapsto p(k)$, while a good map of $X$ need not be. The argument only uses that $\Delta'$ is good.}
\end{remark}

\section{The lower bound}\label{sec:lower}

\begin{proposition}[Two colors]\label{prop:two}
For every $k\ge 3$, $T_k(2)=k^2+k-1$. The coloring
\[
A_1=[1,k-1]\cup[k^2,k^2+k-2],\qquad A_2=[k,k^2-1]
\]
of $[1,k^2+k-2]$ is good.
\end{proposition}

\begin{proof}
This is our result $s_a(2)=a^2+3a+1$ \cite{Khamseh} with $a=k-1$; we include the short argument. 
\emph{Upper bound.} Suppose $\Delta(1)=1$. Then $1+(k-1)\cdot 1=k$ gives $\Delta(k)=2$, $k\cdot k=k^2$ gives $\Delta(k^2)=1$, and $1+(k-1)(k+1)=k^2$ gives $\Delta(k+1)=2$. Now $k^2+k-1=k^2+(k-1)\cdot 1=k+(k-1)(k+1)$ can have neither color. \emph{Lower bound.} In $A_2$ every value $x+(k-1)y$ is at least $k+(k-1)k=k^2>k^2-1$. In $A_1$, if $y\le k-1$ then $x+(k-1)y$ lies in $[k,k^2-k]\cup[k^2+k-1,\infty)$, which misses $A_1$; if $y\ge k^2$ the value exceeds $k^2+k-2$.
\end{proof}

\begin{theorem}[Lower bound]\label{thm:lower}
For every $k\ge 3$ the coloring
\begin{align*}
A_1&=[1,k-1]\cup[k^2,k^2+k-2]\cup[k^3+k^2-k,k^3+k^2-2]\cup[k^3+2k^2-k-1,k^3+2k^2-3],\\
A_2&=[k,k^2-1]\cup[k^3+k^2-1,k^3+2k^2-k-2],\\
A_3&=[k^2+k-1,k^3+k^2-k-1]
\end{align*}
of $[1,k^3+2k^2-3]$ is good. Hence $T_k(3)\ge k^3+2k^2-2$.
\end{theorem}

The construction is the case $r=3$ of Khamseh's recursive lower bound \cite[Theorem~8]{Khamseh}.\footnote{The same doubling-with-a-middle-block idea goes back to the classical lower bound $S(r)\ge 3S(r-1)+1$ for Schur numbers \cite{Schur} and to Znám's bound for \eqref{eq:kterm} \cite{Znam}. For \eqref{eq:main} it gives $T_k(r)\ge N_r$ with $N_r=(k+1)N_{r-1}-1$ and $N_1=k$, that is $N_r=((k-1)(k+1)^r+1)/k$. This agrees with $(a(a+2)^r+1)/(a+1)$ in \cite{Khamseh} for $a=k-1$.} Put $N_2=k^2+k-1$ and $s=kN_2-1=k^3+k^2-k-1$. The interval $L=[1,N_2-1]$ carries the two-coloring of Proposition~\ref{prop:two}, the middle block $M=[N_2,s]$ is entirely color $3$, and $U=L+s$ is a translate of $L$. The total length is $(k+1)N_2-2=k^3+2k^2-3$.

\begin{proof}
Suppose that $x+(k-1)y=z$ with $x,y,z\in[1,N-1]$ all of the same color.
\emph{Color $3$.} If $x,y\in M$ then $z\ge N_2+(k-1)N_2=kN_2=s+1$, so $z\notin M$.
\emph{Colors $1$ and $2$.} These colors occur only in $L\cup U$. If $y\in U$ then $y\ge s+1=kN_2$ and $z>(k-1)kN_2>(k+1)N_2>N$, which is impossible, so $y\in L$. If also $x\in L$, then $z\le k(N_2-1)=s+1-k$, so $z\in L\cup M$. If $z\in M$ it has color $3$, which does not occur in $L$. If $z\in L$, the triple lies in $L$, where the coloring is good by Proposition~\ref{prop:two}. Finally, let $x=x_0+s\in U$ with $x_0\in L$. Then $z=(x_0+(k-1)y)+s$. If $x_0+(k-1)y\le N_2-1$, then $z\in U$ and, because $U$ is a translate of $L$, the triple is monochromatic exactly when $(x_0,y,x_0+(k-1)y)$ is monochromatic in $L$, which is impossible. Otherwise $z\ge N_2+s=N$, which is out of range.
\end{proof}

We also verified by exhaustive enumeration of all solutions that this coloring is good for $3\le k\le 12$.

\begin{example}\label{ex:k3lower}
For $k=3$ the equation is $x+2y=z$, $N_2=11$, $s=32$ and $N=43$. The coloring of $[1,42]$ is
\[
\mathtt{11\,222222\,11\,3333333333333333333333\,11\,222222\,11}.
\]
The integer $43$ has no admissible color, since $43=41+2\cdot 1=35+2\cdot 4=21+2\cdot 11$ with left-hand sides in $A_1$, $A_2$ and $A_3$ respectively.
\end{example}

\section{The critical set $\Xmin$}\label{sec:set}

Table~\ref{tab:X} lists the $33$ polynomials of $\Xmin$: $9$ linear, $13$ quadratic and $11$ cubic. Each element $ak^3+bk^2+ck+d$ is also recorded by its coefficient vector $(a,b,c,d)$, together with its values for $k=3,4,5,6$.

\begin{longtable}{>{$}l<{$}ll}
\caption{The critical set $\Xmin$ ($33$ elements).}\label{tab:X}\\
\toprule
\text{element} & $(a,b,c,d)$ & values at $k=3,4,5,6$\\
\midrule\endfirsthead
\toprule
\text{element} & $(a,b,c,d)$ & values at $k=3,4,5,6$\\
\midrule\endhead
\bottomrule\endfoot
1 & $(0,0,0,1)$ & 1, 1, 1, 1 \\
2 & $(0,0,0,2)$ & 2, 2, 2, 2 \\
k - 1 & $(0,0,1,-1)$ & 2, 3, 4, 5 \\
k & $(0,0,1,0)$ & 3, 4, 5, 6 \\
k + 1 & $(0,0,1,1)$ & 4, 5, 6, 7 \\
k + 2 & $(0,0,1,2)$ & 5, 6, 7, 8 \\
2k - 1 & $(0,0,2,-1)$ & 5, 7, 9, 11 \\
2k & $(0,0,2,0)$ & 6, 8, 10, 12 \\
3k - 2 & $(0,0,3,-2)$ & 7, 10, 13, 16 \\
k^{2} - 2 & $(0,1,0,-2)$ & 7, 14, 23, 34 \\
k^{2} - 1 & $(0,1,0,-1)$ & 8, 15, 24, 35 \\
k^{2} & $(0,1,0,0)$ & 9, 16, 25, 36 \\
k^{2} + k - 2 & $(0,1,1,-2)$ & 10, 18, 28, 40 \\
k^{2} + k - 1 & $(0,1,1,-1)$ & 11, 19, 29, 41 \\
k^{2} + k & $(0,1,1,0)$ & 12, 20, 30, 42 \\
k^{2} + k + 1 & $(0,1,1,1)$ & 13, 21, 31, 43 \\
k^{2} + 2k - 2 & $(0,1,2,-2)$ & 13, 22, 33, 46 \\
k^{2} + 2k - 1 & $(0,1,2,-1)$ & 14, 23, 34, 47 \\
k^{2} + 2k & $(0,1,2,0)$ & 15, 24, 35, 48 \\
2k^{2} - 2k + 1 & $(0,2,-2,1)$ & 13, 25, 41, 61 \\
2k^{2} - 1 & $(0,2,0,-1)$ & 17, 31, 49, 71 \\
2k^{2} + k - 2 & $(0,2,1,-2)$ & 19, 34, 53, 76 \\
k^{3} - k^{2} + 1 & $(1,-1,0,1)$ & 19, 49, 101, 181 \\
k^{3} - k + 1 & $(1,0,-1,1)$ & 25, 61, 121, 211 \\
k^{3} & $(1,0,0,0)$ & 27, 64, 125, 216 \\
k^{3} + k - 1 & $(1,0,1,-1)$ & 29, 67, 129, 221 \\
k^{3} + k^{2} - 2k + 1 & $(1,1,-2,1)$ & 31, 73, 141, 241 \\
k^{3} + k^{2} - k & $(1,1,-1,0)$ & 33, 76, 145, 246 \\
k^{3} + k^{2} - 1 & $(1,1,0,-1)$ & 35, 79, 149, 251 \\
k^{3} + 2k^{2} - 3k + 1 & $(1,2,-3,1)$ & 37, 85, 161, 271 \\
k^{3} + 2k^{2} - 2k & $(1,2,-2,0)$ & 39, 88, 165, 276 \\
k^{3} + 2k^{2} - k - 1 & $(1,2,-1,-1)$ & 41, 91, 169, 281 \\
k^{3} + 2k^{2} - 2 & $(1,2,0,-2)$ & 43, 94, 173, 286 \\
\end{longtable}

\begin{theorem}\label{thm:critical}
The set $\Xmin$ has the following properties.
\begin{enumerate}
\item For every $k\ge 3$ and every $p\in\Xmin$, $1\le p(k)\le k^3+2k^2-2$.
\item Exactly $144$ triples $(x,y,z)\in(\Xmin)^3$ satisfy the identity $z=x+(k-1)y$ in $\mathbb{Z}[k]$. They are listed in Appendix~\ref{app:identities}.
\item $\Xmin$ is critical: every map $\Xmin\to\{1,2,3\}$ makes one of these $144$ identities monochromatic.
\item $\Xmin$ is minimally critical: for every $p\in\Xmin$, the set $\Xmin\setminus\{p\}$, with the identities among its elements, has a good map.
\end{enumerate}
\end{theorem}

\begin{proof}
(1) The largest element is $N$ itself, and $1,2$ and $k-1$ are the smallest. For each $p$, the polynomials $p-1$ and $N-p$ have integer coefficients, and a direct computation shows that $1+\max_{i<n}|c_i|/|c_n|\le 5$ for each of them, where $c_n$ is the leading coefficient. By Cauchy's bound on the roots of a polynomial \cite{RahmanSchmeisser}, neither polynomial changes sign for $k\ge 5$. Checking $k=3,4,5$ directly completes the proof. (2) This is a finite computation on coefficient vectors; the list is in Appendix~\ref{app:identities}. (3) This is proved by the case analysis of Section~\ref{sec:upper} and Appendix~\ref{app:cases}, and independently by the SAT computation of Section~\ref{sec:verify}. (4) Appendix~\ref{app:witness} gives, for each $p$, an explicit good map of $\Xmin\setminus\{p\}$. Each of them can be checked against the $144$ identities by hand.
\end{proof}

\begin{proof}[Proof of Theorem~\ref{thm:main}]
The lower bound is Theorem~\ref{thm:lower}. For the upper bound, apply Lemma~\ref{lem:symbolic} to $X=\Xmin$ and $M=N_3(k)$, using parts (1) and (3) of Theorem~\ref{thm:critical}.
\end{proof}

\subsection{Where the set comes from}\label{subsec:origin}
The set $\Xmin$ was not guessed. We obtained it in three steps. Let us describe them, since they explain why the elements look the way they do. Note that none of the steps is needed to \emph{check} the proof.

\emph{Step 1: a pool.} Since $N_3(k)$ is cubic, we considered all polynomials $ak^3+bk^2+ck+d$ with $a\in\{0,1\}$, $|b|\le 3$ and $|c|,|d|\le 5$ that satisfy $1\le p(k)\le N_3(k)$ for every $k\ge3$. There are $910$ such polynomials, and $20\,615$ identities $z=x+(k-1)y$ hold among them. This pool is already critical.\footnote{A second, larger pool with $|b|\le4$ and $|c|,|d|\le 7$ ($1770$ polynomials) gives the same results in every step below. The range condition was checked for $3\le k<400$; by Cauchy's bound, $k<9$ already suffices for these coefficient ranges.}

\emph{Step 2: a short case analysis.} Using the branching rule of Section~\ref{sec:upper}, an exhaustive branch-and-bound search over the pool, with memoization of repeated subproblems, finds a smallest case analysis. It has $35$ cases. Keeping only the elements on which some contradiction actually depends leaves a set $X_k$ of $51$ polynomials with $205$ identities.\footnote{Hence $35$ is the least number of cases that the rules of Section~\ref{sec:prelim} need over the whole pool. The present paper deliberately uses the smaller set $\Xmin$, at the cost of a longer case analysis ($62$ cases instead of $35$). A critical set of minimum size is the more canonical object, and it is the one that a reader or a proof checker has to trust.}

\emph{Step 3: a smallest critical subset.} A smallest critical subset of $X_k$ was computed by implicit hitting-set dualization \cite{Reiter,LiffitonSakallah,IgnatievPLM}. Every critical subset must meet the complement of every maximal good subset (a \emph{correction set}). A minimum hitting set of the correction sets found so far is therefore a lower bound, and when such a hitting set is itself critical it is optimal. The procedure stopped after $360$ iterations with the $33$-element set $\Xmin$. Thus $\Xmin$ is a critical subset of $X_k$ of minimum cardinality, not merely a minimal one. Running the same procedure on another pool of $60$ candidate polynomials also gives the minimum $33$, attained by a different set.

\begin{remark}[The spine]\label{rem:spine}
Eight elements of the $910$-element pool are \emph{necessary}: removing any one of them makes the pool $3$-colorable, so they belong to every critical subset. They are
\[
1,\ k,\ k^2,\ k^2+k-1,\ k^3+k^2-k,\ k^3+k^2-1,\ k^3+2k^2-k-1,\ k^3+2k^2-2 .
\]
They form a chain in which each element arises from the previous one by $\times k$, $+(k-1)\cdot 1$ or $+(k-1)\cdot k$. They are exactly the left end-points of the blocks of the extremal coloring (Theorem~\ref{thm:lower}), together with $N$. All eight lie in $\Xmin$. For $k=4,5$ the necessary elements of the numerical instance $[1,N_3(k)]$ are exactly these eight numbers, namely $1,4,16,19,76,79,91,94$ and $1,5,25,29,145,149,169,173$.
\end{remark}

\begin{remark}[How small can a universal set be?]\label{rem:m}
Call a critical set $X\subset\mathbb{Z}[k]$ satisfying part~(1) of Theorem~\ref{thm:critical} \emph{universal}, and let $m$ be the least size of a universal set. By the proof of Lemma~\ref{lem:symbolic}, $X(k)=\{p(k):p\in X\}$ is, for each $k$, a set of at most $|X|$ integers in $[1,N_3(k)]$ that has no good coloring. Hence $m\ge m_k$, where $m_k$ is the least size of such a set of integers. Hitting-set computations give $m_3=21$ exactly (for example $\{1,2,3,4,5,7,9,11,12,13,17,19,23,25,27,29,31,33,35,41,43\}$ and $m_4\ge 24$. Therefore
\[
24\le m\le 33 .
\]
Whether $\Xmin$ is optimal among \emph{all} universal sets is open.\footnote{The value $m_3=21$ is small because many coincidences occur at $k=3$, for example $k^2=3k$ and $k+2=2k-1$, and these create extra relations. At $k=3$ the $33$ elements of $\Xmin$ take only $27$ distinct values; for every $k\ge 4$ they take $33$ distinct values.}
\end{remark}

\section{The upper bound: a $62$-case analysis}\label{sec:upper}

\subsection{The branching rule}
The case analysis starts from $\Delta(1)=1$ and $\Delta(k)=2$ (Lemma~\ref{lem:norm}) and uses only elements of $\Xmin$ and the $144$ identities. At each node, Rule~\ref{rule2} is applied repeatedly until nothing new is forced. If Rule~\ref{rule3} then applies, or a relation has become monochromatic, the node is a leaf. Otherwise we choose an element $e$ for which exactly one color $c$ is excluded by a relation whose other members already have color $c$, and we split into the two remaining colors. In the text of Appendix~\ref{app:cases} this reads: ``Since \dots\ must not be monochromatic in $A_c$, either $\Delta(e)=i$ or $\Delta(e)=j$.'' Every split therefore has an explicit reason and exactly two branches.

\begin{figure}[ht]
\centering
\begin{tikzpicture}[
  grow'=right, level distance=26mm, edge from parent path={(\tikzparentnode.east) -- (\tikzchildnode.west)},
  every node/.style={font=\footnotesize, anchor=west, inner sep=1.5pt},
  level 1/.style={sibling distance=34mm},
  level 2/.style={sibling distance=17mm},
  level 3/.style={sibling distance=8.5mm}]
\node {$\Delta(1)=1,\ \Delta(k)=2$}
  child {node {$\Delta(k^2)=1$\ \ [23]}
    child {node {$\Delta(k^2+k-1)=2$\ \ [2]}}
    child {node {$\Delta(k^2+k-1)=3$\ \ [21]}
      child {node {$\Delta(k^3+k^2-k)=1$\ \ [10]}}
      child {node {$\Delta(k^3+k^2-k)=2$\ \ [11]}}}}
  child {node {$\Delta(k^2)=3$\ \ [39]}
    child {node {$\Delta(k^3)=1$\ \ [21]}
      child {node {$\Delta(k^3+k-1)=2$\ \ [11]}}
      child {node {$\Delta(k^3+k-1)=3$\ \ [10]}}}
    child {node {$\Delta(k^3)=2$\ \ [18]}
      child {node {$\Delta(k^3+k^2-k)=1$\ \ [9]}}
      child {node {$\Delta(k^3+k^2-k)=3$\ \ [9]}}}};
\end{tikzpicture}
\caption{The first three levels of the case analysis on $\Xmin$. The number of leaves below each node is shown in brackets; there are $62$ in total.}\label{fig:tree}
\end{figure}
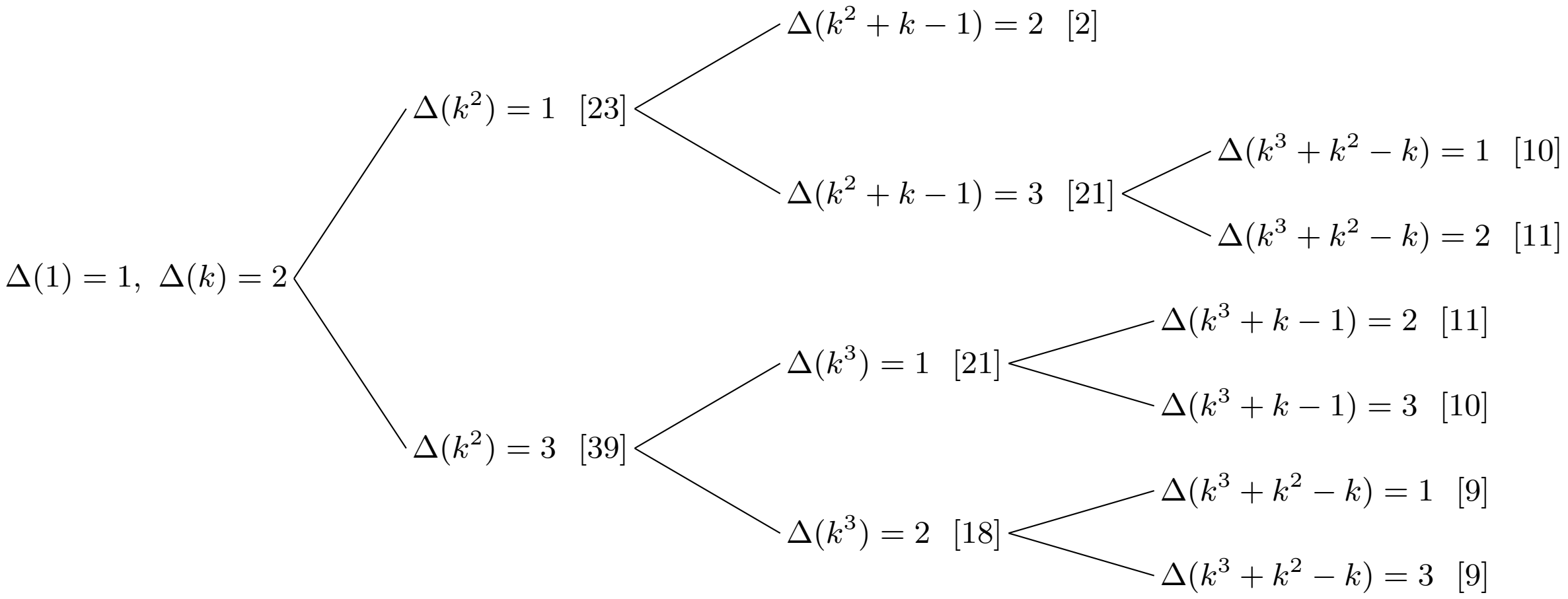

The first split is on $k^2$: since $k\cdot k=k^2$ and $\Delta(k)=2$, $\Delta(k^2)\in\{1,3\}$. The resulting tree has $62$ leaves, $61$ splits, depth $8$ and $374$ applications of Rule~\ref{rule2}, and it uses $142$ of the $144$ identities and all $33$ elements. Figure~\ref{fig:tree} shows the first three levels. The complete scheme, with the reason for every leaf, is printed at the start of Appendix~\ref{app:cases}.

\subsection{The case analysis read as four propositions}
The tree splits naturally into four statements. In each, $\Delta$ is a good map of $\Xmin$ with $\Delta(1)=1$ and $\Delta(k)=2$. The case numbers refer to Appendix~\ref{app:cases}.

\begin{proposition}\label{prop:k2}
$\Delta(k^2)=1$. \textup{(Cases 2.*, $39$ leaves.)}
\end{proposition}

\begin{proof}[Idea of proof]
Color $2$ is excluded by $k\cdot k=k^2$. Color $3$ is refuted as follows. Since $k\cdot k^2=k^3$, the element $k^3$ has color $1$ or $2$. In the first case we split on $k^3+k-1$ (excluded from color $1$ by $k^3+(k-1)\cdot1$); in the second, on $k^3+k^2-k$ (excluded from color $2$ by $k^3+(k-1)\cdot k$). Refuting this ``wrong'' color of $k^2$ is the largest part of the proof, because many nearly good maps have to be excluded.
\end{proof}

\begin{proposition}\label{prop:k2k1}
If $\Delta(k^2)=1$ then $\Delta(k^2+k-1)=3$. \textup{(Cases 1.1.*, $2$ leaves.)}
\end{proposition}

\begin{proof}[Idea of proof]
Color $1$ is excluded by $k^2+(k-1)\cdot 1=k^2+k-1$. If $\Delta(k^2+k-1)=2$, Rule~\ref{rule2} applied to $1+(k-1)(k+1)=k^2$ and $k+(k-1)(k+1)=k^2+k-1$ forces $\Delta(k+1)=3$. A split on $k^2+k$ then closes both branches, at $k^3+k^2-k$ and at $k^2-1$ respectively.
\end{proof}

\begin{proposition}\label{prop:B8a}
If $\Delta(k^2)=1$ and $\Delta(k^2+k-1)=3$, then $\Delta(k^3+k^2-k)\ne 1$. \textup{(Cases 1.2.1.*, $10$ leaves.)}
\end{proposition}

\begin{proposition}\label{prop:B8b}
If $\Delta(k^2)=1$ and $\Delta(k^2+k-1)=3$, then $\Delta(k^3+k^2-k)\ne 2$. \textup{(Cases 1.2.2.*, $11$ leaves.)}
\end{proposition}

\begin{theorem}\label{thm:upper}
$\Xmin$ is critical. Consequently $T_k(3)\le k^3+2k^2-2$ for every $k\ge 3$.
\end{theorem}

\begin{proof}
Let $\Delta$ be a good map of $\Xmin$, normalised by Lemma~\ref{lem:norm}. By Propositions~\ref{prop:k2} and~\ref{prop:k2k1}, $\Delta(k^2)=1$ and $\Delta(k^2+k-1)=3$. Then $k\cdot(k^2+k-1)=k^3+k^2-k$ excludes color $3$ for $k^3+k^2-k$, and Propositions~\ref{prop:B8a} and~\ref{prop:B8b} exclude colors $1$ and $2$. This is a contradiction. The second statement follows from Lemma~\ref{lem:symbolic} and Theorem~\ref{thm:critical}(1).
\end{proof}

\begin{remark}\label{rem:roleN}
In every good coloring of $[1,N-1]$ the element $k^3+k^2-k$ has color $1$ (Section~\ref{sec:skeleton}). Proposition~\ref{prop:B8a} can refute this color only because $\Xmin$ contains $N=k^3+2k^2-2$ and its neighbours $k^3+2k^2-k-1$ and $k^3+2k^2-2k$. This is where the exact value of the upper bound enters the proof: without $N$ the set would be $3$-colorable (Remark~\ref{rem:spine}).
\end{remark}

\section{Machine verification}\label{sec:verify}

All statements of Theorem~\ref{thm:critical} and the whole case analysis were checked by computer. The programs are independent of the search that produced $\Xmin$.\footnote{All computations were run in Python~3 with PySAT \cite{PySAT}, using the CaDiCaL solver \cite{CaDiCaL} (version 1.5.3 as bundled in PySAT). Each check takes at most a few seconds on a laptop, except the parametric runs of Section~\ref{sec:skeleton}, which take about two and a half minutes in total.}

\subsection*{SAT encoding}
For each $p\in\Xmin$ we use three Boolean variables $v_{p,1},v_{p,2},v_{p,3}$ with the clause $v_{p,1}\vee v_{p,2}\vee v_{p,3}$. For each identity $(x,y,z)$ and each color $c$ we add the clause $\neg v_{x,c}\vee\neg v_{y,c}\vee\neg v_{z,c}$, which says that the triple is not entirely of color $c$. (Clauses stating ``at most one color'' are not needed for unsatisfiability.) The resulting formula has $99$ variables and $33+3\cdot144=465$ clauses, and it is unsatisfiable. This proves Theorem~\ref{thm:critical}(3) a second time.

\subsection*{Checks performed}
\begin{enumerate}
\item \emph{Identities.} All triples $(x,y,z)\in(\Xmin)^3$ with $z=x+(k-1)y$ in $\mathbb{Z}[k]$ were enumerated; there are exactly $144$.
\item \emph{Range.} $1\le p(k)\le N_3(k)$ for all $p\in\Xmin$ and $3\le k<400$, and the Cauchy bound argument covers all larger $k$.
\item \emph{Criticality.} The formula above is unsatisfiable.
\item \emph{Minimality.} For each of the $33$ elements $p$, the formula for $\Xmin\setminus\{p\}$ is satisfiable, and the satisfying assignment was checked directly against the identities (Appendix~\ref{app:witness}).
\item \emph{Numerical instances.} For every $k=3,\dots,20$ the values $\Xmin(k)$ were treated as integers, \emph{all} numerical solutions of \eqref{eq:main} among them were used (so coincidences are included), and the resulting formula is unsatisfiable. The number of distinct values is $27$ for $k=3$ and $33$ for $4\le k\le 20$.
\item \emph{The case analysis.} The text of Appendix~\ref{app:cases} was parsed back into a tree. For every node the checker verifies the following: each relation quoted is a polynomial identity between elements of $\Xmin$; in each application of Rule~\ref{rule2} the other members of the two relations have two different known colors; each split excludes exactly the stated color and has exactly the two remaining colors as branches; and each leaf is closed by Rule~\ref{rule3} (the other members of the three relations cover all three colors) or by a monochromatic relation. The checker confirms $62$ leaves, $61$ splits, $374$ applications of Rule~\ref{rule2} and depth $8$.
\item \emph{Lower bound.} The coloring of Theorem~\ref{thm:lower} was checked against all solutions in $[1,N_3(k)-1]$ for $3\le k\le 12$.
\end{enumerate}

\section{The structure of extremal colorings}\label{sec:skeleton}

The coloring in Theorem~\ref{thm:lower} is not the only good coloring of $[1,N-1]$. However, it is almost forced. This explains why we cannot avoid a case analysis: unlike the two-color case, some integers are genuinely free, so no chain of forced colors alone can reach $N$.\footnote{For two colors the extremal coloring of $[1,k^2+k-2]$ is unique up to renaming the colors, and four forced steps suffice (Proposition~\ref{prop:two}). With three colors there are $16$ good colorings of $[1,42]$ for $k=3$ with $\Delta(1)=1$ and $\Delta(3)=2$, and many more for larger $k$.}

\begin{table}[ht]\centering\small
\caption{The ten forced blocks of good colorings of $[1,N-1]$ with $\Delta(1)=1$ and $\Delta(k)=2$.}\label{tab:blocks}
\begin{tabular}{llcccc}
\toprule block & interval & color & $k=3$ & $k=4$ & $k=5$\\ \midrule
$B_1$ & $[1,k-1]$ & 1 & $[1,2]$ & $[1,3]$ & $[1,4]$\\
$B_2$ & $[k,k^2-1]$ & 2 & $[3,8]$ & $[4,15]$ & $[5,24]$\\
$B_3$ & $[k^2,k^2+k-2]$ & 1 & $[9,10]$ & $[16,18]$ & $[25,28]$\\
$B_4$ & $[k^2+k-1,2k^2-2]$ & 3 & $[11,16]$ & $[19,30]$ & $[29,48]$\\
$B_5$ & $[2k^2+2k-3,3k^2-3]$ & 3 & $[21,24]$ & $[37,45]$ & $[57,72]$\\
$B_6$ & $[k^3-k^2+1,k^3-2k+1]$ & 3 & $[19,22]$ & $[49,57]$ & $[101,116]$\\
$B_7$ & $[k^3,k^3+k^2-k-1]$ & 3 & $[27,32]$ & $[64,75]$ & $[125,144]$\\
$B_8$ & $[k^3+k^2-k,k^3+k^2-2]$ & 1 & $[33,34]$ & $[76,78]$ & $[145,148]$\\
$B_9$ & $[k^3+k^2-1,k^3+2k^2-k-2]$ & 2 & $[35,40]$ & $[79,90]$ & $[149,168]$\\
$B_{10}$ & $[k^3+2k^2-k-1,k^3+2k^2-3]$ & 1 & $[41,42]$ & $[91,93]$ & $[169,172]$\\
\bottomrule
\end{tabular}
\end{table}

Let $F_1=[2k^2-1,2k^2+2k-4]$, $F_2=[3k^2-2,k^3-k^2]$ and $F_3=[k^3-2k+2,k^3-1]$. For $k\ge4$ the blocks $B_1,\dots,B_{10}$ and the regions $F_1,F_2,F_3$ partition $[1,N-1]$.

\begin{theorem}[Parametric forcing]\label{thm:param}
Let $\Delta$ be a good coloring of $[1,N-1]$ with $\Delta(1)=1$ and $\Delta(k)=2$. Then $B_1\cup B_3\cup B_8\cup B_{10}\subseteq A_1$, color $2$ does not occur in $B_4\cup B_5\cup B_6\cup B_7$, and color $3$ does not occur in $B_2\cup B_9$.
\end{theorem}

\begin{proof}
A statement such as ``no $t\in B_4$ has color $2$'' concerns a generic element $t$ of an interval whose end-points are polynomials in $k$. We apply the symbolic principle with two parameters. The elements are expressions $P(k)+Q(k)\,t$, where $P$ has small coefficients and
\[
Q\in\{0,\,1,\,2,\,k-1,\,k,\,k+1,\,2k-2,\,k^2-k,\,k^2-1,\,k^2,\,k^2-2k+1,\,k^2+k,\,k^2+k-1\}.
\]
The relations are identities in both $k$ and $t$. The value of each element is linear in $t$, so checking that it lies in $[1,N-1]$ at the two end-points of the block suffices for every $t$ in the block. The Cauchy bound, as in Theorem~\ref{thm:critical}(1), reduces the range condition to finitely many values of $k$. For each of the fourteen pairs (block, excluded color)
\begin{gather*}
(B_1,2),\ (B_1,3),\ (B_2,3),\ (B_3,2),\ (B_3,3),\ (B_4,2),\ (B_5,2),\\
(B_6,2),\ (B_7,2),\ (B_8,2),\ (B_8,3),\ (B_9,3),\ (B_{10},2),\ (B_{10},3),
\end{gather*}
the SAT solver shows that the assumptions $\Delta(1)=1$, $\Delta(k)=2$ and $\Delta(t)=c$ admit no good map. All fourteen runs are unsatisfiable.
\end{proof}

The statements ``color $1$ does not occur in $B_2,B_4,\dots,B_9$'' cannot be proved from any fixed family $P+Qt$, because each needs about $k$ elements that depend on $t$. They follow instead from a counting argument on intervals.

\begin{lemma}\label{lem:interval}
If $X$ is an interval of $k-1$ consecutive integers and $Y=[c,d]$, then $X+(k-1)Y=[\min X+(k-1)c,\ \max X+(k-1)d]$. The same holds for $X-(k-1)Y$.
\end{lemma}

\begin{proof}
$X$ meets every residue class modulo $k-1$, so every integer between the two bounds has the form $x+(k-1)y$ with $x\in X$ and $y\in Y$.
\end{proof}

\begin{theorem}[Rigid skeleton]\label{thm:skeleton}
Under the hypotheses of Theorem~\ref{thm:param}, $A_2=B_2\cup B_9$ exactly, $A_1\supseteq B_1\cup B_3\cup B_8\cup B_{10}$, $A_3\supseteq B_4\cup B_5\cup B_6\cup B_7$, and $F_1\cup F_2\cup F_3\subseteq A_1\cup A_3$.
\end{theorem}

\begin{proof}
Write $K_1=B_1$, $K_3=B_3$, $K_8=B_8$ and $K_{10}=B_{10}$. By Theorem~\ref{thm:param} these lie in $A_1$, and each has exactly $k-1$ elements. An integer $t$ cannot have color $1$ if $t\in X+(k-1)Y$ or $t\in X-(k-1)Y$ with $X,Y\subseteq A_1$. By Lemma~\ref{lem:interval}:
\begin{align*}
B_2&\subseteq (K_1+(k-1)K_1)\cup(K_3-(k-1)\cdot1)=[k,k^2-k]\cup[k^2-k+1,k^2-1],\\
B_4&\subseteq (K_3+(k-1)K_1)\cup(K_8-(k-1)k^2)\\
&=[k^2+k-1,2k^2-k-1]\cup[2k^2-k,2k^2-2],\\
B_5&=K_{10}-(k-1)K_3,\qquad B_6=K_1+(k-1)K_3,\\
B_7&\subseteq (K_3+(k-1)K_3)\cup(K_8-(k-1)K_1)\\
&=[k^3,k^3+k^2-2k]\cup[k^3+k-1,k^3+k^2-k-1],\\
B_9&\subseteq (K_8+(k-1)K_1)\cup(K_{10}-(k-1)K_1)\\
&=[k^3+k^2-1,k^3+2k^2-2k-1]\cup[k^3+k^2+k-2,k^3+2k^2-k-2].
\end{align*}
Hence color $1$ occurs in none of $B_2,B_4,\dots,B_7,B_9$, and together with Theorem~\ref{thm:param} this determines their colors. Finally $B_2+(k-1)B_2=[k^2,k^3-k]$ and $B_9-(k-1)B_2=[2k^2+k-2,k^3+k^2-2]$ cover $[k^2,k^3+k^2-2]$, so color $2$ occurs nowhere in this interval outside $B_9$. In particular it does not occur in $F_1\cup F_2\cup F_3$, and $A_2=B_2\cup B_9$. Each interval identity above follows by computing the end-points; we also checked them by brute force for $3\le k\le 60$.
\end{proof}

\begin{corollary}[Second proof of the upper bound]\label{cor:second}
$[1,N]$ has no good coloring.
\end{corollary}

\begin{proof}
A good coloring of $[1,N]$ restricts to one of $[1,N-1]$, so the rigid skeleton applies. But
\[
N=(k^3+2k^2-k-1)+(k-1)\cdot 1=(k^3+k^2-1)+(k-1)\cdot k=(2k^2+2k-3)+(k-1)(k^2+k-1),
\]
with the other members in $A_1$, $A_2$ and $A_3$ respectively, so Rule~\ref{rule3} applies at $N$.
\end{proof}

We also computed, for $4\le k\le 7$, the set of colors possible for each integer in a good coloring of $[1,N-1]$ with $\Delta(1)=1$ and $\Delta(k)=2$. The free integers, those admitting two colors, are exactly $F_1\cup F_2\cup F_3$, and each admits colors $1$ and $3$. For $k=3$ the blocks overlap, the free integers are $17,18,25,26$, and there are $2^4=16$ good colorings.\footnote{This second proof relies on fourteen solver runs with two parameters. The proof via $\Xmin$ relies on a single set of $33$ polynomials and a case analysis that can be read line by line. This is why we regard the latter as the main proof.}

\section{A worked example: $k=3$}\label{sec:k3}

For $k=3$ the equation is $x+2y=z$ and $N=43$. The $33$ polynomials of $\Xmin$ take the $27$ distinct values
\[
1,2,3,4,5,6,7,8,9,10,11,12,13,14,15,17,19,25,27,29,31,33,35,37,39,41,43.
\]
Coincidences such as $k+2=2k-1=5$ and $k^2-2=3k-2=7$ are harmless (Remark~\ref{rem:coinc}).

We follow Case~1.1.1 of Appendix~\ref{app:cases} numerically. The assumptions are $\Delta(1)=1$, $\Delta(3)=2$, $\Delta(9)=1$ (that is, $\Delta(k^2)=1$), $\Delta(11)=2$ (that is, $\Delta(k^2+k-1)=2$) and $\Delta(12)=1$ (that is, $\Delta(k^2+k)=1$).
\begin{itemize}
\item Rule~\ref{rule2} on $1+2\cdot4=9$ (color $1$) and $3+2\cdot4=11$ (color $2$) gives $\Delta(4)=3$. This is $\Delta(k+1)=3$, forced already in Case~1.1.
\item Rule~\ref{rule2} on $1+2\cdot12=25$ (color $1$) and $3+2\cdot11=25$ (color $2$) gives $\Delta(25)=3$. This is $\Delta(k^3-k+1)=3$.
\item Contradiction at $33=k^3+k^2-k$: $9+2\cdot12=33$ has its other members in $A_1$, $3\cdot 11=33$ in $A_2$, and $25+2\cdot4=33$ in $A_3$. So $33$ can have no color.
\end{itemize}
The symbolic version of these steps is identical, with $k^2+k$, $k^3-k+1$ and $k^3+k^2-k$ in place of $12$, $25$ and $33$. This is why one argument serves every $k$.

\section{Concluding remarks}\label{sec:concl}

The formula $T_k(3)=k^3+2k^2-2$ is known from the equivalent Rado-number result of Chang, De Loera, and Wesley \cite{ChangDeLoeraWesley}. The contribution of this paper is an independent proof based on a fixed critical set of $33$ polynomials, together with a complete $62$-case argument using only two elementary rules. The set is critical, and it has minimum cardinality among the critical subsets of the $51$-element set from which it was extracted. We also describe the rigid structure shared by all extremal colorings of $[1,N_3(k)-1]$. These results provide a smaller and more explicit certificate for the known formula, as well as structural information about the extremal colorings. We end with some questions that remain open.
\begin{enumerate}
\item Determine $m$, the least size of a universal critical set. We know $24\le m\le 33$ (Remark~\ref{rem:m}). Computing $m_k$ exactly for $k=4,5,6$ would narrow the gap, and if some optimal numerical set is ``generic'' it may lift to a smaller universal set.
\item Find a case analysis on $\Xmin$ with fewer than $62$ cases, or a different framework that proves Proposition~\ref{prop:k2} (the $39$ cases refuting $\Delta(k^2)=3$) more conceptually. The interval argument of Section~\ref{sec:skeleton} is a candidate.
\item The construction of Theorem~\ref{thm:lower} gives $T_k(r)\ge ((k-1)(k+1)^r+1)/k$ for every $r$ \cite{Khamseh}. It would be interesting to know whether equality holds for $r\ge 4$, but this is beyond the scope of this paper.
\end{enumerate}

\subsection*{Competing Interests}
The authors report there are no competing interests to declare.



\newpage
\appendix

\section{The complete $62$-case analysis}\label{app:cases}

\subsection*{Scheme}
Each line shows a split, with the case number and the number of leaves below it. Each leaf shows the element at which Rule~\ref{rule3} applies, or the relation that becomes monochromatic. Throughout, $\Delta(1)=1$ and $\Delta(k)=2$ (Lemma~\ref{lem:norm}).

\medskip
{\small
\schemeitem{0}{$\Delta(k^{2}) = 1$}{1}{23}
\schemeitem{1}{$\Delta(k^{2} + k - 1) = 2$}{1.1}{2}
\schemeitem{2}{$\Delta(k^{2} + k) = 1$}{1.1.1}{1}
\schemeleaf{3}{contradiction at $k^{3} + k^{2} - k$}
\schemeitem{2}{$\Delta(k^{2} + k) = 2$}{1.1.2}{1}
\schemeleaf{3}{contradiction at $k^{2} - 1$}
\schemeitem{1}{$\Delta(k^{2} + k - 1) = 3$}{1.2}{21}
\schemeitem{2}{$\Delta(k^{3} + k^{2} - k) = 1$}{1.2.1}{10}
\schemeitem{3}{$\Delta(k + 1) = 2$}{1.2.1.1}{6}
\schemeitem{4}{$\Delta(k^{3} + k^{2} - 2k + 1) = 2$}{1.2.1.1.1}{2}
\schemeitem{5}{$\Delta(k^{2} + 2k - 1) = 1$}{1.2.1.1.1.1}{1}
\schemeleaf{6}{contradiction at $k^{3} + 2k^{2} - k - 1$}
\schemeitem{5}{$\Delta(k^{2} + 2k - 1) = 3$}{1.2.1.1.1.2}{1}
\schemeleaf{6}{contradiction at $k^{2} + 2k$}
\schemeitem{4}{$\Delta(k^{3} + k^{2} - 2k + 1) = 3$}{1.2.1.1.2}{4}
\schemeitem{5}{$\Delta(k^{2} - 1) = 1$}{1.2.1.1.2.1}{1}
\schemeleaf{6}{contradiction at $3k - 2$}
\schemeitem{5}{$\Delta(k^{2} - 1) = 2$}{1.2.1.1.2.2}{3}
\schemeitem{6}{$\Delta(k^{3} + 2k^{2} - k - 1) = 1$}{1.2.1.1.2.2.1}{2}
\schemeitem{7}{$\Delta(2k) = 1$}{1.2.1.1.2.2.1.1}{1}
\schemeleaf{8}{contradiction at $2$}
\schemeitem{7}{$\Delta(2k) = 2$}{1.2.1.1.2.2.1.2}{1}
\schemeleaf{8}{contradiction at $2$}
\schemeitem{6}{$\Delta(k^{3} + 2k^{2} - k - 1) = 3$}{1.2.1.1.2.2.2}{1}
\schemeleaf{7}{contradiction at $3k - 2$}
\schemeitem{3}{$\Delta(k + 1) = 3$}{1.2.1.2}{4}
\schemeitem{4}{$\Delta(k^{3} + k^{2} - 2k + 1) = 2$}{1.2.1.2.1}{2}
\schemeitem{5}{$\Delta(2) = 1$}{1.2.1.2.1.1}{1}
\schemeleaf{6}{contradiction at $k^{2} + 2k - 1$}
\schemeitem{5}{$\Delta(2) = 3$}{1.2.1.2.1.2}{1}
\schemeleaf{6}{contradiction at $k^{3} + 2k^{2} - 2k$}
\schemeitem{4}{$\Delta(k^{3} + k^{2} - 2k + 1) = 3$}{1.2.1.2.2}{2}
\schemeitem{5}{$\Delta(k - 1) = 1$}{1.2.1.2.2.1}{1}
\schemeleaf{6}{contradiction at $k^{2} + 2k$}
\schemeitem{5}{$\Delta(k - 1) = 2$}{1.2.1.2.2.2}{1}
\schemeleaf{6}{contradiction at $2$}
\schemeitem{2}{$\Delta(k^{3} + k^{2} - k) = 2$}{1.2.2}{11}
\schemeitem{3}{$\Delta(k^{3} + 2k^{2} - 2k) = 1$}{1.2.2.1}{2}
\schemeitem{4}{$\Delta(k^{2} + 2k - 2) = 1$}{1.2.2.1.1}{1}
\schemeleaf{5}{$k\cdot (k^{2} + 2k - 2) = k^{3} + 2k^{2} - 2k$ monochromatic}
\schemeitem{4}{$\Delta(k^{2} + 2k - 2) = 3$}{1.2.2.1.2}{1}
\schemeleaf{5}{contradiction at $k^{2} + k - 2$}
\schemeitem{3}{$\Delta(k^{3} + 2k^{2} - 2k) = 3$}{1.2.2.2}{9}
\schemeitem{4}{$\Delta(2k - 1) = 1$}{1.2.2.2.1}{4}
\schemeitem{5}{$\Delta(3k - 2) = 2$}{1.2.2.2.1.1}{3}
\schemeitem{6}{$\Delta(k^{2} + 2k - 2) = 1$}{1.2.2.2.1.1.1}{2}
\schemeitem{7}{$\Delta(k - 1) = 2$}{1.2.2.2.1.1.1.1}{1}
\schemeleaf{8}{contradiction at $k^{3} + k^{2} - 2k + 1$}
\schemeitem{7}{$\Delta(k - 1) = 3$}{1.2.2.2.1.1.1.2}{1}
\schemeleaf{8}{contradiction at $k^{2} + k$}
\schemeitem{6}{$\Delta(k^{2} + 2k - 2) = 3$}{1.2.2.2.1.1.2}{1}
\schemeleaf{7}{$k\cdot (k^{2} + 2k - 2) = k^{3} + 2k^{2} - 2k$ monochromatic}
\schemeitem{5}{$\Delta(3k - 2) = 3$}{1.2.2.2.1.2}{1}
\schemeleaf{6}{$k\cdot (k^{2} + 2k - 2) = k^{3} + 2k^{2} - 2k$ monochromatic}
\schemeitem{4}{$\Delta(2k - 1) = 2$}{1.2.2.2.2}{5}
\schemeitem{5}{$\Delta(k + 1) = 2$}{1.2.2.2.2.1}{3}
\schemeitem{6}{$\Delta(k^{2} + 2k - 2) = 1$}{1.2.2.2.2.1.1}{2}
\schemeitem{7}{$\Delta(2) = 2$}{1.2.2.2.2.1.1.1}{1}
\schemeleaf{8}{$k\cdot (k + 2) = k^{2} + 2k$ monochromatic}
\schemeitem{7}{$\Delta(2) = 3$}{1.2.2.2.2.1.1.2}{1}
\schemeleaf{8}{$k\cdot (k + 2) = k^{2} + 2k$ monochromatic}
\schemeitem{6}{$\Delta(k^{2} + 2k - 2) = 3$}{1.2.2.2.2.1.2}{1}
\schemeleaf{7}{$k\cdot (k^{2} + 2k - 2) = k^{3} + 2k^{2} - 2k$ monochromatic}
\schemeitem{5}{$\Delta(k + 1) = 3$}{1.2.2.2.2.2}{2}
\schemeitem{6}{$\Delta(k^{3} + k^{2} - 2k + 1) = 1$}{1.2.2.2.2.2.1}{1}
\schemeleaf{7}{contradiction at $k^{2} + 2k - 2$}
\schemeitem{6}{$\Delta(k^{3} + k^{2} - 2k + 1) = 2$}{1.2.2.2.2.2.2}{1}
\schemeleaf{7}{contradiction at $3k - 2$}
\schemeitem{0}{$\Delta(k^{2}) = 3$}{2}{39}
\schemeitem{1}{$\Delta(k^{3}) = 1$}{2.1}{21}
\schemeitem{2}{$\Delta(k^{3} + k - 1) = 2$}{2.1.1}{11}
\schemeitem{3}{$\Delta(k^{3} - k + 1) = 2$}{2.1.1.1}{3}
\schemeitem{4}{$\Delta(k + 1) = 2$}{2.1.1.1.1}{1}
\schemeleaf{5}{$k\cdot (k^{2} + k - 1) = k^{3} + k^{2} - k$ monochromatic}
\schemeitem{4}{$\Delta(k + 1) = 3$}{2.1.1.1.2}{2}
\schemeitem{5}{$\Delta(k^{2} + k - 1) = 1$}{2.1.1.1.2.1}{1}
\schemeleaf{6}{contradiction at $k^{2} + k$}
\schemeitem{5}{$\Delta(k^{2} + k - 1) = 3$}{2.1.1.1.2.2}{1}
\schemeleaf{6}{contradiction at $k^{3} + 2k^{2} - k - 1$}
\schemeitem{3}{$\Delta(k^{3} - k + 1) = 3$}{2.1.1.2}{8}
\schemeitem{4}{$\Delta(k^{3} + k^{2} - k) = 2$}{2.1.1.2.1}{7}
\schemeitem{5}{$\Delta(k^{2} - 1) = 1$}{2.1.1.2.1.1}{3}
\schemeitem{6}{$\Delta(k^{2} + k - 1) = 2$}{2.1.1.2.1.1.1}{1}
\schemeleaf{7}{$k\cdot (k^{2} + k - 1) = k^{3} + k^{2} - k$ monochromatic}
\schemeitem{6}{$\Delta(k^{2} + k - 1) = 3$}{2.1.1.2.1.1.2}{2}
\schemeitem{7}{$\Delta(2k^{2} + k - 2) = 2$}{2.1.1.2.1.1.2.1}{1}
\schemeleaf{8}{contradiction at $k^{3} + 2k^{2} - k - 1$}
\schemeitem{7}{$\Delta(2k^{2} + k - 2) = 3$}{2.1.1.2.1.1.2.2}{1}
\schemeleaf{8}{$k\cdot (k + 2) = k^{2} + 2k$ monochromatic}
\schemeitem{5}{$\Delta(k^{2} - 1) = 2$}{2.1.1.2.1.2}{4}
\schemeitem{6}{$\Delta(k^{2} + 2k - 2) = 1$}{2.1.1.2.1.2.1}{2}
\schemeitem{7}{$\Delta(k - 1) = 1$}{2.1.1.2.1.2.1.1}{1}
\schemeleaf{8}{$k\cdot (k^{2} + k - 1) = k^{3} + k^{2} - k$ monochromatic}
\schemeitem{7}{$\Delta(k - 1) = 3$}{2.1.1.2.1.2.1.2}{1}
\schemeleaf{8}{contradiction at $k^{3} + 2k^{2} - k - 1$}
\schemeitem{6}{$\Delta(k^{2} + 2k - 2) = 3$}{2.1.1.2.1.2.2}{2}
\schemeitem{7}{$\Delta(2k - 1) = 1$}{2.1.1.2.1.2.2.1}{1}
\schemeleaf{8}{$k\cdot (k^{2} + k - 1) = k^{3} + k^{2} - k$ monochromatic}
\schemeitem{7}{$\Delta(2k - 1) = 2$}{2.1.1.2.1.2.2.2}{1}
\schemeleaf{8}{contradiction at $3k - 2$}
\schemeitem{4}{$\Delta(k^{3} + k^{2} - k) = 3$}{2.1.1.2.2}{1}
\schemeleaf{5}{$k\cdot (k^{2} + k - 1) = k^{3} + k^{2} - k$ monochromatic}
\schemeitem{2}{$\Delta(k^{3} + k - 1) = 3$}{2.1.2}{10}
\schemeitem{3}{$\Delta(k^{2} + k - 1) = 1$}{2.1.2.1}{6}
\schemeitem{4}{$\Delta(k^{2} + 2k - 2) = 2$}{2.1.2.1.1}{3}
\schemeitem{5}{$\Delta(k^{3} + 2k^{2} - 2k) = 1$}{2.1.2.1.1.1}{2}
\schemeitem{6}{$\Delta(k^{2} + 2k - 1) = 2$}{2.1.2.1.1.1.1}{1}
\schemeleaf{7}{contradiction at $k^{3} + k^{2} - k$}
\schemeitem{6}{$\Delta(k^{2} + 2k - 1) = 3$}{2.1.2.1.1.1.2}{1}
\schemeleaf{7}{contradiction at $k^{3} - k + 1$}
\schemeitem{5}{$\Delta(k^{3} + 2k^{2} - 2k) = 3$}{2.1.2.1.1.2}{1}
\schemeleaf{6}{contradiction at $k^{2} + 2k - 1$}
\schemeitem{4}{$\Delta(k^{2} + 2k - 2) = 3$}{2.1.2.1.2}{3}
\schemeitem{5}{$\Delta(k^{3} + k^{2} - k) = 2$}{2.1.2.1.2.1}{1}
\schemeleaf{6}{$k\cdot 2 = 2k$ monochromatic}
\schemeitem{5}{$\Delta(k^{3} + k^{2} - k) = 3$}{2.1.2.1.2.2}{2}
\schemeitem{6}{$\Delta(k^{3} + k^{2} - 2k + 1) = 2$}{2.1.2.1.2.2.1}{1}
\schemeleaf{7}{contradiction at $k^{3} + 2k^{2} - 2k$}
\schemeitem{6}{$\Delta(k^{3} + k^{2} - 2k + 1) = 3$}{2.1.2.1.2.2.2}{1}
\schemeleaf{7}{contradiction at $k^{2} + k + 1$}
\schemeitem{3}{$\Delta(k^{2} + k - 1) = 2$}{2.1.2.2}{4}
\schemeitem{4}{$\Delta(k + 1) = 1$}{2.1.2.2.1}{2}
\schemeitem{5}{$\Delta(2k - 1) = 1$}{2.1.2.2.1.1}{1}
\schemeleaf{6}{contradiction at $k^{3} + k^{2} - k$}
\schemeitem{5}{$\Delta(2k - 1) = 3$}{2.1.2.2.1.2}{1}
\schemeleaf{6}{contradiction at $k^{2} + 2k - 1$}
\schemeitem{4}{$\Delta(k + 1) = 3$}{2.1.2.2.2}{2}
\schemeitem{5}{$\Delta(2k - 1) = 1$}{2.1.2.2.2.1}{1}
\schemeleaf{6}{contradiction at $k^{3} + 2k^{2} - 2k$}
\schemeitem{5}{$\Delta(2k - 1) = 3$}{2.1.2.2.2.2}{1}
\schemeleaf{6}{$k\cdot 2 = 2k$ monochromatic}
\schemeitem{1}{$\Delta(k^{3}) = 2$}{2.2}{18}
\schemeitem{2}{$\Delta(k^{3} + k^{2} - k) = 1$}{2.2.1}{9}
\schemeitem{3}{$\Delta(k^{3} + k^{2} - 2k + 1) = 2$}{2.2.1.1}{6}
\schemeitem{4}{$\Delta(k^{3} + k^{2} - 1) = 2$}{2.2.1.1.1}{3}
\schemeitem{5}{$\Delta(k + 1) = 1$}{2.2.1.1.1.1}{2}
\schemeitem{6}{$\Delta(k + 2) = 1$}{2.2.1.1.1.1.1}{1}
\schemeleaf{7}{contradiction at $k^{2} - 1$}
\schemeitem{6}{$\Delta(k + 2) = 2$}{2.2.1.1.1.1.2}{1}
\schemeleaf{7}{contradiction at $2k^{2} + k - 2$}
\schemeitem{5}{$\Delta(k + 1) = 3$}{2.2.1.1.1.2}{1}
\schemeleaf{6}{contradiction at $k^{3} + k - 1$}
\schemeitem{4}{$\Delta(k^{3} + k^{2} - 1) = 3$}{2.2.1.1.2}{3}
\schemeitem{5}{$\Delta(k^{2} + k) = 1$}{2.2.1.1.2.1}{1}
\schemeleaf{6}{contradiction at $k^{2} + k - 2$}
\schemeitem{5}{$\Delta(k^{2} + k) = 3$}{2.2.1.1.2.2}{2}
\schemeitem{6}{$\Delta(k^{2} + k - 1) = 1$}{2.2.1.1.2.2.1}{1}
\schemeleaf{7}{$k\cdot (k^{2} + k - 1) = k^{3} + k^{2} - k$ monochromatic}
\schemeitem{6}{$\Delta(k^{2} + k - 1) = 2$}{2.2.1.1.2.2.2}{1}
\schemeleaf{7}{$k\cdot (k + 1) = k^{2} + k$ monochromatic}
\schemeitem{3}{$\Delta(k^{3} + k^{2} - 2k + 1) = 3$}{2.2.1.2}{3}
\schemeitem{4}{$\Delta(k^{2} + k - 1) = 1$}{2.2.1.2.1}{1}
\schemeleaf{5}{$k\cdot (k^{2} + k - 1) = k^{3} + k^{2} - k$ monochromatic}
\schemeitem{4}{$\Delta(k^{2} + k - 1) = 2$}{2.2.1.2.2}{2}
\schemeitem{5}{$\Delta(k + 1) = 1$}{2.2.1.2.2.1}{1}
\schemeleaf{6}{contradiction at $k^{3} + k^{2} - 1$}
\schemeitem{5}{$\Delta(k + 1) = 3$}{2.2.1.2.2.2}{1}
\schemeleaf{6}{contradiction at $k^{3} - k + 1$}
\schemeitem{2}{$\Delta(k^{3} + k^{2} - k) = 3$}{2.2.2}{9}
\schemeitem{3}{$\Delta(k^{3} - k + 1) = 2$}{2.2.2.1}{5}
\schemeitem{4}{$\Delta(k^{2} + k - 1) = 1$}{2.2.2.1.1}{4}
\schemeitem{5}{$\Delta(k^{2} + 2k - 2) = 2$}{2.2.2.1.1.1}{1}
\schemeleaf{6}{contradiction at $k^{2} + k + 1$}
\schemeitem{5}{$\Delta(k^{2} + 2k - 2) = 3$}{2.2.2.1.1.2}{3}
\schemeitem{6}{$\Delta(k^{2} + k - 2) = 1$}{2.2.2.1.1.2.1}{2}
\schemeitem{7}{$\Delta(3k - 2) = 1$}{2.2.2.1.1.2.1.1}{1}
\schemeleaf{8}{contradiction at $k^{2} + k + 1$}
\schemeitem{7}{$\Delta(3k - 2) = 2$}{2.2.2.1.1.2.1.2}{1}
\schemeleaf{8}{contradiction at $k + 1$}
\schemeitem{6}{$\Delta(k^{2} + k - 2) = 2$}{2.2.2.1.1.2.2}{1}
\schemeleaf{7}{$k\cdot (k + 1) = k^{2} + k$ monochromatic}
\schemeitem{4}{$\Delta(k^{2} + k - 1) = 3$}{2.2.2.1.2}{1}
\schemeleaf{5}{$k\cdot (k^{2} + k - 1) = k^{3} + k^{2} - k$ monochromatic}
\schemeitem{3}{$\Delta(k^{3} - k + 1) = 3$}{2.2.2.2}{4}
\schemeitem{4}{$\Delta(k + 1) = 1$}{2.2.2.2.1}{1}
\schemeleaf{5}{$k\cdot (k + 1) = k^{2} + k$ monochromatic}
\schemeitem{4}{$\Delta(k + 1) = 2$}{2.2.2.2.2}{3}
\schemeitem{5}{$\Delta(k^{2} + k - 1) = 1$}{2.2.2.2.2.1}{2}
\schemeitem{6}{$\Delta(2) = 1$}{2.2.2.2.2.1.1}{1}
\schemeleaf{7}{contradiction at $k - 1$}
\schemeitem{6}{$\Delta(2) = 2$}{2.2.2.2.2.1.2}{1}
\schemeleaf{7}{contradiction at $k^{2} - 1$}
\schemeitem{5}{$\Delta(k^{2} + k - 1) = 3$}{2.2.2.2.2.2}{1}
\schemeleaf{6}{$k\cdot (k^{2} + k - 1) = k^{3} + k^{2} - k$ monochromatic}
\par}

\subsection*{The cases in full}
Throughout, $\Delta(1)=1$ and $\Delta(k)=2$. Since $k\cdot k=k^2$ must not be monochromatic in $A_2$, either $\Delta(k^2)=1$ or $\Delta(k^2)=3$. Every case lists, in parentheses, the elements whose colors are known when the case begins. ``Rule~2 applied to $R_1$ and $R_2$'' means that the element in question lies in both relations and that the other members of $R_1$ and of $R_2$ have two different colors (Rule~\ref{rule2}).

{\small
\case{1}{$\Delta(k^{2}) = 1$}\quad {\small($A_1\supseteq \{1, k^{2}\},\ A_2\supseteq \{k\}$)}

Since $k^{2} + (k-1)\cdot 1 = k^{2} + k - 1$ must not be monochromatic in $A_1$, either $\Delta(k^{2} + k - 1) = 2$ or $\Delta(k^{2} + k - 1) = 3$.

\case{1.1}{$\Delta(k^{2} + k - 1) = 2$}\quad {\small($A_1\supseteq \{1, k^{2}\},\ A_2\supseteq \{k, k^{2} + k - 1\}$)}

By Rule~2 applied to $1 + (k-1)\cdot (k + 1) = k^{2}$ and $k + (k-1)\cdot (k + 1) = k^{2} + k - 1$, $\Delta(k + 1) = 3$.

Since $k\cdot (k + 1) = k^{2} + k$ must not be monochromatic in $A_3$, either $\Delta(k^{2} + k) = 1$ or $\Delta(k^{2} + k) = 2$.

\case{1.1.1}{$\Delta(k^{2} + k) = 1$}\quad {\small($A_1\supseteq \{1, k^{2}, k^{2} + k\},\ A_2\supseteq \{k, k^{2} + k - 1\},\ A_3\supseteq \{k + 1\}$)}

By Rule~2 applied to $1 + (k-1)\cdot (k^{2} + k) = k^{3} - k + 1$ and $k + (k-1)\cdot (k^{2} + k - 1) = k^{3} - k + 1$, $\Delta(k^{3} - k + 1) = 3$.

\emph{Contradiction (Rule~3) at $k^{3} + k^{2} - k$:} if $\Delta(k^{3} + k^{2} - k)=1$ then $k^{2} + (k-1)\cdot (k^{2} + k) = k^{3} + k^{2} - k$ is monochromatic; if $\Delta(k^{3} + k^{2} - k)=2$ then $k\cdot (k^{2} + k - 1) = k^{3} + k^{2} - k$ is monochromatic; if $\Delta(k^{3} + k^{2} - k)=3$ then $(k^{3} - k + 1) + (k-1)\cdot (k + 1) = k^{3} + k^{2} - k$ is monochromatic.

\case{1.1.2}{$\Delta(k^{2} + k) = 2$}\quad {\small($A_1\supseteq \{1, k^{2}\},\ A_2\supseteq \{k, k^{2} + k - 1, k^{2} + k\},\ A_3\supseteq \{k + 1\}$)}

By Rule~2 applied to $k\cdot k^{2} = k^{3}$ and $k + (k-1)\cdot (k^{2} + k) = k^{3}$, $\Delta(k^{3}) = 3$.

By Rule~2 applied to $(k^{2} + k - 1) + (k-1)\cdot (k^{2} + k) = k^{3} + k^{2} - 1$ and $k^{3} + (k-1)\cdot (k + 1) = k^{3} + k^{2} - 1$, $\Delta(k^{3} + k^{2} - 1) = 1$.

By Rule~2 applied to $(2k^{2} - 1) + (k-1)\cdot k^{2} = k^{3} + k^{2} - 1$ and $(k^{2} + k - 1) + (k-1)\cdot k = 2k^{2} - 1$, $\Delta(2k^{2} - 1) = 3$.

By Rule~2 applied to $(k^{3} + k^{2} - k) + (k-1)\cdot 1 = k^{3} + k^{2} - 1$ and $k\cdot (k^{2} + k - 1) = k^{3} + k^{2} - k$, $\Delta(k^{3} + k^{2} - k) = 3$.

By Rule~2 applied to $k + (k-1)\cdot (k^{2} + k - 1) = k^{3} - k + 1$ and $(k^{3} - k + 1) + (k-1)\cdot (k + 1) = k^{3} + k^{2} - k$, $\Delta(k^{3} - k + 1) = 1$.

By Rule~2 applied to $1 + (k-1)\cdot k^{2} = k^{3} - k^{2} + 1$ and $(k^{3} - k^{2} + 1) + (k-1)\cdot (k + 1) = k^{3}$, $\Delta(k^{3} - k^{2} + 1) = 2$.

\emph{Contradiction (Rule~3) at $k^{2} - 1$:} if $\Delta(k^{2} - 1)=1$ then $k^{2} + (k-1)\cdot (k^{2} - 1) = k^{3} - k + 1$ is monochromatic; if $\Delta(k^{2} - 1)=2$ then $k + (k-1)\cdot (k^{2} - 1) = k^{3} - k^{2} + 1$ is monochromatic; if $\Delta(k^{2} - 1)=3$ then $(2k^{2} - 1) + (k-1)\cdot (k^{2} - 1) = k^{3} + k^{2} - k$ is monochromatic.

\case{1.2}{$\Delta(k^{2} + k - 1) = 3$}\quad {\small($A_1\supseteq \{1, k^{2}\},\ A_2\supseteq \{k\},\ A_3\supseteq \{k^{2} + k - 1\}$)}

Since $k\cdot (k^{2} + k - 1) = k^{3} + k^{2} - k$ must not be monochromatic in $A_3$, either $\Delta(k^{3} + k^{2} - k) = 1$ or $\Delta(k^{3} + k^{2} - k) = 2$.

\case{1.2.1}{$\Delta(k^{3} + k^{2} - k) = 1$}\quad {\small($A_1\supseteq \{1, k^{2}, k^{3} + k^{2} - k\},\ A_2\supseteq \{k\},\ A_3\supseteq \{k^{2} + k - 1\}$)}

Since $1 + (k-1)\cdot (k + 1) = k^{2}$ must not be monochromatic in $A_1$, either $\Delta(k + 1) = 2$ or $\Delta(k + 1) = 3$.

\case{1.2.1.1}{$\Delta(k + 1) = 2$}\quad {\small($A_1\supseteq \{1, k^{2}, k^{3} + k^{2} - k\},\ A_2\supseteq \{k, k + 1\},\ A_3\supseteq \{k^{2} + k - 1\}$)}

By Rule~2 applied to $k^{2} + (k-1)\cdot (k^{2} + k) = k^{3} + k^{2} - k$ and $k\cdot (k + 1) = k^{2} + k$, $\Delta(k^{2} + k) = 3$.

By Rule~2 applied to $(k^{3} + k^{2} - k) + (k-1)\cdot 1 = k^{3} + k^{2} - 1$ and $(k^{2} + k - 1) + (k-1)\cdot (k^{2} + k) = k^{3} + k^{2} - 1$, $\Delta(k^{3} + k^{2} - 1) = 2$.

By Rule~2 applied to $k\cdot k^{2} = k^{3}$ and $k^{3} + (k-1)\cdot (k + 1) = k^{3} + k^{2} - 1$, $\Delta(k^{3}) = 3$.

Since $(k^{3} + k^{2} - 2k + 1) + (k-1)\cdot 1 = k^{3} + k^{2} - k$ must not be monochromatic in $A_1$, either $\Delta(k^{3} + k^{2} - 2k + 1) = 2$ or $\Delta(k^{3} + k^{2} - 2k + 1) = 3$.

\case{1.2.1.1.1}{$\Delta(k^{3} + k^{2} - 2k + 1) = 2$}\quad {\small($A_1\supseteq \{1, k^{2}, k^{3} + k^{2} - k\},\ A_2\supseteq \{k, k + 1, k^{3} + k^{2} - 2k + 1, k^{3} + k^{2} - 1\},\ A_3\supseteq \{k^{2} + k - 1, k^{2} + k, k^{3}\}$)}

Since $k + (k-1)\cdot (k^{2} + 2k - 1) = k^{3} + k^{2} - 2k + 1$ must not be monochromatic in $A_2$, either $\Delta(k^{2} + 2k - 1) = 1$ or $\Delta(k^{2} + 2k - 1) = 3$.

\case{1.2.1.1.1.1}{$\Delta(k^{2} + 2k - 1) = 1$}\quad {\small($A_1\supseteq \{1, k^{2}, k^{2} + 2k - 1, k^{3} + k^{2} - k\},\ A_2\supseteq \{k, k + 1, k^{3} + k^{2} - 2k + 1, k^{3} + k^{2} - 1\},\ A_3\supseteq \{k^{2} + k - 1, k^{2} + k, k^{3}\}$)}

By Rule~2 applied to $(2k - 1) + (k-1)\cdot (k^{2} + 2k - 1) = k^{3} + k^{2} - k$ and $(2k - 1) + (k-1)\cdot (k^{2} + k - 1) = k^{3}$, $\Delta(2k - 1) = 2$.

By Rule~2 applied to $k^{2} + (k-1)\cdot (k^{2} + 2k - 1) = k^{3} + 2k^{2} - 3k + 1$ and $(k^{3} + k^{2} - 2k + 1) + (k-1)\cdot k = k^{3} + 2k^{2} - 3k + 1$, $\Delta(k^{3} + 2k^{2} - 3k + 1) = 3$.

By Rule~2 applied to $(2k - 1) + (k-1)\cdot (k + 1) = k^{2} + 2k - 2$ and $(k^{2} + k - 1) + (k-1)\cdot (k^{2} + 2k - 2) = k^{3} + 2k^{2} - 3k + 1$, $\Delta(k^{2} + 2k - 2) = 1$.

By Rule~2 applied to $k^{2} + (k-1)\cdot 2 = k^{2} + 2k - 2$ and $(k^{3} + k^{2} - 2k + 1) + (k-1)\cdot 2 = k^{3} + k^{2} - 1$, $\Delta(2) = 3$.

\emph{Contradiction (Rule~3) at $k^{3} + 2k^{2} - k - 1$:} if $\Delta(k^{3} + 2k^{2} - k - 1)=1$ then $(k^{2} + 2k - 2) + (k-1)\cdot (k^{2} + 2k - 1) = k^{3} + 2k^{2} - k - 1$ is monochromatic; if $\Delta(k^{3} + 2k^{2} - k - 1)=2$ then $(k^{3} + k^{2} - 1) + (k-1)\cdot k = k^{3} + 2k^{2} - k - 1$ is monochromatic; if $\Delta(k^{3} + 2k^{2} - k - 1)=3$ then $(k^{3} + 2k^{2} - 3k + 1) + (k-1)\cdot 2 = k^{3} + 2k^{2} - k - 1$ is monochromatic.

\case{1.2.1.1.1.2}{$\Delta(k^{2} + 2k - 1) = 3$}\quad {\small($A_1\supseteq \{1, k^{2}, k^{3} + k^{2} - k\},\ A_2\supseteq \{k, k + 1, k^{3} + k^{2} - 2k + 1, k^{3} + k^{2} - 1\},\ A_3\supseteq \{k^{2} + k - 1, k^{2} + k, k^{2} + 2k - 1, k^{3}\}$)}

By Rule~2 applied to $(k^{3} + k^{2} - 2k + 1) + (k-1)\cdot (k + 1) = k^{3} + 2k^{2} - 2k$ and $(k^{2} + k - 1) + (k-1)\cdot (k^{2} + 2k - 1) = k^{3} + 2k^{2} - 2k$, $\Delta(k^{3} + 2k^{2} - 2k) = 1$.

By Rule~2 applied to $(k^{3} + 2k^{2} - 3k + 1) + (k-1)\cdot 1 = k^{3} + 2k^{2} - 2k$ and $(k^{3} + k^{2} - 2k + 1) + (k-1)\cdot k = k^{3} + 2k^{2} - 3k + 1$, $\Delta(k^{3} + 2k^{2} - 3k + 1) = 3$.

By Rule~2 applied to $(k^{3} + 2k^{2} - 2k) + (k-1)\cdot 1 = k^{3} + 2k^{2} - k - 1$ and $(k^{3} + k^{2} - 1) + (k-1)\cdot k = k^{3} + 2k^{2} - k - 1$, $\Delta(k^{3} + 2k^{2} - k - 1) = 3$.

By Rule~2 applied to $(k^{3} + k^{2} - 2k + 1) + (k-1)\cdot 2 = k^{3} + k^{2} - 1$ and $(k^{3} + 2k^{2} - 3k + 1) + (k-1)\cdot 2 = k^{3} + 2k^{2} - k - 1$, $\Delta(2) = 1$.

By Rule~2 applied to $1 + (k-1)\cdot 2 = 2k - 1$ and $(2k - 1) + (k-1)\cdot (k^{2} + k - 1) = k^{3}$, $\Delta(2k - 1) = 2$.

\emph{Contradiction (Rule~3) at $k^{2} + 2k$:} if $\Delta(k^{2} + 2k)=1$ then $k^{2} + (k-1)\cdot (k^{2} + 2k) = k^{3} + 2k^{2} - 2k$ is monochromatic; if $\Delta(k^{2} + 2k)=2$ then $(2k - 1) + (k-1)\cdot (k^{2} + 2k) = k^{3} + k^{2} - 1$ is monochromatic; if $\Delta(k^{2} + 2k)=3$ then $(k^{2} + k - 1) + (k-1)\cdot (k^{2} + 2k) = k^{3} + 2k^{2} - k - 1$ is monochromatic.

\case{1.2.1.1.2}{$\Delta(k^{3} + k^{2} - 2k + 1) = 3$}\quad {\small($A_1\supseteq \{1, k^{2}, k^{3} + k^{2} - k\},\ A_2\supseteq \{k, k + 1, k^{3} + k^{2} - 1\},\ A_3\supseteq \{k^{2} + k - 1, k^{2} + k, k^{3}, k^{3} + k^{2} - 2k + 1\}$)}

Since $(k^{2} + k - 1) + (k-1)\cdot (k^{2} - 1) = k^{3}$ must not be monochromatic in $A_3$, either $\Delta(k^{2} - 1) = 1$ or $\Delta(k^{2} - 1) = 2$.

\case{1.2.1.1.2.1}{$\Delta(k^{2} - 1) = 1$}\quad {\small($A_1\supseteq \{1, k^{2} - 1, k^{2}, k^{3} + k^{2} - k\},\ A_2\supseteq \{k, k + 1, k^{3} + k^{2} - 1\},\ A_3\supseteq \{k^{2} + k - 1, k^{2} + k, k^{3}, k^{3} + k^{2} - 2k + 1\}$)}

By Rule~2 applied to $(k^{2} - 1) + (k-1)\cdot 1 = k^{2} + k - 2$ and $(k^{2} + k - 1) + (k-1)\cdot (k^{2} + k - 2) = k^{3} + k^{2} - 2k + 1$, $\Delta(k^{2} + k - 2) = 2$.

By Rule~2 applied to $(k - 1) + (k-1)\cdot (k + 1) = k^{2} + k - 2$ and $k^{3} + (k-1)\cdot (k - 1) = k^{3} + k^{2} - 2k + 1$, $\Delta(k - 1) = 1$.

By Rule~2 applied to $(2k - 1) + (k-1)\cdot (k - 1) = k^{2}$ and $(2k - 1) + (k-1)\cdot (k^{2} + k - 1) = k^{3}$, $\Delta(2k - 1) = 2$.

By Rule~2 applied to $k^{2} + (k-1)\cdot (k^{2} - 1) = k^{3} - k + 1$ and $(2k - 1) + (k-1)\cdot (k^{2} + k - 2) = k^{3} - k + 1$, $\Delta(k^{3} - k + 1) = 3$.

By Rule~2 applied to $(k^{3} + k - 1) + (k-1)\cdot (k - 1) = k^{3} + k^{2} - k$ and $(k^{3} + k - 1) + (k-1)\cdot k = k^{3} + k^{2} - 1$, $\Delta(k^{3} + k - 1) = 3$.

By Rule~2 applied to $k^{2} + (k-1)\cdot (k - 1) = 2k^{2} - 2k + 1$ and $k + (k-1)\cdot (2k - 1) = 2k^{2} - 2k + 1$, $\Delta(2k^{2} - 2k + 1) = 3$.

By Rule~2 applied to $(k^{3} + k^{2} - k) + (k-1)\cdot (k - 1) = k^{3} + 2k^{2} - 3k + 1$ and $(2k^{2} - 2k + 1) + (k-1)\cdot (k^{2} + k) = k^{3} + 2k^{2} - 3k + 1$, $\Delta(k^{3} + 2k^{2} - 3k + 1) = 2$.

By Rule~2 applied to $(2k^{2} - 1) + (k-1)\cdot (k^{2} - 1) = k^{3} + k^{2} - k$ and $(2k^{2} - 1) + (k-1)\cdot (k^{2} + k - 2) = k^{3} + 2k^{2} - 3k + 1$, $\Delta(2k^{2} - 1) = 3$.

By Rule~2 applied to $(k^{3} + k^{2} - 1) + (k-1)\cdot k = k^{3} + 2k^{2} - k - 1$ and $(2k^{2} - 1) + (k-1)\cdot (k^{2} + k) = k^{3} + 2k^{2} - k - 1$, $\Delta(k^{3} + 2k^{2} - k - 1) = 1$.

By Rule~2 applied to $(k^{3} + 2k^{2} - k - 1) + (k-1)\cdot 1 = k^{3} + 2k^{2} - 2$ and $(k^{3} + k^{2} - 1) + (k-1)\cdot (k + 1) = k^{3} + 2k^{2} - 2$, $\Delta(k^{3} + 2k^{2} - 2) = 3$.

By Rule~2 applied to $(k^{3} + 2k^{2} - 2k) + (k-1)\cdot 1 = k^{3} + 2k^{2} - k - 1$ and $(2k^{2} - 1) + (k-1)\cdot (k^{2} + k - 1) = k^{3} + 2k^{2} - 2k$, $\Delta(k^{3} + 2k^{2} - 2k) = 2$.

By Rule~2 applied to $(2k^{2} + k - 2) + (k-1)\cdot (k^{2} + k - 2) = k^{3} + 2k^{2} - 2k$ and $(2k^{2} + k - 2) + (k-1)\cdot (k^{2} + k) = k^{3} + 2k^{2} - 2$, $\Delta(2k^{2} + k - 2) = 1$.

By Rule~2 applied to $(2k^{2} + k - 2) + (k-1)\cdot (k^{2} - 2) = k^{3} + k^{2} - k$ and $(k^{2} + k - 1) + (k-1)\cdot (k^{2} - 2) = k^{3} - k + 1$, $\Delta(k^{2} - 2) = 2$.

By Rule~2 applied to $1 + (k-1)\cdot k^{2} = k^{3} - k^{2} + 1$ and $(2k - 1) + (k-1)\cdot (k^{2} - 2) = k^{3} - k^{2} + 1$, $\Delta(k^{3} - k^{2} + 1) = 3$.

By Rule~2 applied to $k^{2} + (k-1)\cdot (k + 2) = 2k^{2} + k - 2$ and $(k^{3} - k^{2} + 1) + (k-1)\cdot (k + 2) = k^{3} + k - 1$, $\Delta(k + 2) = 2$.

By Rule~2 applied to $(k + 2) + (k-1)\cdot (k + 1) = k^{2} + k + 1$ and $(2k^{2} - 1) + (k-1)\cdot (k^{2} + k + 1) = k^{3} + 2k^{2} - 2$, $\Delta(k^{2} + k + 1) = 1$.

By Rule~2 applied to $(k^{2} + k + 1) + (k-1)\cdot 1 = k^{2} + 2k$ and $(2k - 1) + (k-1)\cdot (k^{2} + 2k) = k^{3} + k^{2} - 1$, $\Delta(k^{2} + 2k) = 3$.

By Rule~2 applied to $(2k - 1) + (k-1)\cdot (k + 1) = k^{2} + 2k - 2$ and $(k^{2} + 2k - 2) + (k-1)\cdot (k^{2} + 2k) = k^{3} + 2k^{2} - 2$, $\Delta(k^{2} + 2k - 2) = 1$.

By Rule~2 applied to $k^{2} + (k-1)\cdot 2 = k^{2} + 2k - 2$ and $(k^{3} - k + 1) + (k-1)\cdot 2 = k^{3} + k - 1$, $\Delta(2) = 2$.

\emph{Contradiction (Rule~3) at $3k - 2$:} if $\Delta(3k - 2)=1$ then $(3k - 2) + (k-1)\cdot (k^{2} + 2k - 2) = k^{3} + k^{2} - k$ is monochromatic; if $\Delta(3k - 2)=2$ then $k + (k-1)\cdot 2 = 3k - 2$ is monochromatic; if $\Delta(3k - 2)=3$ then $(3k - 2) + (k-1)\cdot (k^{2} + k - 1) = k^{3} + k - 1$ is monochromatic.

\case{1.2.1.1.2.2}{$\Delta(k^{2} - 1) = 2$}\quad {\small($A_1\supseteq \{1, k^{2}, k^{3} + k^{2} - k\},\ A_2\supseteq \{k, k + 1, k^{2} - 1, k^{3} + k^{2} - 1\},\ A_3\supseteq \{k^{2} + k - 1, k^{2} + k, k^{3}, k^{3} + k^{2} - 2k + 1\}$)}

By Rule~2 applied to $(k - 1) + (k-1)\cdot k = k^{2} - 1$ and $k^{3} + (k-1)\cdot (k - 1) = k^{3} + k^{2} - 2k + 1$, $\Delta(k - 1) = 1$.

By Rule~2 applied to $1 + (k-1)\cdot k^{2} = k^{3} - k^{2} + 1$ and $k + (k-1)\cdot (k^{2} - 1) = k^{3} - k^{2} + 1$, $\Delta(k^{3} - k^{2} + 1) = 3$.

By Rule~2 applied to $(2k - 1) + (k-1)\cdot (k - 1) = k^{2}$ and $(2k - 1) + (k-1)\cdot (k^{2} + k - 1) = k^{3}$, $\Delta(2k - 1) = 2$.

By Rule~2 applied to $(k^{3} + k - 1) + (k-1)\cdot (k - 1) = k^{3} + k^{2} - k$ and $(k^{3} + k - 1) + (k-1)\cdot k = k^{3} + k^{2} - 1$, $\Delta(k^{3} + k - 1) = 3$.

By Rule~2 applied to $k^{2} + (k-1)\cdot (k - 1) = 2k^{2} - 2k + 1$ and $k + (k-1)\cdot (2k - 1) = 2k^{2} - 2k + 1$, $\Delta(2k^{2} - 2k + 1) = 3$.

Since $(k^{3} + k^{2} - 1) + (k-1)\cdot k = k^{3} + 2k^{2} - k - 1$ must not be monochromatic in $A_2$, either $\Delta(k^{3} + 2k^{2} - k - 1) = 1$ or $\Delta(k^{3} + 2k^{2} - k - 1) = 3$.

\case{1.2.1.1.2.2.1}{$\Delta(k^{3} + 2k^{2} - k - 1) = 1$}\quad {\small($A_1\supseteq \{1, k - 1, k^{2}, k^{3} + k^{2} - k, k^{3} + 2k^{2} - k - 1\},\ A_2\supseteq \{k, k + 1, 2k - 1, k^{2} - 1, k^{3} + k^{2} - 1\},\ A_3\supseteq \{k^{2} + k - 1, k^{2} + k, 2k^{2} - 2k + 1, k^{3} - k^{2} + 1, k^{3}, k^{3} + k - 1, k^{3} + k^{2} - 2k + 1\}$)}

By Rule~2 applied to $(k^{3} + 2k^{2} - k - 1) + (k-1)\cdot 1 = k^{3} + 2k^{2} - 2$ and $(k^{3} + k^{2} - 1) + (k-1)\cdot (k + 1) = k^{3} + 2k^{2} - 2$, $\Delta(k^{3} + 2k^{2} - 2) = 3$.

By Rule~2 applied to $(2k^{2} + k - 2) + (k-1)\cdot (k^{2} - 1) = k^{3} + k^{2} - 1$ and $(2k^{2} + k - 2) + (k-1)\cdot (k^{2} + k) = k^{3} + 2k^{2} - 2$, $\Delta(2k^{2} + k - 2) = 1$.

By Rule~2 applied to $k^{2} + (k-1)\cdot (k + 2) = 2k^{2} + k - 2$ and $(k^{3} - k^{2} + 1) + (k-1)\cdot (k + 2) = k^{3} + k - 1$, $\Delta(k + 2) = 2$.

Since $(k^{3} - k^{2} + 1) + (k-1)\cdot 2k = k^{3} + k^{2} - 2k + 1$ must not be monochromatic in $A_3$, either $\Delta(2k) = 1$ or $\Delta(2k) = 2$.

\case{1.2.1.1.2.2.1.1}{$\Delta(2k) = 1$}\quad {\small($A_1\supseteq \{1, k - 1, 2k, k^{2}, 2k^{2} + k - 2, k^{3} + k^{2} - k, k^{3} + 2k^{2} - k - 1\},\ A_2\supseteq \{k, k + 1, k + 2, 2k - 1, k^{2} - 1, k^{3} + k^{2} - 1\},\ A_3\supseteq \{k^{2} + k - 1, k^{2} + k, 2k^{2} - 2k + 1, k^{3} - k^{2} + 1, k^{3}, k^{3} + k - 1, k^{3} + k^{2} - 2k + 1, k^{3} + 2k^{2} - 2\}$)}

By Rule~2 applied to $(3k - 2) + (k-1)\cdot 2k = 2k^{2} + k - 2$ and $(3k - 2) + (k-1)\cdot (k^{2} + k - 1) = k^{3} + k - 1$, $\Delta(3k - 2) = 2$.

By Rule~2 applied to $(2k^{2} - 1) + (k-1)\cdot 1 = 2k^{2} + k - 2$ and $(3k - 2) + (k-1)\cdot (2k - 1) = 2k^{2} - 1$, $\Delta(2k^{2} - 1) = 3$.

By Rule~2 applied to $(k + 2) + (k-1)\cdot (k + 1) = k^{2} + k + 1$ and $(2k^{2} - 1) + (k-1)\cdot (k^{2} + k + 1) = k^{3} + 2k^{2} - 2$, $\Delta(k^{2} + k + 1) = 1$.

By Rule~2 applied to $(k^{2} + k + 1) + (k-1)\cdot 1 = k^{2} + 2k$ and $k\cdot (k + 2) = k^{2} + 2k$, $\Delta(k^{2} + 2k) = 3$.

By Rule~2 applied to $k + (k-1)\cdot (k + 2) = k^{2} + 2k - 2$ and $(k^{2} + 2k - 2) + (k-1)\cdot (k^{2} + 2k) = k^{3} + 2k^{2} - 2$, $\Delta(k^{2} + 2k - 2) = 1$.

\emph{Contradiction (Rule~3) at $2$:} if $\Delta(2)=1$ then $k^{2} + (k-1)\cdot 2 = k^{2} + 2k - 2$ is monochromatic; if $\Delta(2)=2$ then $k + (k-1)\cdot 2 = 3k - 2$ is monochromatic; if $\Delta(2)=3$ then $(2k^{2} - 2k + 1) + (k-1)\cdot 2 = 2k^{2} - 1$ is monochromatic.

\case{1.2.1.1.2.2.1.2}{$\Delta(2k) = 2$}\quad {\small($A_1\supseteq \{1, k - 1, k^{2}, 2k^{2} + k - 2, k^{3} + k^{2} - k, k^{3} + 2k^{2} - k - 1\},\ A_2\supseteq \{k, k + 1, k + 2, 2k - 1, 2k, k^{2} - 1, k^{3} + k^{2} - 1\},\ A_3\supseteq \{k^{2} + k - 1, k^{2} + k, 2k^{2} - 2k + 1, k^{3} - k^{2} + 1, k^{3}, k^{3} + k - 1, k^{3} + k^{2} - 2k + 1, k^{3} + 2k^{2} - 2\}$)}

By Rule~2 applied to $(2k^{2} - 1) + (k-1)\cdot 1 = 2k^{2} + k - 2$ and $(2k - 1) + (k-1)\cdot 2k = 2k^{2} - 1$, $\Delta(2k^{2} - 1) = 3$.

By Rule~2 applied to $(k + 2) + (k-1)\cdot (k + 1) = k^{2} + k + 1$ and $(2k^{2} - 1) + (k-1)\cdot (k^{2} + k + 1) = k^{3} + 2k^{2} - 2$, $\Delta(k^{2} + k + 1) = 1$.

By Rule~2 applied to $(k^{2} + k + 1) + (k-1)\cdot 1 = k^{2} + 2k$ and $k\cdot (k + 2) = k^{2} + 2k$, $\Delta(k^{2} + 2k) = 3$.

By Rule~2 applied to $k + (k-1)\cdot (k + 2) = k^{2} + 2k - 2$ and $(k^{2} + 2k - 2) + (k-1)\cdot (k^{2} + 2k) = k^{3} + 2k^{2} - 2$, $\Delta(k^{2} + 2k - 2) = 1$.

By Rule~2 applied to $(3k - 2) + (k-1)\cdot (k^{2} + 2k - 2) = k^{3} + k^{2} - k$ and $(3k - 2) + (k-1)\cdot (k^{2} + k - 1) = k^{3} + k - 1$, $\Delta(3k - 2) = 2$.

\emph{Contradiction (Rule~3) at $2$:} if $\Delta(2)=1$ then $k^{2} + (k-1)\cdot 2 = k^{2} + 2k - 2$ is monochromatic; if $\Delta(2)=2$ then $k + (k-1)\cdot 2 = 3k - 2$ is monochromatic; if $\Delta(2)=3$ then $(2k^{2} - 2k + 1) + (k-1)\cdot 2 = 2k^{2} - 1$ is monochromatic.

\case{1.2.1.1.2.2.2}{$\Delta(k^{3} + 2k^{2} - k - 1) = 3$}\quad {\small($A_1\supseteq \{1, k - 1, k^{2}, k^{3} + k^{2} - k\},\ A_2\supseteq \{k, k + 1, 2k - 1, k^{2} - 1, k^{3} + k^{2} - 1\},\ A_3\supseteq \{k^{2} + k - 1, k^{2} + k, 2k^{2} - 2k + 1, k^{3} - k^{2} + 1, k^{3}, k^{3} + k - 1, k^{3} + k^{2} - 2k + 1, k^{3} + 2k^{2} - k - 1\}$)}

By Rule~2 applied to $(2k^{2} + k - 2) + (k-1)\cdot (k^{2} - 1) = k^{3} + k^{2} - 1$ and $(2k^{2} + k - 2) + (k-1)\cdot (k^{2} + k - 1) = k^{3} + 2k^{2} - k - 1$, $\Delta(2k^{2} + k - 2) = 1$.

By Rule~2 applied to $(2k^{2} - 1) + (k-1)\cdot 1 = 2k^{2} + k - 2$ and $(2k^{2} - 1) + (k-1)\cdot (k^{2} + k) = k^{3} + 2k^{2} - k - 1$, $\Delta(2k^{2} - 1) = 2$.

By Rule~2 applied to $(2k - 1) + (k-1)\cdot 2k = 2k^{2} - 1$ and $(k^{3} + k - 1) + (k-1)\cdot 2k = k^{3} + 2k^{2} - k - 1$, $\Delta(2k) = 1$.

\emph{Contradiction (Rule~3) at $3k - 2$:} if $\Delta(3k - 2)=1$ then $(3k - 2) + (k-1)\cdot 2k = 2k^{2} + k - 2$ is monochromatic; if $\Delta(3k - 2)=2$ then $(3k - 2) + (k-1)\cdot (2k - 1) = 2k^{2} - 1$ is monochromatic; if $\Delta(3k - 2)=3$ then $(3k - 2) + (k-1)\cdot (k^{2} + k - 1) = k^{3} + k - 1$ is monochromatic.

\case{1.2.1.2}{$\Delta(k + 1) = 3$}\quad {\small($A_1\supseteq \{1, k^{2}, k^{3} + k^{2} - k\},\ A_2\supseteq \{k\},\ A_3\supseteq \{k + 1, k^{2} + k - 1\}$)}

By Rule~2 applied to $k^{2} + (k-1)\cdot (k^{2} + k) = k^{3} + k^{2} - k$ and $k\cdot (k + 1) = k^{2} + k$, $\Delta(k^{2} + k) = 2$.

By Rule~2 applied to $k\cdot k^{2} = k^{3}$ and $k + (k-1)\cdot (k^{2} + k) = k^{3}$, $\Delta(k^{3}) = 3$.

By Rule~2 applied to $(k^{3} + k^{2} - k) + (k-1)\cdot 1 = k^{3} + k^{2} - 1$ and $k^{3} + (k-1)\cdot (k + 1) = k^{3} + k^{2} - 1$, $\Delta(k^{3} + k^{2} - 1) = 2$.

By Rule~2 applied to $1 + (k-1)\cdot k^{2} = k^{3} - k^{2} + 1$ and $(k^{3} - k^{2} + 1) + (k-1)\cdot (k + 1) = k^{3}$, $\Delta(k^{3} - k^{2} + 1) = 2$.

By Rule~2 applied to $k + (k-1)\cdot (k^{2} - 1) = k^{3} - k^{2} + 1$ and $(k^{2} + k - 1) + (k-1)\cdot (k^{2} - 1) = k^{3}$, $\Delta(k^{2} - 1) = 1$.

By Rule~2 applied to $k^{2} + (k-1)\cdot (k^{2} - 1) = k^{3} - k + 1$ and $(k^{3} - k^{2} + 1) + (k-1)\cdot k = k^{3} - k + 1$, $\Delta(k^{3} - k + 1) = 3$.

Since $(k^{3} + k^{2} - 2k + 1) + (k-1)\cdot 1 = k^{3} + k^{2} - k$ must not be monochromatic in $A_1$, either $\Delta(k^{3} + k^{2} - 2k + 1) = 2$ or $\Delta(k^{3} + k^{2} - 2k + 1) = 3$.

\case{1.2.1.2.1}{$\Delta(k^{3} + k^{2} - 2k + 1) = 2$}\quad {\small($A_1\supseteq \{1, k^{2} - 1, k^{2}, k^{3} + k^{2} - k\},\ A_2\supseteq \{k, k^{2} + k, k^{3} - k^{2} + 1, k^{3} + k^{2} - 2k + 1, k^{3} + k^{2} - 1\},\ A_3\supseteq \{k + 1, k^{2} + k - 1, k^{3} - k + 1, k^{3}\}$)}

Since $(k^{3} + k^{2} - 2k + 1) + (k-1)\cdot 2 = k^{3} + k^{2} - 1$ must not be monochromatic in $A_2$, either $\Delta(2) = 1$ or $\Delta(2) = 3$.

\case{1.2.1.2.1.1}{$\Delta(2) = 1$}\quad {\small($A_1\supseteq \{1, 2, k^{2} - 1, k^{2}, k^{3} + k^{2} - k\},\ A_2\supseteq \{k, k^{2} + k, k^{3} - k^{2} + 1, k^{3} + k^{2} - 2k + 1, k^{3} + k^{2} - 1\},\ A_3\supseteq \{k + 1, k^{2} + k - 1, k^{3} - k + 1, k^{3}\}$)}

By Rule~2 applied to $k\cdot 2 = 2k$ and $(k^{3} - k^{2} + 1) + (k-1)\cdot 2k = k^{3} + k^{2} - 2k + 1$, $\Delta(2k) = 3$.

By Rule~2 applied to $(k^{3} + k^{2} - 2k + 1) + (k-1)\cdot k = k^{3} + 2k^{2} - 3k + 1$ and $(k^{3} - k + 1) + (k-1)\cdot 2k = k^{3} + 2k^{2} - 3k + 1$, $\Delta(k^{3} + 2k^{2} - 3k + 1) = 1$.

\emph{Contradiction (Rule~3) at $k^{2} + 2k - 1$:} if $\Delta(k^{2} + 2k - 1)=1$ then $k^{2} + (k-1)\cdot (k^{2} + 2k - 1) = k^{3} + 2k^{2} - 3k + 1$ is monochromatic; if $\Delta(k^{2} + 2k - 1)=2$ then $k + (k-1)\cdot (k^{2} + 2k - 1) = k^{3} + k^{2} - 2k + 1$ is monochromatic; if $\Delta(k^{2} + 2k - 1)=3$ then $2k + (k-1)\cdot (k + 1) = k^{2} + 2k - 1$ is monochromatic.

\case{1.2.1.2.1.2}{$\Delta(2) = 3$}\quad {\small($A_1\supseteq \{1, k^{2} - 1, k^{2}, k^{3} + k^{2} - k\},\ A_2\supseteq \{k, k^{2} + k, k^{3} - k^{2} + 1, k^{3} + k^{2} - 2k + 1, k^{3} + k^{2} - 1\},\ A_3\supseteq \{2, k + 1, k^{2} + k - 1, k^{3} - k + 1, k^{3}\}$)}

By Rule~2 applied to $(k^{3} - k^{2} + 1) + (k-1)\cdot 2k = k^{3} + k^{2} - 2k + 1$ and $k\cdot 2 = 2k$, $\Delta(2k) = 1$.

By Rule~2 applied to $(k^{3} + k - 1) + (k-1)\cdot k = k^{3} + k^{2} - 1$ and $(k^{3} - k + 1) + (k-1)\cdot 2 = k^{3} + k - 1$, $\Delta(k^{3} + k - 1) = 1$.

By Rule~2 applied to $(k^{3} + k - 1) + (k-1)\cdot 2k = k^{3} + 2k^{2} - k - 1$ and $(k^{3} + k^{2} - 1) + (k-1)\cdot k = k^{3} + 2k^{2} - k - 1$, $\Delta(k^{3} + 2k^{2} - k - 1) = 3$.

By Rule~2 applied to $(k^{3} + k^{2} - 2k + 1) + (k-1)\cdot k = k^{3} + 2k^{2} - 3k + 1$ and $(k^{3} + 2k^{2} - 3k + 1) + (k-1)\cdot 2 = k^{3} + 2k^{2} - k - 1$, $\Delta(k^{3} + 2k^{2} - 3k + 1) = 1$.

By Rule~2 applied to $k^{2} + (k-1)\cdot (k^{2} + 2k - 1) = k^{3} + 2k^{2} - 3k + 1$ and $k + (k-1)\cdot (k^{2} + 2k - 1) = k^{3} + k^{2} - 2k + 1$, $\Delta(k^{2} + 2k - 1) = 3$.

By Rule~2 applied to $(k^{2} + 2k - 2) + (k-1)\cdot (k^{2} - 1) = k^{3} + k - 1$ and $(k^{2} + 2k - 2) + (k-1)\cdot (k^{2} + 2k - 1) = k^{3} + 2k^{2} - k - 1$, $\Delta(k^{2} + 2k - 2) = 2$.

\emph{Contradiction (Rule~3) at $k^{3} + 2k^{2} - 2k$:} if $\Delta(k^{3} + 2k^{2} - 2k)=1$ then $(k^{3} + 2k^{2} - 3k + 1) + (k-1)\cdot 1 = k^{3} + 2k^{2} - 2k$ is monochromatic; if $\Delta(k^{3} + 2k^{2} - 2k)=2$ then $k\cdot (k^{2} + 2k - 2) = k^{3} + 2k^{2} - 2k$ is monochromatic; if $\Delta(k^{3} + 2k^{2} - 2k)=3$ then $(k^{2} + k - 1) + (k-1)\cdot (k^{2} + 2k - 1) = k^{3} + 2k^{2} - 2k$ is monochromatic.

\case{1.2.1.2.2}{$\Delta(k^{3} + k^{2} - 2k + 1) = 3$}\quad {\small($A_1\supseteq \{1, k^{2} - 1, k^{2}, k^{3} + k^{2} - k\},\ A_2\supseteq \{k, k^{2} + k, k^{3} - k^{2} + 1, k^{3} + k^{2} - 1\},\ A_3\supseteq \{k + 1, k^{2} + k - 1, k^{3} - k + 1, k^{3}, k^{3} + k^{2} - 2k + 1\}$)}

By Rule~2 applied to $(k^{2} - 1) + (k-1)\cdot 1 = k^{2} + k - 2$ and $(k^{2} + k - 1) + (k-1)\cdot (k^{2} + k - 2) = k^{3} + k^{2} - 2k + 1$, $\Delta(k^{2} + k - 2) = 2$.

Since $k^{3} + (k-1)\cdot (k - 1) = k^{3} + k^{2} - 2k + 1$ must not be monochromatic in $A_3$, either $\Delta(k - 1) = 1$ or $\Delta(k - 1) = 2$.

\case{1.2.1.2.2.1}{$\Delta(k - 1) = 1$}\quad {\small($A_1\supseteq \{1, k - 1, k^{2} - 1, k^{2}, k^{3} + k^{2} - k\},\ A_2\supseteq \{k, k^{2} + k - 2, k^{2} + k, k^{3} - k^{2} + 1, k^{3} + k^{2} - 1\},\ A_3\supseteq \{k + 1, k^{2} + k - 1, k^{3} - k + 1, k^{3}, k^{3} + k^{2} - 2k + 1\}$)}

By Rule~2 applied to $(2k - 1) + (k-1)\cdot (k - 1) = k^{2}$ and $(2k - 1) + (k-1)\cdot (k^{2} + k - 1) = k^{3}$, $\Delta(2k - 1) = 2$.

By Rule~2 applied to $(2k - 1) + (k-1)\cdot (k^{2} - 2) = k^{3} - k^{2} + 1$ and $(k^{2} + k - 1) + (k-1)\cdot (k^{2} - 2) = k^{3} - k + 1$, $\Delta(k^{2} - 2) = 1$.

By Rule~2 applied to $(2k^{2} + k - 2) + (k-1)\cdot (k^{2} - 2) = k^{3} + k^{2} - k$ and $(k^{2} + k - 1) + (k-1)\cdot (k + 1) = 2k^{2} + k - 2$, $\Delta(2k^{2} + k - 2) = 2$.

By Rule~2 applied to $(2k^{2} + k - 2) + (k-1)\cdot (k^{2} + k - 2) = k^{3} + 2k^{2} - 2k$ and $(k^{3} + k^{2} - 2k + 1) + (k-1)\cdot (k + 1) = k^{3} + 2k^{2} - 2k$, $\Delta(k^{3} + 2k^{2} - 2k) = 1$.

By Rule~2 applied to $(k^{3} + 2k^{2} - 2k) + (k-1)\cdot 1 = k^{3} + 2k^{2} - k - 1$ and $(k^{3} + k^{2} - 1) + (k-1)\cdot k = k^{3} + 2k^{2} - k - 1$, $\Delta(k^{3} + 2k^{2} - k - 1) = 3$.

\emph{Contradiction (Rule~3) at $k^{2} + 2k$:} if $\Delta(k^{2} + 2k)=1$ then $k^{2} + (k-1)\cdot (k^{2} + 2k) = k^{3} + 2k^{2} - 2k$ is monochromatic; if $\Delta(k^{2} + 2k)=2$ then $(2k - 1) + (k-1)\cdot (k^{2} + 2k) = k^{3} + k^{2} - 1$ is monochromatic; if $\Delta(k^{2} + 2k)=3$ then $(k^{2} + k - 1) + (k-1)\cdot (k^{2} + 2k) = k^{3} + 2k^{2} - k - 1$ is monochromatic.

\case{1.2.1.2.2.2}{$\Delta(k - 1) = 2$}\quad {\small($A_1\supseteq \{1, k^{2} - 1, k^{2}, k^{3} + k^{2} - k\},\ A_2\supseteq \{k - 1, k, k^{2} + k - 2, k^{2} + k, k^{3} - k^{2} + 1, k^{3} + k^{2} - 1\},\ A_3\supseteq \{k + 1, k^{2} + k - 1, k^{3} - k + 1, k^{3}, k^{3} + k^{2} - 2k + 1\}$)}

By Rule~2 applied to $(k^{3} + k^{2} - 1) + (k-1)\cdot (k - 1) = k^{3} + 2k^{2} - 2k$ and $(k^{3} + k^{2} - 2k + 1) + (k-1)\cdot (k + 1) = k^{3} + 2k^{2} - 2k$, $\Delta(k^{3} + 2k^{2} - 2k) = 1$.

By Rule~2 applied to $(k^{3} + 2k^{2} - 2k) + (k-1)\cdot 1 = k^{3} + 2k^{2} - k - 1$ and $(k^{3} + k^{2} - 1) + (k-1)\cdot k = k^{3} + 2k^{2} - k - 1$, $\Delta(k^{3} + 2k^{2} - k - 1) = 3$.

By Rule~2 applied to $k^{2} + (k-1)\cdot (k^{2} + 2k) = k^{3} + 2k^{2} - 2k$ and $(k^{2} + k - 1) + (k-1)\cdot (k^{2} + 2k) = k^{3} + 2k^{2} - k - 1$, $\Delta(k^{2} + 2k) = 2$.

By Rule~2 applied to $(2k - 1) + (k-1)\cdot (k^{2} + 2k) = k^{3} + k^{2} - 1$ and $(2k - 1) + (k-1)\cdot (k^{2} + k - 1) = k^{3}$, $\Delta(2k - 1) = 1$.

By Rule~2 applied to $(k^{3} + k - 1) + (k-1)\cdot (2k - 1) = k^{3} + 2k^{2} - 2k$ and $(k^{3} + k - 1) + (k-1)\cdot k = k^{3} + k^{2} - 1$, $\Delta(k^{3} + k - 1) = 3$.

By Rule~2 applied to $(2k - 1) + (k-1)\cdot 1 = 3k - 2$ and $(3k - 2) + (k-1)\cdot (k^{2} + k - 1) = k^{3} + k - 1$, $\Delta(3k - 2) = 2$.

\emph{Contradiction (Rule~3) at $2$:} if $\Delta(2)=1$ then $1 + (k-1)\cdot 2 = 2k - 1$ is monochromatic; if $\Delta(2)=2$ then $k + (k-1)\cdot 2 = 3k - 2$ is monochromatic; if $\Delta(2)=3$ then $(k^{3} - k + 1) + (k-1)\cdot 2 = k^{3} + k - 1$ is monochromatic.

\case{1.2.2}{$\Delta(k^{3} + k^{2} - k) = 2$}\quad {\small($A_1\supseteq \{1, k^{2}\},\ A_2\supseteq \{k, k^{3} + k^{2} - k\},\ A_3\supseteq \{k^{2} + k - 1\}$)}

By Rule~2 applied to $k\cdot k^{2} = k^{3}$ and $k^{3} + (k-1)\cdot k = k^{3} + k^{2} - k$, $\Delta(k^{3}) = 3$.

Since $(k^{3} + k^{2} - k) + (k-1)\cdot k = k^{3} + 2k^{2} - 2k$ must not be monochromatic in $A_2$, either $\Delta(k^{3} + 2k^{2} - 2k) = 1$ or $\Delta(k^{3} + 2k^{2} - 2k) = 3$.

\case{1.2.2.1}{$\Delta(k^{3} + 2k^{2} - 2k) = 1$}\quad {\small($A_1\supseteq \{1, k^{2}, k^{3} + 2k^{2} - 2k\},\ A_2\supseteq \{k, k^{3} + k^{2} - k\},\ A_3\supseteq \{k^{2} + k - 1, k^{3}\}$)}

By Rule~2 applied to $k^{2} + (k-1)\cdot (k^{2} + 2k) = k^{3} + 2k^{2} - 2k$ and $k + (k-1)\cdot (k^{2} + 2k) = k^{3} + k^{2} - k$, $\Delta(k^{2} + 2k) = 3$.

By Rule~2 applied to $(k^{3} + 2k^{2} - 2k) + (k-1)\cdot 1 = k^{3} + 2k^{2} - k - 1$ and $(k^{2} + k - 1) + (k-1)\cdot (k^{2} + 2k) = k^{3} + 2k^{2} - k - 1$, $\Delta(k^{3} + 2k^{2} - k - 1) = 2$.

By Rule~2 applied to $1 + (k-1)\cdot (k + 1) = k^{2}$ and $(k^{3} + k^{2} - k) + (k-1)\cdot (k + 1) = k^{3} + 2k^{2} - k - 1$, $\Delta(k + 1) = 3$.

By Rule~2 applied to $1 + (k-1)\cdot k^{2} = k^{3} - k^{2} + 1$ and $(k^{3} - k^{2} + 1) + (k-1)\cdot (k + 1) = k^{3}$, $\Delta(k^{3} - k^{2} + 1) = 2$.

By Rule~2 applied to $k + (k-1)\cdot (k^{2} - 1) = k^{3} - k^{2} + 1$ and $(k^{2} + k - 1) + (k-1)\cdot (k^{2} - 1) = k^{3}$, $\Delta(k^{2} - 1) = 1$.

By Rule~2 applied to $(k^{3} + k^{2} - 1) + (k-1)\cdot k = k^{3} + 2k^{2} - k - 1$ and $k^{3} + (k-1)\cdot (k + 1) = k^{3} + k^{2} - 1$, $\Delta(k^{3} + k^{2} - 1) = 1$.

By Rule~2 applied to $(2k^{2} + k - 2) + (k-1)\cdot (k^{2} - 1) = k^{3} + k^{2} - 1$ and $(k^{2} + k - 1) + (k-1)\cdot (k + 1) = 2k^{2} + k - 2$, $\Delta(2k^{2} + k - 2) = 2$.

Since $(k^{2} + 2k - 2) + (k-1)\cdot k = 2k^{2} + k - 2$ must not be monochromatic in $A_2$, either $\Delta(k^{2} + 2k - 2) = 1$ or $\Delta(k^{2} + 2k - 2) = 3$.

\case{1.2.2.1.1}{$\Delta(k^{2} + 2k - 2) = 1$}\quad {\small($A_1\supseteq \{1, k^{2} - 1, k^{2}, k^{2} + 2k - 2, k^{3} + k^{2} - 1, k^{3} + 2k^{2} - 2k\},\ A_2\supseteq \{k, 2k^{2} + k - 2, k^{3} - k^{2} + 1, k^{3} + k^{2} - k, k^{3} + 2k^{2} - k - 1\},\ A_3\supseteq \{k + 1, k^{2} + k - 1, k^{2} + 2k, k^{3}\}$)}

\emph{Contradiction:} $k\cdot (k^{2} + 2k - 2) = k^{3} + 2k^{2} - 2k$ is now monochromatic.

\case{1.2.2.1.2}{$\Delta(k^{2} + 2k - 2) = 3$}\quad {\small($A_1\supseteq \{1, k^{2} - 1, k^{2}, k^{3} + k^{2} - 1, k^{3} + 2k^{2} - 2k\},\ A_2\supseteq \{k, 2k^{2} + k - 2, k^{3} - k^{2} + 1, k^{3} + k^{2} - k, k^{3} + 2k^{2} - k - 1\},\ A_3\supseteq \{k + 1, k^{2} + k - 1, k^{2} + 2k - 2, k^{2} + 2k, k^{3}\}$)}

By Rule~2 applied to $(k^{3} + 2k^{2} - 3k + 1) + (k-1)\cdot 1 = k^{3} + 2k^{2} - 2k$ and $(k^{2} + k - 1) + (k-1)\cdot (k^{2} + 2k - 2) = k^{3} + 2k^{2} - 3k + 1$, $\Delta(k^{3} + 2k^{2} - 3k + 1) = 2$.

By Rule~2 applied to $(k^{3} + k^{2} - 1) + (k-1)\cdot (k - 1) = k^{3} + 2k^{2} - 2k$ and $(k^{3} + k^{2} - k) + (k-1)\cdot (k - 1) = k^{3} + 2k^{2} - 3k + 1$, $\Delta(k - 1) = 3$.

By Rule~2 applied to $(2k^{2} - 1) + (k-1)\cdot k^{2} = k^{3} + k^{2} - 1$ and $(k^{2} + 2k - 2) + (k-1)\cdot (k - 1) = 2k^{2} - 1$, $\Delta(2k^{2} - 1) = 2$.

\emph{Contradiction (Rule~3) at $k^{2} + k - 2$:} if $\Delta(k^{2} + k - 2)=1$ then $(k^{2} - 1) + (k-1)\cdot 1 = k^{2} + k - 2$ is monochromatic; if $\Delta(k^{2} + k - 2)=2$ then $(2k^{2} - 1) + (k-1)\cdot (k^{2} + k - 2) = k^{3} + 2k^{2} - 3k + 1$ is monochromatic; if $\Delta(k^{2} + k - 2)=3$ then $(k - 1) + (k-1)\cdot (k + 1) = k^{2} + k - 2$ is monochromatic.

\case{1.2.2.2}{$\Delta(k^{3} + 2k^{2} - 2k) = 3$}\quad {\small($A_1\supseteq \{1, k^{2}\},\ A_2\supseteq \{k, k^{3} + k^{2} - k\},\ A_3\supseteq \{k^{2} + k - 1, k^{3}, k^{3} + 2k^{2} - 2k\}$)}

Since $(2k - 1) + (k-1)\cdot (k^{2} + k - 1) = k^{3}$ must not be monochromatic in $A_3$, either $\Delta(2k - 1) = 1$ or $\Delta(2k - 1) = 2$.

\case{1.2.2.2.1}{$\Delta(2k - 1) = 1$}\quad {\small($A_1\supseteq \{1, 2k - 1, k^{2}\},\ A_2\supseteq \{k, k^{3} + k^{2} - k\},\ A_3\supseteq \{k^{2} + k - 1, k^{3}, k^{3} + 2k^{2} - 2k\}$)}

Since $(2k - 1) + (k-1)\cdot 1 = 3k - 2$ must not be monochromatic in $A_1$, either $\Delta(3k - 2) = 2$ or $\Delta(3k - 2) = 3$.

\case{1.2.2.2.1.1}{$\Delta(3k - 2) = 2$}\quad {\small($A_1\supseteq \{1, 2k - 1, k^{2}\},\ A_2\supseteq \{k, 3k - 2, k^{3} + k^{2} - k\},\ A_3\supseteq \{k^{2} + k - 1, k^{3}, k^{3} + 2k^{2} - 2k\}$)}

By Rule~2 applied to $1 + (k-1)\cdot 2 = 2k - 1$ and $k + (k-1)\cdot 2 = 3k - 2$, $\Delta(2) = 3$.

Since $(3k - 2) + (k-1)\cdot k = k^{2} + 2k - 2$ must not be monochromatic in $A_2$, either $\Delta(k^{2} + 2k - 2) = 1$ or $\Delta(k^{2} + 2k - 2) = 3$.

\case{1.2.2.2.1.1.1}{$\Delta(k^{2} + 2k - 2) = 1$}\quad {\small($A_1\supseteq \{1, 2k - 1, k^{2}, k^{2} + 2k - 2\},\ A_2\supseteq \{k, 3k - 2, k^{3} + k^{2} - k\},\ A_3\supseteq \{2, k^{2} + k - 1, k^{3}, k^{3} + 2k^{2} - 2k\}$)}

Since $(2k - 1) + (k-1)\cdot (k - 1) = k^{2}$ must not be monochromatic in $A_1$, either $\Delta(k - 1) = 2$ or $\Delta(k - 1) = 3$.

\case{1.2.2.2.1.1.1.1}{$\Delta(k - 1) = 2$}\quad {\small($A_1\supseteq \{1, 2k - 1, k^{2}, k^{2} + 2k - 2\},\ A_2\supseteq \{k - 1, k, 3k - 2, k^{3} + k^{2} - k\},\ A_3\supseteq \{2, k^{2} + k - 1, k^{3}, k^{3} + 2k^{2} - 2k\}$)}

By Rule~2 applied to $(k - 1) + (k-1)\cdot k = k^{2} - 1$ and $(k^{2} + k - 1) + (k-1)\cdot (k^{2} - 1) = k^{3}$, $\Delta(k^{2} - 1) = 1$.

By Rule~2 applied to $(k^{2} + 2k - 2) + (k-1)\cdot (k^{2} - 1) = k^{3} + k - 1$ and $(k^{3} + k - 1) + (k-1)\cdot (k - 1) = k^{3} + k^{2} - k$, $\Delta(k^{3} + k - 1) = 3$.

By Rule~2 applied to $k^{2} + (k-1)\cdot (k^{2} - 1) = k^{3} - k + 1$ and $(k^{3} - k + 1) + (k-1)\cdot 2 = k^{3} + k - 1$, $\Delta(k^{3} - k + 1) = 2$.

By Rule~2 applied to $(2k - 1) + (k-1)\cdot (k + 1) = k^{2} + 2k - 2$ and $(k^{3} - k + 1) + (k-1)\cdot (k + 1) = k^{3} + k^{2} - k$, $\Delta(k + 1) = 3$.

\emph{Contradiction (Rule~3) at $k^{3} + k^{2} - 2k + 1$:} if $\Delta(k^{3} + k^{2} - 2k + 1)=1$ then $(2k - 1) + (k-1)\cdot (k^{2} + 2k - 2) = k^{3} + k^{2} - 2k + 1$ is monochromatic; if $\Delta(k^{3} + k^{2} - 2k + 1)=2$ then $(k^{3} - k + 1) + (k-1)\cdot k = k^{3} + k^{2} - 2k + 1$ is monochromatic; if $\Delta(k^{3} + k^{2} - 2k + 1)=3$ then $(k^{3} + k^{2} - 2k + 1) + (k-1)\cdot (k + 1) = k^{3} + 2k^{2} - 2k$ is monochromatic.

\case{1.2.2.2.1.1.1.2}{$\Delta(k - 1) = 3$}\quad {\small($A_1\supseteq \{1, 2k - 1, k^{2}, k^{2} + 2k - 2\},\ A_2\supseteq \{k, 3k - 2, k^{3} + k^{2} - k\},\ A_3\supseteq \{2, k - 1, k^{2} + k - 1, k^{3}, k^{3} + 2k^{2} - 2k\}$)}

By Rule~2 applied to $(2k - 1) + (k-1)\cdot (k^{2} + 2k - 2) = k^{3} + k^{2} - 2k + 1$ and $k^{3} + (k-1)\cdot (k - 1) = k^{3} + k^{2} - 2k + 1$, $\Delta(k^{3} + k^{2} - 2k + 1) = 2$.

By Rule~2 applied to $k + (k-1)\cdot (k^{2} + 2k - 1) = k^{3} + k^{2} - 2k + 1$ and $(k^{2} + k - 1) + (k-1)\cdot (k^{2} + 2k - 1) = k^{3} + 2k^{2} - 2k$, $\Delta(k^{2} + 2k - 1) = 1$.

By Rule~2 applied to $k^{2} + (k-1)\cdot (k^{2} + 2k - 1) = k^{3} + 2k^{2} - 3k + 1$ and $(k^{3} + k^{2} - 2k + 1) + (k-1)\cdot k = k^{3} + 2k^{2} - 3k + 1$, $\Delta(k^{3} + 2k^{2} - 3k + 1) = 3$.

By Rule~2 applied to $(k^{2} + 2k - 2) + (k-1)\cdot (k^{2} + 2k - 1) = k^{3} + 2k^{2} - k - 1$ and $(k^{3} + 2k^{2} - 3k + 1) + (k-1)\cdot 2 = k^{3} + 2k^{2} - k - 1$, $\Delta(k^{3} + 2k^{2} - k - 1) = 2$.

By Rule~2 applied to $(k^{3} + k^{2} - 1) + (k-1)\cdot k = k^{3} + 2k^{2} - k - 1$ and $(k^{3} + k^{2} - 1) + (k-1)\cdot (k - 1) = k^{3} + 2k^{2} - 2k$, $\Delta(k^{3} + k^{2} - 1) = 1$.

By Rule~2 applied to $(2k - 1) + (k-1)\cdot (k + 1) = k^{2} + 2k - 2$ and $(k^{3} + k^{2} - k) + (k-1)\cdot (k + 1) = k^{3} + 2k^{2} - k - 1$, $\Delta(k + 1) = 3$.

By Rule~2 applied to $(2k^{2} - 1) + (k-1)\cdot k^{2} = k^{3} + k^{2} - 1$ and $(2k^{2} - 1) + (k-1)\cdot (k^{2} + k - 1) = k^{3} + 2k^{2} - 2k$, $\Delta(2k^{2} - 1) = 2$.

\emph{Contradiction (Rule~3) at $k^{2} + k$:} if $\Delta(k^{2} + k)=1$ then $(k^{2} + k) + (k-1)\cdot 1 = k^{2} + 2k - 1$ is monochromatic; if $\Delta(k^{2} + k)=2$ then $(2k^{2} - 1) + (k-1)\cdot (k^{2} + k) = k^{3} + 2k^{2} - k - 1$ is monochromatic; if $\Delta(k^{2} + k)=3$ then $k\cdot (k + 1) = k^{2} + k$ is monochromatic.

\case{1.2.2.2.1.1.2}{$\Delta(k^{2} + 2k - 2) = 3$}\quad {\small($A_1\supseteq \{1, 2k - 1, k^{2}\},\ A_2\supseteq \{k, 3k - 2, k^{3} + k^{2} - k\},\ A_3\supseteq \{2, k^{2} + k - 1, k^{2} + 2k - 2, k^{3}, k^{3} + 2k^{2} - 2k\}$)}

\emph{Contradiction:} $k\cdot (k^{2} + 2k - 2) = k^{3} + 2k^{2} - 2k$ is now monochromatic.

\case{1.2.2.2.1.2}{$\Delta(3k - 2) = 3$}\quad {\small($A_1\supseteq \{1, 2k - 1, k^{2}\},\ A_2\supseteq \{k, k^{3} + k^{2} - k\},\ A_3\supseteq \{3k - 2, k^{2} + k - 1, k^{3}, k^{3} + 2k^{2} - 2k\}$)}

By Rule~2 applied to $(2k - 1) + (k-1)\cdot (k - 1) = k^{2}$ and $(3k - 2) + (k-1)\cdot (k - 1) = k^{2} + k - 1$, $\Delta(k - 1) = 2$.

By Rule~2 applied to $(k^{3} + k - 1) + (k-1)\cdot (k - 1) = k^{3} + k^{2} - k$ and $(3k - 2) + (k-1)\cdot (k^{2} + k - 1) = k^{3} + k - 1$, $\Delta(k^{3} + k - 1) = 1$.

By Rule~2 applied to $(k - 1) + (k-1)\cdot k = k^{2} - 1$ and $(k^{2} + k - 1) + (k-1)\cdot (k^{2} - 1) = k^{3}$, $\Delta(k^{2} - 1) = 1$.

By Rule~2 applied to $(k^{2} - 1) + (k-1)\cdot 1 = k^{2} + k - 2$ and $(3k - 2) + (k-1)\cdot (k^{2} + k - 2) = k^{3}$, $\Delta(k^{2} + k - 2) = 2$.

By Rule~2 applied to $(k^{2} + 2k - 2) + (k-1)\cdot (k^{2} - 1) = k^{3} + k - 1$ and $(k^{2} + 2k - 2) + (k-1)\cdot (k^{2} + k - 2) = k^{3} + k^{2} - k$, $\Delta(k^{2} + 2k - 2) = 3$.

\emph{Contradiction:} $k\cdot (k^{2} + 2k - 2) = k^{3} + 2k^{2} - 2k$ is now monochromatic.

\case{1.2.2.2.2}{$\Delta(2k - 1) = 2$}\quad {\small($A_1\supseteq \{1, k^{2}\},\ A_2\supseteq \{k, 2k - 1, k^{3} + k^{2} - k\},\ A_3\supseteq \{k^{2} + k - 1, k^{3}, k^{3} + 2k^{2} - 2k\}$)}

By Rule~2 applied to $(2k - 1) + (k-1)\cdot (k^{2} + 2k - 1) = k^{3} + k^{2} - k$ and $(k^{2} + k - 1) + (k-1)\cdot (k^{2} + 2k - 1) = k^{3} + 2k^{2} - 2k$, $\Delta(k^{2} + 2k - 1) = 1$.

Since $1 + (k-1)\cdot (k + 1) = k^{2}$ must not be monochromatic in $A_1$, either $\Delta(k + 1) = 2$ or $\Delta(k + 1) = 3$.

\case{1.2.2.2.2.1}{$\Delta(k + 1) = 2$}\quad {\small($A_1\supseteq \{1, k^{2}, k^{2} + 2k - 1\},\ A_2\supseteq \{k, k + 1, 2k - 1, k^{3} + k^{2} - k\},\ A_3\supseteq \{k^{2} + k - 1, k^{3}, k^{3} + 2k^{2} - 2k\}$)}

By Rule~2 applied to $(k^{2} + k) + (k-1)\cdot 1 = k^{2} + 2k - 1$ and $k\cdot (k + 1) = k^{2} + k$, $\Delta(k^{2} + k) = 3$.

Since $(2k - 1) + (k-1)\cdot (k + 1) = k^{2} + 2k - 2$ must not be monochromatic in $A_2$, either $\Delta(k^{2} + 2k - 2) = 1$ or $\Delta(k^{2} + 2k - 2) = 3$.

\case{1.2.2.2.2.1.1}{$\Delta(k^{2} + 2k - 2) = 1$}\quad {\small($A_1\supseteq \{1, k^{2}, k^{2} + 2k - 2, k^{2} + 2k - 1\},\ A_2\supseteq \{k, k + 1, 2k - 1, k^{3} + k^{2} - k\},\ A_3\supseteq \{k^{2} + k - 1, k^{2} + k, k^{3}, k^{3} + 2k^{2} - 2k\}$)}

By Rule~2 applied to $(k^{2} + 2k - 2) + (k-1)\cdot (k^{2} + 2k - 1) = k^{3} + 2k^{2} - k - 1$ and $(k^{3} + k^{2} - k) + (k-1)\cdot (k + 1) = k^{3} + 2k^{2} - k - 1$, $\Delta(k^{3} + 2k^{2} - k - 1) = 3$.

By Rule~2 applied to $k + (k-1)\cdot (k^{2} + 2k) = k^{3} + k^{2} - k$ and $(k^{2} + k - 1) + (k-1)\cdot (k^{2} + 2k) = k^{3} + 2k^{2} - k - 1$, $\Delta(k^{2} + 2k) = 1$.

Since $k^{2} + (k-1)\cdot 2 = k^{2} + 2k - 2$ must not be monochromatic in $A_1$, either $\Delta(2) = 2$ or $\Delta(2) = 3$.

\case{1.2.2.2.2.1.1.1}{$\Delta(2) = 2$}\quad {\small($A_1\supseteq \{1, k^{2}, k^{2} + 2k - 2, k^{2} + 2k - 1, k^{2} + 2k\},\ A_2\supseteq \{2, k, k + 1, 2k - 1, k^{3} + k^{2} - k\},\ A_3\supseteq \{k^{2} + k - 1, k^{2} + k, k^{3}, k^{3} + 2k^{2} - 2k, k^{3} + 2k^{2} - k - 1\}$)}

By Rule~2 applied to $k\cdot 2 = 2k$ and $k^{3} + (k-1)\cdot 2k = k^{3} + 2k^{2} - 2k$, $\Delta(2k) = 1$.

By Rule~2 applied to $1 + (k-1)\cdot 2k = 2k^{2} - 2k + 1$ and $k + (k-1)\cdot (2k - 1) = 2k^{2} - 2k + 1$, $\Delta(2k^{2} - 2k + 1) = 3$.

By Rule~2 applied to $(k^{2} + k + 1) + (k-1)\cdot 1 = k^{2} + 2k$ and $(2k^{2} - 2k + 1) + (k-1)\cdot (k^{2} + k + 1) = k^{3} + 2k^{2} - 2k$, $\Delta(k^{2} + k + 1) = 2$.

By Rule~2 applied to $k^{2} + (k-1)\cdot (k^{2} + 2k - 1) = k^{3} + 2k^{2} - 3k + 1$ and $(2k^{2} - 2k + 1) + (k-1)\cdot (k^{2} + k) = k^{3} + 2k^{2} - 3k + 1$, $\Delta(k^{3} + 2k^{2} - 3k + 1) = 2$.

By Rule~2 applied to $1 + (k-1)\cdot (k^{2} + 2k) = k^{3} + k^{2} - 2k + 1$ and $(k^{3} + k^{2} - 2k + 1) + (k-1)\cdot k = k^{3} + 2k^{2} - 3k + 1$, $\Delta(k^{3} + k^{2} - 2k + 1) = 3$.

By Rule~2 applied to $(k + 2) + (k-1)\cdot (k + 1) = k^{2} + k + 1$ and $(k^{3} + k^{2} - 2k + 1) + (k-1)\cdot (k + 2) = k^{3} + 2k^{2} - k - 1$, $\Delta(k + 2) = 1$.

\emph{Contradiction:} $k\cdot (k + 2) = k^{2} + 2k$ is now monochromatic.

\case{1.2.2.2.2.1.1.2}{$\Delta(2) = 3$}\quad {\small($A_1\supseteq \{1, k^{2}, k^{2} + 2k - 2, k^{2} + 2k - 1, k^{2} + 2k\},\ A_2\supseteq \{k, k + 1, 2k - 1, k^{3} + k^{2} - k\},\ A_3\supseteq \{2, k^{2} + k - 1, k^{2} + k, k^{3}, k^{3} + 2k^{2} - 2k, k^{3} + 2k^{2} - k - 1\}$)}

By Rule~2 applied to $(k^{2} + 2k - 2) + (k-1)\cdot (k^{2} + 2k) = k^{3} + 2k^{2} - 2$ and $(k^{3} + 2k^{2} - 2k) + (k-1)\cdot 2 = k^{3} + 2k^{2} - 2$, $\Delta(k^{3} + 2k^{2} - 2) = 2$.

By Rule~2 applied to $(k^{3} + k^{2} - k) + (k-1)\cdot (k + 2) = k^{3} + 2k^{2} - 2$ and $2 + (k-1)\cdot (k + 2) = k^{2} + k$, $\Delta(k + 2) = 1$.

\emph{Contradiction:} $k\cdot (k + 2) = k^{2} + 2k$ is now monochromatic.

\case{1.2.2.2.2.1.2}{$\Delta(k^{2} + 2k - 2) = 3$}\quad {\small($A_1\supseteq \{1, k^{2}, k^{2} + 2k - 1\},\ A_2\supseteq \{k, k + 1, 2k - 1, k^{3} + k^{2} - k\},\ A_3\supseteq \{k^{2} + k - 1, k^{2} + k, k^{2} + 2k - 2, k^{3}, k^{3} + 2k^{2} - 2k\}$)}

\emph{Contradiction:} $k\cdot (k^{2} + 2k - 2) = k^{3} + 2k^{2} - 2k$ is now monochromatic.

\case{1.2.2.2.2.2}{$\Delta(k + 1) = 3$}\quad {\small($A_1\supseteq \{1, k^{2}, k^{2} + 2k - 1\},\ A_2\supseteq \{k, 2k - 1, k^{3} + k^{2} - k\},\ A_3\supseteq \{k + 1, k^{2} + k - 1, k^{3}, k^{3} + 2k^{2} - 2k\}$)}

By Rule~2 applied to $(k^{2} + k) + (k-1)\cdot 1 = k^{2} + 2k - 1$ and $k\cdot (k + 1) = k^{2} + k$, $\Delta(k^{2} + k) = 2$.

By Rule~2 applied to $1 + (k-1)\cdot k^{2} = k^{3} - k^{2} + 1$ and $(k^{3} - k^{2} + 1) + (k-1)\cdot (k + 1) = k^{3}$, $\Delta(k^{3} - k^{2} + 1) = 2$.

By Rule~2 applied to $2k + (k-1)\cdot k = k^{2} + k$ and $k^{3} + (k-1)\cdot 2k = k^{3} + 2k^{2} - 2k$, $\Delta(2k) = 1$.

By Rule~2 applied to $k + (k-1)\cdot (k^{2} - 1) = k^{3} - k^{2} + 1$ and $(k^{2} + k - 1) + (k-1)\cdot (k^{2} - 1) = k^{3}$, $\Delta(k^{2} - 1) = 1$.

By Rule~2 applied to $k^{2} + (k-1)\cdot (k^{2} - 1) = k^{3} - k + 1$ and $(k^{3} - k^{2} + 1) + (k-1)\cdot k = k^{3} - k + 1$, $\Delta(k^{3} - k + 1) = 3$.

By Rule~2 applied to $(2k - 1) + (k-1)\cdot (k^{2} - 2) = k^{3} - k^{2} + 1$ and $(k^{2} + k - 1) + (k-1)\cdot (k^{2} - 2) = k^{3} - k + 1$, $\Delta(k^{2} - 2) = 1$.

Since $(k^{3} + k^{2} - 2k + 1) + (k-1)\cdot (k + 1) = k^{3} + 2k^{2} - 2k$ must not be monochromatic in $A_3$, either $\Delta(k^{3} + k^{2} - 2k + 1) = 1$ or $\Delta(k^{3} + k^{2} - 2k + 1) = 2$.

\case{1.2.2.2.2.2.1}{$\Delta(k^{3} + k^{2} - 2k + 1) = 1$}\quad {\small($A_1\supseteq \{1, 2k, k^{2} - 2, k^{2} - 1, k^{2}, k^{2} + 2k - 1, k^{3} + k^{2} - 2k + 1\},\ A_2\supseteq \{k, 2k - 1, k^{2} + k, k^{3} - k^{2} + 1, k^{3} + k^{2} - k\},\ A_3\supseteq \{k + 1, k^{2} + k - 1, k^{3} - k + 1, k^{3}, k^{3} + 2k^{2} - 2k\}$)}

By Rule~2 applied to $(2k^{2} - 1) + (k-1)\cdot (k^{2} - 2) = k^{3} + k^{2} - 2k + 1$ and $(2k^{2} - 1) + (k-1)\cdot (k^{2} + k - 1) = k^{3} + 2k^{2} - 2k$, $\Delta(2k^{2} - 1) = 2$.

By Rule~2 applied to $1 + (k-1)\cdot (k^{2} + 2k) = k^{3} + k^{2} - 2k + 1$ and $k + (k-1)\cdot (k^{2} + 2k) = k^{3} + k^{2} - k$, $\Delta(k^{2} + 2k) = 3$.

By Rule~2 applied to $(2k^{2} - 1) + (k-1)\cdot (k^{2} + k) = k^{3} + 2k^{2} - k - 1$ and $(k^{2} + k - 1) + (k-1)\cdot (k^{2} + 2k) = k^{3} + 2k^{2} - k - 1$, $\Delta(k^{3} + 2k^{2} - k - 1) = 1$.

By Rule~2 applied to $(k^{3} + k - 1) + (k-1)\cdot 2k = k^{3} + 2k^{2} - k - 1$ and $(2k - 1) + (k-1)\cdot (k^{2} + k) = k^{3} + k - 1$, $\Delta(k^{3} + k - 1) = 3$.

By Rule~2 applied to $(3k - 2) + (k-1)\cdot (2k - 1) = 2k^{2} - 1$ and $(3k - 2) + (k-1)\cdot (k^{2} + k - 1) = k^{3} + k - 1$, $\Delta(3k - 2) = 1$.

By Rule~2 applied to $(3k - 2) + (k-1)\cdot 2k = 2k^{2} + k - 2$ and $(k^{2} + k - 1) + (k-1)\cdot (k + 1) = 2k^{2} + k - 2$, $\Delta(2k^{2} + k - 2) = 2$.

By Rule~2 applied to $(k^{3} + 2k^{2} - k - 1) + (k-1)\cdot 1 = k^{3} + 2k^{2} - 2$ and $(2k^{2} + k - 2) + (k-1)\cdot (k^{2} + k) = k^{3} + 2k^{2} - 2$, $\Delta(k^{3} + 2k^{2} - 2) = 3$.

\emph{Contradiction (Rule~3) at $k^{2} + 2k - 2$:} if $\Delta(k^{2} + 2k - 2)=1$ then $(k^{2} + 2k - 2) + (k-1)\cdot (k^{2} + 2k - 1) = k^{3} + 2k^{2} - k - 1$ is monochromatic; if $\Delta(k^{2} + 2k - 2)=2$ then $(k^{2} + 2k - 2) + (k-1)\cdot k = 2k^{2} + k - 2$ is monochromatic; if $\Delta(k^{2} + 2k - 2)=3$ then $(k^{2} + 2k - 2) + (k-1)\cdot (k^{2} + 2k) = k^{3} + 2k^{2} - 2$ is monochromatic.

\case{1.2.2.2.2.2.2}{$\Delta(k^{3} + k^{2} - 2k + 1) = 2$}\quad {\small($A_1\supseteq \{1, 2k, k^{2} - 2, k^{2} - 1, k^{2}, k^{2} + 2k - 1\},\ A_2\supseteq \{k, 2k - 1, k^{2} + k, k^{3} - k^{2} + 1, k^{3} + k^{2} - 2k + 1, k^{3} + k^{2} - k\},\ A_3\supseteq \{k + 1, k^{2} + k - 1, k^{3} - k + 1, k^{3}, k^{3} + 2k^{2} - 2k\}$)}

By Rule~2 applied to $k^{2} + (k-1)\cdot (k^{2} + 2k - 1) = k^{3} + 2k^{2} - 3k + 1$ and $(k^{3} + k^{2} - 2k + 1) + (k-1)\cdot k = k^{3} + 2k^{2} - 3k + 1$, $\Delta(k^{3} + 2k^{2} - 3k + 1) = 3$.

By Rule~2 applied to $(2k - 1) + (k-1)\cdot (k^{2} + 2k - 2) = k^{3} + k^{2} - 2k + 1$ and $(k^{2} + k - 1) + (k-1)\cdot (k^{2} + 2k - 2) = k^{3} + 2k^{2} - 3k + 1$, $\Delta(k^{2} + 2k - 2) = 1$.

By Rule~2 applied to $(k^{2} + 2k - 2) + (k-1)\cdot (k^{2} - 1) = k^{3} + k - 1$ and $(2k - 1) + (k-1)\cdot (k^{2} + k) = k^{3} + k - 1$, $\Delta(k^{3} + k - 1) = 3$.

By Rule~2 applied to $k^{2} + (k-1)\cdot 2 = k^{2} + 2k - 2$ and $(k^{3} - k + 1) + (k-1)\cdot 2 = k^{3} + k - 1$, $\Delta(2) = 2$.

By Rule~2 applied to $(k^{3} + k^{2} - 2k + 1) + (k-1)\cdot 2 = k^{3} + k^{2} - 1$ and $k^{3} + (k-1)\cdot (k + 1) = k^{3} + k^{2} - 1$, $\Delta(k^{3} + k^{2} - 1) = 1$.

\emph{Contradiction (Rule~3) at $3k - 2$:} if $\Delta(3k - 2)=1$ then $(3k - 2) + (k-1)\cdot (k^{2} + 2k - 1) = k^{3} + k^{2} - 1$ is monochromatic; if $\Delta(3k - 2)=2$ then $k + (k-1)\cdot 2 = 3k - 2$ is monochromatic; if $\Delta(3k - 2)=3$ then $(3k - 2) + (k-1)\cdot (k^{2} + k - 1) = k^{3} + k - 1$ is monochromatic.

\case{2}{$\Delta(k^{2}) = 3$}\quad {\small($A_1\supseteq \{1\},\ A_2\supseteq \{k\},\ A_3\supseteq \{k^{2}\}$)}

Since $k\cdot k^{2} = k^{3}$ must not be monochromatic in $A_3$, either $\Delta(k^{3}) = 1$ or $\Delta(k^{3}) = 2$.

\case{2.1}{$\Delta(k^{3}) = 1$}\quad {\small($A_1\supseteq \{1, k^{3}\},\ A_2\supseteq \{k\},\ A_3\supseteq \{k^{2}\}$)}

Since $k^{3} + (k-1)\cdot 1 = k^{3} + k - 1$ must not be monochromatic in $A_1$, either $\Delta(k^{3} + k - 1) = 2$ or $\Delta(k^{3} + k - 1) = 3$.

\case{2.1.1}{$\Delta(k^{3} + k - 1) = 2$}\quad {\small($A_1\supseteq \{1, k^{3}\},\ A_2\supseteq \{k, k^{3} + k - 1\},\ A_3\supseteq \{k^{2}\}$)}

By Rule~2 applied to $1 + (k-1)\cdot (k^{2} + k + 1) = k^{3}$ and $k + (k-1)\cdot (k^{2} + k + 1) = k^{3} + k - 1$, $\Delta(k^{2} + k + 1) = 3$.

By Rule~2 applied to $(k^{3} + k - 1) + (k-1)\cdot k = k^{3} + k^{2} - 1$ and $k^{2} + (k-1)\cdot (k^{2} + k + 1) = k^{3} + k^{2} - 1$, $\Delta(k^{3} + k^{2} - 1) = 1$.

Since $(k^{3} - k + 1) + (k-1)\cdot 1 = k^{3}$ must not be monochromatic in $A_1$, either $\Delta(k^{3} - k + 1) = 2$ or $\Delta(k^{3} - k + 1) = 3$.

\case{2.1.1.1}{$\Delta(k^{3} - k + 1) = 2$}\quad {\small($A_1\supseteq \{1, k^{3}, k^{3} + k^{2} - 1\},\ A_2\supseteq \{k, k^{3} - k + 1, k^{3} + k - 1\},\ A_3\supseteq \{k^{2}, k^{2} + k + 1\}$)}

Since $k^{3} + (k-1)\cdot (k + 1) = k^{3} + k^{2} - 1$ must not be monochromatic in $A_1$, either $\Delta(k + 1) = 2$ or $\Delta(k + 1) = 3$.

\case{2.1.1.1.1}{$\Delta(k + 1) = 2$}\quad {\small($A_1\supseteq \{1, k^{3}, k^{3} + k^{2} - 1\},\ A_2\supseteq \{k, k + 1, k^{3} - k + 1, k^{3} + k - 1\},\ A_3\supseteq \{k^{2}, k^{2} + k + 1\}$)}

By Rule~2 applied to $(k^{3} + k^{2} - k) + (k-1)\cdot 1 = k^{3} + k^{2} - 1$ and $(k^{3} - k + 1) + (k-1)\cdot (k + 1) = k^{3} + k^{2} - k$, $\Delta(k^{3} + k^{2} - k) = 3$.

By Rule~2 applied to $k\cdot (k + 1) = k^{2} + k$ and $k^{2} + (k-1)\cdot (k^{2} + k) = k^{3} + k^{2} - k$, $\Delta(k^{2} + k) = 1$.

By Rule~2 applied to $(k^{2} + k - 1) + (k-1)\cdot (k^{2} + k) = k^{3} + k^{2} - 1$ and $k + (k-1)\cdot (k + 1) = k^{2} + k - 1$, $\Delta(k^{2} + k - 1) = 3$.

\emph{Contradiction:} $k\cdot (k^{2} + k - 1) = k^{3} + k^{2} - k$ is now monochromatic.

\case{2.1.1.1.2}{$\Delta(k + 1) = 3$}\quad {\small($A_1\supseteq \{1, k^{3}, k^{3} + k^{2} - 1\},\ A_2\supseteq \{k, k^{3} - k + 1, k^{3} + k - 1\},\ A_3\supseteq \{k + 1, k^{2}, k^{2} + k + 1\}$)}

Since $k + (k-1)\cdot (k^{2} + k - 1) = k^{3} - k + 1$ must not be monochromatic in $A_2$, either $\Delta(k^{2} + k - 1) = 1$ or $\Delta(k^{2} + k - 1) = 3$.

\case{2.1.1.1.2.1}{$\Delta(k^{2} + k - 1) = 1$}\quad {\small($A_1\supseteq \{1, k^{2} + k - 1, k^{3}, k^{3} + k^{2} - 1\},\ A_2\supseteq \{k, k^{3} - k + 1, k^{3} + k - 1\},\ A_3\supseteq \{k + 1, k^{2}, k^{2} + k + 1\}$)}

By Rule~2 applied to $1 + (k-1)\cdot (k + 2) = k^{2} + k - 1$ and $(k + 2) + (k-1)\cdot (k + 1) = k^{2} + k + 1$, $\Delta(k + 2) = 2$.

By Rule~2 applied to $(k^{2} + k - 1) + (k-1)\cdot 1 = k^{2} + 2k - 2$ and $k + (k-1)\cdot (k + 2) = k^{2} + 2k - 2$, $\Delta(k^{2} + 2k - 2) = 3$.

By Rule~2 applied to $(2k - 1) + (k-1)\cdot (k^{2} + k - 1) = k^{3}$ and $(2k - 1) + (k-1)\cdot (k + 1) = k^{2} + 2k - 2$, $\Delta(2k - 1) = 2$.

\emph{Contradiction (Rule~3) at $k^{2} + k$:} if $\Delta(k^{2} + k)=1$ then $(k^{2} + k - 1) + (k-1)\cdot (k^{2} + k) = k^{3} + k^{2} - 1$ is monochromatic; if $\Delta(k^{2} + k)=2$ then $(2k - 1) + (k-1)\cdot (k^{2} + k) = k^{3} + k - 1$ is monochromatic; if $\Delta(k^{2} + k)=3$ then $k\cdot (k + 1) = k^{2} + k$ is monochromatic.

\case{2.1.1.1.2.2}{$\Delta(k^{2} + k - 1) = 3$}\quad {\small($A_1\supseteq \{1, k^{3}, k^{3} + k^{2} - 1\},\ A_2\supseteq \{k, k^{3} - k + 1, k^{3} + k - 1\},\ A_3\supseteq \{k + 1, k^{2}, k^{2} + k - 1, k^{2} + k + 1\}$)}

By Rule~2 applied to $(k^{3} - k + 1) + (k-1)\cdot k = k^{3} + k^{2} - 2k + 1$ and $k^{2} + (k-1)\cdot (k^{2} + k - 1) = k^{3} + k^{2} - 2k + 1$, $\Delta(k^{3} + k^{2} - 2k + 1) = 1$.

By Rule~2 applied to $(k^{3} + k^{2} - k) + (k-1)\cdot 1 = k^{3} + k^{2} - 1$ and $k\cdot (k^{2} + k - 1) = k^{3} + k^{2} - k$, $\Delta(k^{3} + k^{2} - k) = 2$.

By Rule~2 applied to $(k^{3} + k^{2} - 2k + 1) + (k-1)\cdot 2 = k^{3} + k^{2} - 1$ and $(k^{3} - k + 1) + (k-1)\cdot 2 = k^{3} + k - 1$, $\Delta(2) = 3$.

By Rule~2 applied to $1 + (k-1)\cdot (k^{2} + 2k) = k^{3} + k^{2} - 2k + 1$ and $k + (k-1)\cdot (k^{2} + 2k) = k^{3} + k^{2} - k$, $\Delta(k^{2} + 2k) = 3$.

By Rule~2 applied to $(k^{3} + k^{2} - k) + (k-1)\cdot k = k^{3} + 2k^{2} - 2k$ and $k^{2} + (k-1)\cdot (k^{2} + 2k) = k^{3} + 2k^{2} - 2k$, $\Delta(k^{3} + 2k^{2} - 2k) = 1$.

By Rule~2 applied to $k^{3} + (k-1)\cdot 2k = k^{3} + 2k^{2} - 2k$ and $k\cdot 2 = 2k$, $\Delta(2k) = 2$.

\emph{Contradiction (Rule~3) at $k^{3} + 2k^{2} - k - 1$:} if $\Delta(k^{3} + 2k^{2} - k - 1)=1$ then $(k^{3} + 2k^{2} - 2k) + (k-1)\cdot 1 = k^{3} + 2k^{2} - k - 1$ is monochromatic; if $\Delta(k^{3} + 2k^{2} - k - 1)=2$ then $(k^{3} + k - 1) + (k-1)\cdot 2k = k^{3} + 2k^{2} - k - 1$ is monochromatic; if $\Delta(k^{3} + 2k^{2} - k - 1)=3$ then $(k^{2} + k - 1) + (k-1)\cdot (k^{2} + 2k) = k^{3} + 2k^{2} - k - 1$ is monochromatic.

\case{2.1.1.2}{$\Delta(k^{3} - k + 1) = 3$}\quad {\small($A_1\supseteq \{1, k^{3}, k^{3} + k^{2} - 1\},\ A_2\supseteq \{k, k^{3} + k - 1\},\ A_3\supseteq \{k^{2}, k^{2} + k + 1, k^{3} - k + 1\}$)}

Since $(k^{3} + k^{2} - k) + (k-1)\cdot 1 = k^{3} + k^{2} - 1$ must not be monochromatic in $A_1$, either $\Delta(k^{3} + k^{2} - k) = 2$ or $\Delta(k^{3} + k^{2} - k) = 3$.

\case{2.1.1.2.1}{$\Delta(k^{3} + k^{2} - k) = 2$}\quad {\small($A_1\supseteq \{1, k^{3}, k^{3} + k^{2} - 1\},\ A_2\supseteq \{k, k^{3} + k - 1, k^{3} + k^{2} - k\},\ A_3\supseteq \{k^{2}, k^{2} + k + 1, k^{3} - k + 1\}$)}

Since $k^{2} + (k-1)\cdot (k^{2} - 1) = k^{3} - k + 1$ must not be monochromatic in $A_3$, either $\Delta(k^{2} - 1) = 1$ or $\Delta(k^{2} - 1) = 2$.

\case{2.1.1.2.1.1}{$\Delta(k^{2} - 1) = 1$}\quad {\small($A_1\supseteq \{1, k^{2} - 1, k^{3}, k^{3} + k^{2} - 1\},\ A_2\supseteq \{k, k^{3} + k - 1, k^{3} + k^{2} - k\},\ A_3\supseteq \{k^{2}, k^{2} + k + 1, k^{3} - k + 1\}$)}

Since $(k^{2} + k - 1) + (k-1)\cdot (k^{2} - 1) = k^{3}$ must not be monochromatic in $A_1$, either $\Delta(k^{2} + k - 1) = 2$ or $\Delta(k^{2} + k - 1) = 3$.

\case{2.1.1.2.1.1.1}{$\Delta(k^{2} + k - 1) = 2$}\quad {\small($A_1\supseteq \{1, k^{2} - 1, k^{3}, k^{3} + k^{2} - 1\},\ A_2\supseteq \{k, k^{2} + k - 1, k^{3} + k - 1, k^{3} + k^{2} - k\},\ A_3\supseteq \{k^{2}, k^{2} + k + 1, k^{3} - k + 1\}$)}

\emph{Contradiction:} $k\cdot (k^{2} + k - 1) = k^{3} + k^{2} - k$ is now monochromatic.

\case{2.1.1.2.1.1.2}{$\Delta(k^{2} + k - 1) = 3$}\quad {\small($A_1\supseteq \{1, k^{2} - 1, k^{3}, k^{3} + k^{2} - 1\},\ A_2\supseteq \{k, k^{3} + k - 1, k^{3} + k^{2} - k\},\ A_3\supseteq \{k^{2}, k^{2} + k - 1, k^{2} + k + 1, k^{3} - k + 1\}$)}

Since $(2k^{2} + k - 2) + (k-1)\cdot (k^{2} - 1) = k^{3} + k^{2} - 1$ must not be monochromatic in $A_1$, either $\Delta(2k^{2} + k - 2) = 2$ or $\Delta(2k^{2} + k - 2) = 3$.

\case{2.1.1.2.1.1.2.1}{$\Delta(2k^{2} + k - 2) = 2$}\quad {\small($A_1\supseteq \{1, k^{2} - 1, k^{3}, k^{3} + k^{2} - 1\},\ A_2\supseteq \{k, 2k^{2} + k - 2, k^{3} + k - 1, k^{3} + k^{2} - k\},\ A_3\supseteq \{k^{2}, k^{2} + k - 1, k^{2} + k + 1, k^{3} - k + 1\}$)}

By Rule~2 applied to $(2k^{2} + k - 2) + (k-1)\cdot (k^{2} - 2) = k^{3} + k^{2} - k$ and $(k^{2} + k - 1) + (k-1)\cdot (k^{2} - 2) = k^{3} - k + 1$, $\Delta(k^{2} - 2) = 1$.

By Rule~2 applied to $(k^{2} + 2k - 2) + (k-1)\cdot (k^{2} - 2) = k^{3}$ and $(k^{2} + 2k - 2) + (k-1)\cdot k = 2k^{2} + k - 2$, $\Delta(k^{2} + 2k - 2) = 3$.

By Rule~2 applied to $(k^{3} + k^{2} - k) + (k-1)\cdot k = k^{3} + 2k^{2} - 2k$ and $k\cdot (k^{2} + 2k - 2) = k^{3} + 2k^{2} - 2k$, $\Delta(k^{3} + 2k^{2} - 2k) = 1$.

By Rule~2 applied to $(k^{3} + 2k^{2} - 3k + 1) + (k-1)\cdot 1 = k^{3} + 2k^{2} - 2k$ and $(k^{2} + k - 1) + (k-1)\cdot (k^{2} + 2k - 2) = k^{3} + 2k^{2} - 3k + 1$, $\Delta(k^{3} + 2k^{2} - 3k + 1) = 2$.

By Rule~2 applied to $(k^{3} + k^{2} - 2k + 1) + (k-1)\cdot k = k^{3} + 2k^{2} - 3k + 1$ and $k^{2} + (k-1)\cdot (k^{2} + k - 1) = k^{3} + k^{2} - 2k + 1$, $\Delta(k^{3} + k^{2} - 2k + 1) = 1$.

By Rule~2 applied to $1 + (k-1)\cdot (k^{2} + 2k) = k^{3} + k^{2} - 2k + 1$ and $k + (k-1)\cdot (k^{2} + 2k) = k^{3} + k^{2} - k$, $\Delta(k^{2} + 2k) = 3$.

By Rule~2 applied to $(k^{3} + k^{2} - 2k + 1) + (k-1)\cdot 2 = k^{3} + k^{2} - 1$ and $k^{2} + (k-1)\cdot 2 = k^{2} + 2k - 2$, $\Delta(2) = 2$.

\emph{Contradiction (Rule~3) at $k^{3} + 2k^{2} - k - 1$:} if $\Delta(k^{3} + 2k^{2} - k - 1)=1$ then $(k^{3} + 2k^{2} - 2k) + (k-1)\cdot 1 = k^{3} + 2k^{2} - k - 1$ is monochromatic; if $\Delta(k^{3} + 2k^{2} - k - 1)=2$ then $(k^{3} + 2k^{2} - 3k + 1) + (k-1)\cdot 2 = k^{3} + 2k^{2} - k - 1$ is monochromatic; if $\Delta(k^{3} + 2k^{2} - k - 1)=3$ then $(k^{2} + k - 1) + (k-1)\cdot (k^{2} + 2k) = k^{3} + 2k^{2} - k - 1$ is monochromatic.

\case{2.1.1.2.1.1.2.2}{$\Delta(2k^{2} + k - 2) = 3$}\quad {\small($A_1\supseteq \{1, k^{2} - 1, k^{3}, k^{3} + k^{2} - 1\},\ A_2\supseteq \{k, k^{3} + k - 1, k^{3} + k^{2} - k\},\ A_3\supseteq \{k^{2}, k^{2} + k - 1, k^{2} + k + 1, 2k^{2} + k - 2, k^{3} - k + 1\}$)}

By Rule~2 applied to $k^{3} + (k-1)\cdot (k + 1) = k^{3} + k^{2} - 1$ and $(k^{2} + k - 1) + (k-1)\cdot (k + 1) = 2k^{2} + k - 2$, $\Delta(k + 1) = 2$.

By Rule~2 applied to $(k^{3} + k^{2} - k) + (k-1)\cdot (k + 1) = k^{3} + 2k^{2} - k - 1$ and $(2k^{2} + k - 2) + (k-1)\cdot (k^{2} + k - 1) = k^{3} + 2k^{2} - k - 1$, $\Delta(k^{3} + 2k^{2} - k - 1) = 1$.

By Rule~2 applied to $(k^{3} + 2k^{2} - 2k) + (k-1)\cdot 1 = k^{3} + 2k^{2} - k - 1$ and $(k^{3} + k^{2} - k) + (k-1)\cdot k = k^{3} + 2k^{2} - 2k$, $\Delta(k^{3} + 2k^{2} - 2k) = 3$.

By Rule~2 applied to $k + (k-1)\cdot (k^{2} + 2k) = k^{3} + k^{2} - k$ and $k^{2} + (k-1)\cdot (k^{2} + 2k) = k^{3} + 2k^{2} - 2k$, $\Delta(k^{2} + 2k) = 1$.

By Rule~2 applied to $1 + (k-1)\cdot (k^{2} + 2k) = k^{3} + k^{2} - 2k + 1$ and $k^{2} + (k-1)\cdot (k^{2} + k - 1) = k^{3} + k^{2} - 2k + 1$, $\Delta(k^{3} + k^{2} - 2k + 1) = 2$.

By Rule~2 applied to $k + (k-1)\cdot (k^{2} + 2k - 1) = k^{3} + k^{2} - 2k + 1$ and $(k^{2} + k - 1) + (k-1)\cdot (k^{2} + 2k - 1) = k^{3} + 2k^{2} - 2k$, $\Delta(k^{2} + 2k - 1) = 1$.

By Rule~2 applied to $(k^{2} + k) + (k-1)\cdot 1 = k^{2} + 2k - 1$ and $k\cdot (k + 1) = k^{2} + k$, $\Delta(k^{2} + k) = 3$.

By Rule~2 applied to $(k^{3} + 2k^{2} - k - 1) + (k-1)\cdot 1 = k^{3} + 2k^{2} - 2$ and $(2k^{2} + k - 2) + (k-1)\cdot (k^{2} + k) = k^{3} + 2k^{2} - 2$, $\Delta(k^{3} + 2k^{2} - 2) = 2$.

By Rule~2 applied to $(k^{3} + k^{2} - k) + (k-1)\cdot (k + 2) = k^{3} + 2k^{2} - 2$ and $k^{2} + (k-1)\cdot (k + 2) = 2k^{2} + k - 2$, $\Delta(k + 2) = 1$.

\emph{Contradiction:} $k\cdot (k + 2) = k^{2} + 2k$ is now monochromatic.

\case{2.1.1.2.1.2}{$\Delta(k^{2} - 1) = 2$}\quad {\small($A_1\supseteq \{1, k^{3}, k^{3} + k^{2} - 1\},\ A_2\supseteq \{k, k^{2} - 1, k^{3} + k - 1, k^{3} + k^{2} - k\},\ A_3\supseteq \{k^{2}, k^{2} + k + 1, k^{3} - k + 1\}$)}

Since $(k^{2} + 2k - 2) + (k-1)\cdot (k^{2} - 1) = k^{3} + k - 1$ must not be monochromatic in $A_2$, either $\Delta(k^{2} + 2k - 2) = 1$ or $\Delta(k^{2} + 2k - 2) = 3$.

\case{2.1.1.2.1.2.1}{$\Delta(k^{2} + 2k - 2) = 1$}\quad {\small($A_1\supseteq \{1, k^{2} + 2k - 2, k^{3}, k^{3} + k^{2} - 1\},\ A_2\supseteq \{k, k^{2} - 1, k^{3} + k - 1, k^{3} + k^{2} - k\},\ A_3\supseteq \{k^{2}, k^{2} + k + 1, k^{3} - k + 1\}$)}

By Rule~2 applied to $k\cdot (k^{2} + 2k - 2) = k^{3} + 2k^{2} - 2k$ and $(k^{3} + k^{2} - k) + (k-1)\cdot k = k^{3} + 2k^{2} - 2k$, $\Delta(k^{3} + 2k^{2} - 2k) = 3$.

By Rule~2 applied to $k + (k-1)\cdot (k^{2} + 2k) = k^{3} + k^{2} - k$ and $k^{2} + (k-1)\cdot (k^{2} + 2k) = k^{3} + 2k^{2} - 2k$, $\Delta(k^{2} + 2k) = 1$.

Since $(k - 1) + (k-1)\cdot k = k^{2} - 1$ must not be monochromatic in $A_2$, either $\Delta(k - 1) = 1$ or $\Delta(k - 1) = 3$.

\case{2.1.1.2.1.2.1.1}{$\Delta(k - 1) = 1$}\quad {\small($A_1\supseteq \{1, k - 1, k^{2} + 2k - 2, k^{2} + 2k, k^{3}, k^{3} + k^{2} - 1\},\ A_2\supseteq \{k, k^{2} - 1, k^{3} + k - 1, k^{3} + k^{2} - k\},\ A_3\supseteq \{k^{2}, k^{2} + k + 1, k^{3} - k + 1, k^{3} + 2k^{2} - 2k\}$)}

By Rule~2 applied to $(k^{2} + 2k - 2) + (k-1)\cdot (k - 1) = 2k^{2} - 1$ and $(2k^{2} - 1) + (k-1)\cdot (k^{2} - 1) = k^{3} + k^{2} - k$, $\Delta(2k^{2} - 1) = 3$.

By Rule~2 applied to $(k^{2} + k - 1) + (k-1)\cdot 1 = k^{2} + 2k - 2$ and $(2k^{2} - 1) + (k-1)\cdot (k^{2} + k - 1) = k^{3} + 2k^{2} - 2k$, $\Delta(k^{2} + k - 1) = 2$.

\emph{Contradiction:} $k\cdot (k^{2} + k - 1) = k^{3} + k^{2} - k$ is now monochromatic.

\case{2.1.1.2.1.2.1.2}{$\Delta(k - 1) = 3$}\quad {\small($A_1\supseteq \{1, k^{2} + 2k - 2, k^{2} + 2k, k^{3}, k^{3} + k^{2} - 1\},\ A_2\supseteq \{k, k^{2} - 1, k^{3} + k - 1, k^{3} + k^{2} - k\},\ A_3\supseteq \{k - 1, k^{2}, k^{2} + k + 1, k^{3} - k + 1, k^{3} + 2k^{2} - 2k\}$)}

By Rule~2 applied to $(2k - 1) + (k-1)\cdot (k^{2} + 2k) = k^{3} + k^{2} - 1$ and $(2k - 1) + (k-1)\cdot (k - 1) = k^{2}$, $\Delta(2k - 1) = 2$.

By Rule~2 applied to $(k^{2} + k - 1) + (k-1)\cdot 1 = k^{2} + 2k - 2$ and $(2k - 1) + (k-1)\cdot k = k^{2} + k - 1$, $\Delta(k^{2} + k - 1) = 3$.

By Rule~2 applied to $(k^{2} + 2k - 2) + (k-1)\cdot (k^{2} - 2) = k^{3}$ and $(k^{2} + k - 1) + (k-1)\cdot (k^{2} - 2) = k^{3} - k + 1$, $\Delta(k^{2} - 2) = 2$.

By Rule~2 applied to $(2k^{2} - 1) + (k-1)\cdot (k^{2} - 1) = k^{3} + k^{2} - k$ and $(2k^{2} - 1) + (k-1)\cdot (k^{2} + k - 1) = k^{3} + 2k^{2} - 2k$, $\Delta(2k^{2} - 1) = 1$.

By Rule~2 applied to $(2k - 1) + (k-1)\cdot (k^{2} + 2k - 1) = k^{3} + k^{2} - k$ and $(k^{2} + k - 1) + (k-1)\cdot (k^{2} + 2k - 1) = k^{3} + 2k^{2} - 2k$, $\Delta(k^{2} + 2k - 1) = 1$.

By Rule~2 applied to $(2k^{2} - 1) + (k-1)\cdot 1 = 2k^{2} + k - 2$ and $(2k^{2} + k - 2) + (k-1)\cdot (k^{2} - 2) = k^{3} + k^{2} - k$, $\Delta(2k^{2} + k - 2) = 3$.

By Rule~2 applied to $k^{3} + (k-1)\cdot (k + 1) = k^{3} + k^{2} - 1$ and $(k^{2} + k - 1) + (k-1)\cdot (k + 1) = 2k^{2} + k - 2$, $\Delta(k + 1) = 2$.

\emph{Contradiction (Rule~3) at $k^{3} + 2k^{2} - k - 1$:} if $\Delta(k^{3} + 2k^{2} - k - 1)=1$ then $(k^{2} + 2k - 2) + (k-1)\cdot (k^{2} + 2k - 1) = k^{3} + 2k^{2} - k - 1$ is monochromatic; if $\Delta(k^{3} + 2k^{2} - k - 1)=2$ then $(k^{3} + k^{2} - k) + (k-1)\cdot (k + 1) = k^{3} + 2k^{2} - k - 1$ is monochromatic; if $\Delta(k^{3} + 2k^{2} - k - 1)=3$ then $(2k^{2} + k - 2) + (k-1)\cdot (k^{2} + k - 1) = k^{3} + 2k^{2} - k - 1$ is monochromatic.

\case{2.1.1.2.1.2.2}{$\Delta(k^{2} + 2k - 2) = 3$}\quad {\small($A_1\supseteq \{1, k^{3}, k^{3} + k^{2} - 1\},\ A_2\supseteq \{k, k^{2} - 1, k^{3} + k - 1, k^{3} + k^{2} - k\},\ A_3\supseteq \{k^{2}, k^{2} + k + 1, k^{2} + 2k - 2, k^{3} - k + 1\}$)}

By Rule~2 applied to $(k^{3} + k^{2} - k) + (k-1)\cdot k = k^{3} + 2k^{2} - 2k$ and $k\cdot (k^{2} + 2k - 2) = k^{3} + 2k^{2} - 2k$, $\Delta(k^{3} + 2k^{2} - 2k) = 1$.

By Rule~2 applied to $(k^{3} + k^{2} - 1) + (k-1)\cdot (k - 1) = k^{3} + 2k^{2} - 2k$ and $(k - 1) + (k-1)\cdot k = k^{2} - 1$, $\Delta(k - 1) = 3$.

By Rule~2 applied to $(2k^{2} - 1) + (k-1)\cdot (k^{2} - 1) = k^{3} + k^{2} - k$ and $(k^{2} + 2k - 2) + (k-1)\cdot (k - 1) = 2k^{2} - 1$, $\Delta(2k^{2} - 1) = 1$.

Since $(2k - 1) + (k-1)\cdot (k - 1) = k^{2}$ must not be monochromatic in $A_3$, either $\Delta(2k - 1) = 1$ or $\Delta(2k - 1) = 2$.

\case{2.1.1.2.1.2.2.1}{$\Delta(2k - 1) = 1$}\quad {\small($A_1\supseteq \{1, 2k - 1, 2k^{2} - 1, k^{3}, k^{3} + k^{2} - 1, k^{3} + 2k^{2} - 2k\},\ A_2\supseteq \{k, k^{2} - 1, k^{3} + k - 1, k^{3} + k^{2} - k\},\ A_3\supseteq \{k - 1, k^{2}, k^{2} + k + 1, k^{2} + 2k - 2, k^{3} - k + 1\}$)}

By Rule~2 applied to $1 + (k-1)\cdot 2 = 2k - 1$ and $k^{2} + (k-1)\cdot 2 = k^{2} + 2k - 2$, $\Delta(2) = 2$.

By Rule~2 applied to $(2k - 1) + (k-1)\cdot 1 = 3k - 2$ and $k + (k-1)\cdot 2 = 3k - 2$, $\Delta(3k - 2) = 3$.

By Rule~2 applied to $(2k - 1) + (k-1)\cdot (k^{2} + k - 1) = k^{3}$ and $(3k - 2) + (k-1)\cdot (k - 1) = k^{2} + k - 1$, $\Delta(k^{2} + k - 1) = 2$.

\emph{Contradiction:} $k\cdot (k^{2} + k - 1) = k^{3} + k^{2} - k$ is now monochromatic.

\case{2.1.1.2.1.2.2.2}{$\Delta(2k - 1) = 2$}\quad {\small($A_1\supseteq \{1, 2k^{2} - 1, k^{3}, k^{3} + k^{2} - 1, k^{3} + 2k^{2} - 2k\},\ A_2\supseteq \{k, 2k - 1, k^{2} - 1, k^{3} + k - 1, k^{3} + k^{2} - k\},\ A_3\supseteq \{k - 1, k^{2}, k^{2} + k + 1, k^{2} + 2k - 2, k^{3} - k + 1\}$)}

By Rule~2 applied to $(2k^{2} - 1) + (k-1)\cdot (k^{2} + k - 1) = k^{3} + 2k^{2} - 2k$ and $(2k - 1) + (k-1)\cdot k = k^{2} + k - 1$, $\Delta(k^{2} + k - 1) = 3$.

By Rule~2 applied to $k + (k-1)\cdot (2k - 1) = 2k^{2} - 2k + 1$ and $k^{2} + (k-1)\cdot (k - 1) = 2k^{2} - 2k + 1$, $\Delta(2k^{2} - 2k + 1) = 1$.

By Rule~2 applied to $(k^{3} + 2k^{2} - 3k + 1) + (k-1)\cdot 1 = k^{3} + 2k^{2} - 2k$ and $(k^{2} + k - 1) + (k-1)\cdot (k^{2} + 2k - 2) = k^{3} + 2k^{2} - 3k + 1$, $\Delta(k^{3} + 2k^{2} - 3k + 1) = 2$.

By Rule~2 applied to $(2k^{2} - 2k + 1) + (k-1)\cdot 2 = 2k^{2} - 1$ and $k^{2} + (k-1)\cdot 2 = k^{2} + 2k - 2$, $\Delta(2) = 2$.

By Rule~2 applied to $(k^{3} + 2k^{2} - 2k) + (k-1)\cdot 1 = k^{3} + 2k^{2} - k - 1$ and $(k^{3} + 2k^{2} - 3k + 1) + (k-1)\cdot 2 = k^{3} + 2k^{2} - k - 1$, $\Delta(k^{3} + 2k^{2} - k - 1) = 3$.

By Rule~2 applied to $(2k - 1) + (k-1)\cdot (k^{2} + 2k - 1) = k^{3} + k^{2} - k$ and $(k^{2} + 2k - 2) + (k-1)\cdot (k^{2} + 2k - 1) = k^{3} + 2k^{2} - k - 1$, $\Delta(k^{2} + 2k - 1) = 1$.

\emph{Contradiction (Rule~3) at $3k - 2$:} if $\Delta(3k - 2)=1$ then $(3k - 2) + (k-1)\cdot (k^{2} + 2k - 1) = k^{3} + k^{2} - 1$ is monochromatic; if $\Delta(3k - 2)=2$ then $k + (k-1)\cdot 2 = 3k - 2$ is monochromatic; if $\Delta(3k - 2)=3$ then $(3k - 2) + (k-1)\cdot (k - 1) = k^{2} + k - 1$ is monochromatic.

\case{2.1.1.2.2}{$\Delta(k^{3} + k^{2} - k) = 3$}\quad {\small($A_1\supseteq \{1, k^{3}, k^{3} + k^{2} - 1\},\ A_2\supseteq \{k, k^{3} + k - 1\},\ A_3\supseteq \{k^{2}, k^{2} + k + 1, k^{3} - k + 1, k^{3} + k^{2} - k\}$)}

By Rule~2 applied to $k^{3} + (k-1)\cdot (k + 1) = k^{3} + k^{2} - 1$ and $(k^{3} - k + 1) + (k-1)\cdot (k + 1) = k^{3} + k^{2} - k$, $\Delta(k + 1) = 2$.

By Rule~2 applied to $k\cdot (k + 1) = k^{2} + k$ and $k^{2} + (k-1)\cdot (k^{2} + k) = k^{3} + k^{2} - k$, $\Delta(k^{2} + k) = 1$.

By Rule~2 applied to $(k^{2} + k - 1) + (k-1)\cdot (k^{2} + k) = k^{3} + k^{2} - 1$ and $k + (k-1)\cdot (k + 1) = k^{2} + k - 1$, $\Delta(k^{2} + k - 1) = 3$.

\emph{Contradiction:} $k\cdot (k^{2} + k - 1) = k^{3} + k^{2} - k$ is now monochromatic.

\case{2.1.2}{$\Delta(k^{3} + k - 1) = 3$}\quad {\small($A_1\supseteq \{1, k^{3}\},\ A_2\supseteq \{k\},\ A_3\supseteq \{k^{2}, k^{3} + k - 1\}$)}

Since $(k^{2} + k - 1) + (k-1)\cdot k^{2} = k^{3} + k - 1$ must not be monochromatic in $A_3$, either $\Delta(k^{2} + k - 1) = 1$ or $\Delta(k^{2} + k - 1) = 2$.

\case{2.1.2.1}{$\Delta(k^{2} + k - 1) = 1$}\quad {\small($A_1\supseteq \{1, k^{2} + k - 1, k^{3}\},\ A_2\supseteq \{k\},\ A_3\supseteq \{k^{2}, k^{3} + k - 1\}$)}

Since $(k^{2} + k - 1) + (k-1)\cdot 1 = k^{2} + 2k - 2$ must not be monochromatic in $A_1$, either $\Delta(k^{2} + 2k - 2) = 2$ or $\Delta(k^{2} + 2k - 2) = 3$.

\case{2.1.2.1.1}{$\Delta(k^{2} + 2k - 2) = 2$}\quad {\small($A_1\supseteq \{1, k^{2} + k - 1, k^{3}\},\ A_2\supseteq \{k, k^{2} + 2k - 2\},\ A_3\supseteq \{k^{2}, k^{3} + k - 1\}$)}

By Rule~2 applied to $1 + (k-1)\cdot (k + 2) = k^{2} + k - 1$ and $k + (k-1)\cdot (k + 2) = k^{2} + 2k - 2$, $\Delta(k + 2) = 3$.

By Rule~2 applied to $(k^{2} + 2k - 2) + (k-1)\cdot k = 2k^{2} + k - 2$ and $k^{2} + (k-1)\cdot (k + 2) = 2k^{2} + k - 2$, $\Delta(2k^{2} + k - 2) = 1$.

Since $k\cdot (k^{2} + 2k - 2) = k^{3} + 2k^{2} - 2k$ must not be monochromatic in $A_2$, either $\Delta(k^{3} + 2k^{2} - 2k) = 1$ or $\Delta(k^{3} + 2k^{2} - 2k) = 3$.

\case{2.1.2.1.1.1}{$\Delta(k^{3} + 2k^{2} - 2k) = 1$}\quad {\small($A_1\supseteq \{1, k^{2} + k - 1, 2k^{2} + k - 2, k^{3}, k^{3} + 2k^{2} - 2k\},\ A_2\supseteq \{k, k^{2} + 2k - 2\},\ A_3\supseteq \{k + 2, k^{2}, k^{3} + k - 1\}$)}

Since $(k^{2} + k - 1) + (k-1)\cdot (k^{2} + 2k - 1) = k^{3} + 2k^{2} - 2k$ must not be monochromatic in $A_1$, either $\Delta(k^{2} + 2k - 1) = 2$ or $\Delta(k^{2} + 2k - 1) = 3$.

\case{2.1.2.1.1.1.1}{$\Delta(k^{2} + 2k - 1) = 2$}\quad {\small($A_1\supseteq \{1, k^{2} + k - 1, 2k^{2} + k - 2, k^{3}, k^{3} + 2k^{2} - 2k\},\ A_2\supseteq \{k, k^{2} + 2k - 2, k^{2} + 2k - 1\},\ A_3\supseteq \{k + 2, k^{2}, k^{3} + k - 1\}$)}

By Rule~2 applied to $(k^{3} + 2k^{2} - 2k) + (k-1)\cdot 1 = k^{3} + 2k^{2} - k - 1$ and $(k^{2} + 2k - 2) + (k-1)\cdot (k^{2} + 2k - 1) = k^{3} + 2k^{2} - k - 1$, $\Delta(k^{3} + 2k^{2} - k - 1) = 3$.

By Rule~2 applied to $k + (k-1)\cdot (k^{2} + 2k - 1) = k^{3} + k^{2} - 2k + 1$ and $(k^{3} + k^{2} - 2k + 1) + (k-1)\cdot (k + 2) = k^{3} + 2k^{2} - k - 1$, $\Delta(k^{3} + k^{2} - 2k + 1) = 1$.

By Rule~2 applied to $k^{3} + (k-1)\cdot 2k = k^{3} + 2k^{2} - 2k$ and $(k^{3} + k - 1) + (k-1)\cdot 2k = k^{3} + 2k^{2} - k - 1$, $\Delta(2k) = 2$.

By Rule~2 applied to $(k^{2} + k - 1) + (k-1)\cdot (k + 1) = 2k^{2} + k - 2$ and $2k + (k-1)\cdot (k + 1) = k^{2} + 2k - 1$, $\Delta(k + 1) = 3$.

By Rule~2 applied to $1 + (k-1)\cdot (k^{2} + 2k) = k^{3} + k^{2} - 2k + 1$ and $k\cdot (k + 2) = k^{2} + 2k$, $\Delta(k^{2} + 2k) = 2$.

\emph{Contradiction (Rule~3) at $k^{3} + k^{2} - k$:} if $\Delta(k^{3} + k^{2} - k)=1$ then $k\cdot (k^{2} + k - 1) = k^{3} + k^{2} - k$ is monochromatic; if $\Delta(k^{3} + k^{2} - k)=2$ then $k + (k-1)\cdot (k^{2} + 2k) = k^{3} + k^{2} - k$ is monochromatic; if $\Delta(k^{3} + k^{2} - k)=3$ then $(k^{3} + k^{2} - k) + (k-1)\cdot (k + 1) = k^{3} + 2k^{2} - k - 1$ is monochromatic.

\case{2.1.2.1.1.1.2}{$\Delta(k^{2} + 2k - 1) = 3$}\quad {\small($A_1\supseteq \{1, k^{2} + k - 1, 2k^{2} + k - 2, k^{3}, k^{3} + 2k^{2} - 2k\},\ A_2\supseteq \{k, k^{2} + 2k - 2\},\ A_3\supseteq \{k + 2, k^{2}, k^{2} + 2k - 1, k^{3} + k - 1\}$)}

By Rule~2 applied to $(k^{2} + k - 1) + (k-1)\cdot (k + 1) = 2k^{2} + k - 2$ and $(k + 1) + (k-1)\cdot (k + 2) = k^{2} + 2k - 1$, $\Delta(k + 1) = 2$.

By Rule~2 applied to $(2k - 1) + (k-1)\cdot (k^{2} + k - 1) = k^{3}$ and $(2k - 1) + (k-1)\cdot (k + 1) = k^{2} + 2k - 2$, $\Delta(2k - 1) = 3$.

By Rule~2 applied to $k\cdot (k^{2} + k - 1) = k^{3} + k^{2} - k$ and $(2k - 1) + (k-1)\cdot (k^{2} + 2k - 1) = k^{3} + k^{2} - k$, $\Delta(k^{3} + k^{2} - k) = 2$.

By Rule~2 applied to $(2k^{2} + k - 2) + (k-1)\cdot (k^{2} + k - 2) = k^{3} + 2k^{2} - 2k$ and $(k^{2} + 2k - 2) + (k-1)\cdot (k^{2} + k - 2) = k^{3} + k^{2} - k$, $\Delta(k^{2} + k - 2) = 3$.

\emph{Contradiction (Rule~3) at $k^{3} - k + 1$:} if $\Delta(k^{3} - k + 1)=1$ then $(k^{3} - k + 1) + (k-1)\cdot 1 = k^{3}$ is monochromatic; if $\Delta(k^{3} - k + 1)=2$ then $(k^{3} - k + 1) + (k-1)\cdot (k + 1) = k^{3} + k^{2} - k$ is monochromatic; if $\Delta(k^{3} - k + 1)=3$ then $(2k - 1) + (k-1)\cdot (k^{2} + k - 2) = k^{3} - k + 1$ is monochromatic.

\case{2.1.2.1.1.2}{$\Delta(k^{3} + 2k^{2} - 2k) = 3$}\quad {\small($A_1\supseteq \{1, k^{2} + k - 1, 2k^{2} + k - 2, k^{3}\},\ A_2\supseteq \{k, k^{2} + 2k - 2\},\ A_3\supseteq \{k + 2, k^{2}, k^{3} + k - 1, k^{3} + 2k^{2} - 2k\}$)}

By Rule~2 applied to $(2k - 1) + (k-1)\cdot (k^{2} + k - 1) = k^{3}$ and $(k^{3} + k - 1) + (k-1)\cdot (2k - 1) = k^{3} + 2k^{2} - 2k$, $\Delta(2k - 1) = 2$.

By Rule~2 applied to $(k^{2} + k - 1) + (k-1)\cdot (k + 1) = 2k^{2} + k - 2$ and $(2k - 1) + (k-1)\cdot (k + 1) = k^{2} + 2k - 2$, $\Delta(k + 1) = 3$.

By Rule~2 applied to $(2k^{2} - 1) + (k-1)\cdot 1 = 2k^{2} + k - 2$ and $k^{2} + (k-1)\cdot (k + 1) = 2k^{2} - 1$, $\Delta(2k^{2} - 1) = 2$.

By Rule~2 applied to $(2k - 1) + (k-1)\cdot (k^{2} + 2k - 2) = k^{3} + k^{2} - 2k + 1$ and $(k^{3} + k^{2} - 2k + 1) + (k-1)\cdot (k + 1) = k^{3} + 2k^{2} - 2k$, $\Delta(k^{3} + k^{2} - 2k + 1) = 1$.

By Rule~2 applied to $1 + (k-1)\cdot (k^{2} + 2k) = k^{3} + k^{2} - 2k + 1$ and $k^{2} + (k-1)\cdot (k^{2} + 2k) = k^{3} + 2k^{2} - 2k$, $\Delta(k^{2} + 2k) = 2$.

By Rule~2 applied to $k\cdot (k^{2} + k - 1) = k^{3} + k^{2} - k$ and $k + (k-1)\cdot (k^{2} + 2k) = k^{3} + k^{2} - k$, $\Delta(k^{3} + k^{2} - k) = 3$.

By Rule~2 applied to $(2k^{2} + k - 2) + (k-1)\cdot (k^{2} + k - 1) = k^{3} + 2k^{2} - k - 1$ and $(k^{3} + k^{2} - k) + (k-1)\cdot (k + 1) = k^{3} + 2k^{2} - k - 1$, $\Delta(k^{3} + 2k^{2} - k - 1) = 2$.

By Rule~2 applied to $(2k^{2} - 1) + (k-1)\cdot (k^{2} + k) = k^{3} + 2k^{2} - k - 1$ and $k\cdot (k + 1) = k^{2} + k$, $\Delta(k^{2} + k) = 1$.

\emph{Contradiction (Rule~3) at $k^{2} + 2k - 1$:} if $\Delta(k^{2} + 2k - 1)=1$ then $(k^{2} + k) + (k-1)\cdot 1 = k^{2} + 2k - 1$ is monochromatic; if $\Delta(k^{2} + 2k - 1)=2$ then $(k^{2} + 2k - 2) + (k-1)\cdot (k^{2} + 2k - 1) = k^{3} + 2k^{2} - k - 1$ is monochromatic; if $\Delta(k^{2} + 2k - 1)=3$ then $(k + 1) + (k-1)\cdot (k + 2) = k^{2} + 2k - 1$ is monochromatic.

\case{2.1.2.1.2}{$\Delta(k^{2} + 2k - 2) = 3$}\quad {\small($A_1\supseteq \{1, k^{2} + k - 1, k^{3}\},\ A_2\supseteq \{k\},\ A_3\supseteq \{k^{2}, k^{2} + 2k - 2, k^{3} + k - 1\}$)}

By Rule~2 applied to $(k^{2} + k - 1) + (k-1)\cdot (k^{2} - 1) = k^{3}$ and $(k^{2} + 2k - 2) + (k-1)\cdot (k^{2} - 1) = k^{3} + k - 1$, $\Delta(k^{2} - 1) = 2$.

Since $k\cdot (k^{2} + k - 1) = k^{3} + k^{2} - k$ must not be monochromatic in $A_1$, either $\Delta(k^{3} + k^{2} - k) = 2$ or $\Delta(k^{3} + k^{2} - k) = 3$.

\case{2.1.2.1.2.1}{$\Delta(k^{3} + k^{2} - k) = 2$}\quad {\small($A_1\supseteq \{1, k^{2} + k - 1, k^{3}\},\ A_2\supseteq \{k, k^{2} - 1, k^{3} + k^{2} - k\},\ A_3\supseteq \{k^{2}, k^{2} + 2k - 2, k^{3} + k - 1\}$)}

By Rule~2 applied to $(k^{3} + k^{2} - k) + (k-1)\cdot k = k^{3} + 2k^{2} - 2k$ and $k\cdot (k^{2} + 2k - 2) = k^{3} + 2k^{2} - 2k$, $\Delta(k^{3} + 2k^{2} - 2k) = 1$.

By Rule~2 applied to $(2k^{2} - 1) + (k-1)\cdot (k^{2} + k - 1) = k^{3} + 2k^{2} - 2k$ and $(2k^{2} - 1) + (k-1)\cdot (k^{2} - 1) = k^{3} + k^{2} - k$, $\Delta(2k^{2} - 1) = 3$.

By Rule~2 applied to $(k - 1) + (k-1)\cdot k = k^{2} - 1$ and $(k^{2} + 2k - 2) + (k-1)\cdot (k - 1) = 2k^{2} - 1$, $\Delta(k - 1) = 1$.

By Rule~2 applied to $(k^{3} + k^{2} - 1) + (k-1)\cdot (k - 1) = k^{3} + 2k^{2} - 2k$ and $(2k^{2} - 1) + (k-1)\cdot k^{2} = k^{3} + k^{2} - 1$, $\Delta(k^{3} + k^{2} - 1) = 2$.

By Rule~2 applied to $(k^{3} + 2k^{2} - 2k) + (k-1)\cdot 1 = k^{3} + 2k^{2} - k - 1$ and $(k^{3} + k^{2} - 1) + (k-1)\cdot k = k^{3} + 2k^{2} - k - 1$, $\Delta(k^{3} + 2k^{2} - k - 1) = 3$.

By Rule~2 applied to $k^{3} + (k-1)\cdot 2k = k^{3} + 2k^{2} - 2k$ and $(k^{3} + k - 1) + (k-1)\cdot 2k = k^{3} + 2k^{2} - k - 1$, $\Delta(2k) = 2$.

By Rule~2 applied to $(k^{2} + k - 1) + (k-1)\cdot (k^{2} + 2k - 1) = k^{3} + 2k^{2} - 2k$ and $(k^{2} + 2k - 2) + (k-1)\cdot (k^{2} + 2k - 1) = k^{3} + 2k^{2} - k - 1$, $\Delta(k^{2} + 2k - 1) = 2$.

By Rule~2 applied to $2k + (k-1)\cdot (k + 1) = k^{2} + 2k - 1$ and $k^{2} + (k-1)\cdot (k + 1) = 2k^{2} - 1$, $\Delta(k + 1) = 1$.

By Rule~2 applied to $2 + (k-1)\cdot 1 = k + 1$ and $k^{2} + (k-1)\cdot 2 = k^{2} + 2k - 2$, $\Delta(2) = 2$.

\emph{Contradiction:} $k\cdot 2 = 2k$ is now monochromatic.

\case{2.1.2.1.2.2}{$\Delta(k^{3} + k^{2} - k) = 3$}\quad {\small($A_1\supseteq \{1, k^{2} + k - 1, k^{3}\},\ A_2\supseteq \{k, k^{2} - 1\},\ A_3\supseteq \{k^{2}, k^{2} + 2k - 2, k^{3} + k - 1, k^{3} + k^{2} - k\}$)}

By Rule~2 applied to $(k - 1) + (k-1)\cdot k = k^{2} - 1$ and $(k^{3} + k - 1) + (k-1)\cdot (k - 1) = k^{3} + k^{2} - k$, $\Delta(k - 1) = 1$.

By Rule~2 applied to $(3k - 2) + (k-1)\cdot (k - 1) = k^{2} + k - 1$ and $(3k - 2) + (k-1)\cdot (k^{2} + 2k - 2) = k^{3} + k^{2} - k$, $\Delta(3k - 2) = 2$.

By Rule~2 applied to $k + (k-1)\cdot 2 = 3k - 2$ and $k^{2} + (k-1)\cdot 2 = k^{2} + 2k - 2$, $\Delta(2) = 1$.

Since $k^{3} + (k-1)\cdot (k - 1) = k^{3} + k^{2} - 2k + 1$ must not be monochromatic in $A_1$, either $\Delta(k^{3} + k^{2} - 2k + 1) = 2$ or $\Delta(k^{3} + k^{2} - 2k + 1) = 3$.

\case{2.1.2.1.2.2.1}{$\Delta(k^{3} + k^{2} - 2k + 1) = 2$}\quad {\small($A_1\supseteq \{1, 2, k - 1, k^{2} + k - 1, k^{3}\},\ A_2\supseteq \{k, 3k - 2, k^{2} - 1, k^{3} + k^{2} - 2k + 1\},\ A_3\supseteq \{k^{2}, k^{2} + 2k - 2, k^{3} + k - 1, k^{3} + k^{2} - k\}$)}

By Rule~2 applied to $(k^{3} - k + 1) + (k-1)\cdot 1 = k^{3}$ and $(k^{3} - k + 1) + (k-1)\cdot k = k^{3} + k^{2} - 2k + 1$, $\Delta(k^{3} - k + 1) = 3$.

By Rule~2 applied to $2 + (k-1)\cdot 1 = k + 1$ and $(k^{3} - k + 1) + (k-1)\cdot (k + 1) = k^{3} + k^{2} - k$, $\Delta(k + 1) = 2$.

By Rule~2 applied to $k\cdot (k + 1) = k^{2} + k$ and $k^{2} + (k-1)\cdot (k^{2} + k) = k^{3} + k^{2} - k$, $\Delta(k^{2} + k) = 1$.

By Rule~2 applied to $(k^{2} + k) + (k-1)\cdot 1 = k^{2} + 2k - 1$ and $k + (k-1)\cdot (k^{2} + 2k - 1) = k^{3} + k^{2} - 2k + 1$, $\Delta(k^{2} + 2k - 1) = 3$.

By Rule~2 applied to $(k^{3} + k^{2} - 2k + 1) + (k-1)\cdot k = k^{3} + 2k^{2} - 3k + 1$ and $k^{2} + (k-1)\cdot (k^{2} + 2k - 1) = k^{3} + 2k^{2} - 3k + 1$, $\Delta(k^{3} + 2k^{2} - 3k + 1) = 1$.

\emph{Contradiction (Rule~3) at $k^{3} + 2k^{2} - 2k$:} if $\Delta(k^{3} + 2k^{2} - 2k)=1$ then $(k^{3} + 2k^{2} - 3k + 1) + (k-1)\cdot 1 = k^{3} + 2k^{2} - 2k$ is monochromatic; if $\Delta(k^{3} + 2k^{2} - 2k)=2$ then $(k^{3} + k^{2} - 2k + 1) + (k-1)\cdot (k + 1) = k^{3} + 2k^{2} - 2k$ is monochromatic; if $\Delta(k^{3} + 2k^{2} - 2k)=3$ then $k\cdot (k^{2} + 2k - 2) = k^{3} + 2k^{2} - 2k$ is monochromatic.

\case{2.1.2.1.2.2.2}{$\Delta(k^{3} + k^{2} - 2k + 1) = 3$}\quad {\small($A_1\supseteq \{1, 2, k - 1, k^{2} + k - 1, k^{3}\},\ A_2\supseteq \{k, 3k - 2, k^{2} - 1\},\ A_3\supseteq \{k^{2}, k^{2} + 2k - 2, k^{3} + k - 1, k^{3} + k^{2} - 2k + 1, k^{3} + k^{2} - k\}$)}

By Rule~2 applied to $1 + (k-1)\cdot 2 = 2k - 1$ and $(2k - 1) + (k-1)\cdot (k^{2} + 2k - 2) = k^{3} + k^{2} - 2k + 1$, $\Delta(2k - 1) = 2$.

By Rule~2 applied to $k + (k-1)\cdot (2k - 1) = 2k^{2} - 2k + 1$ and $(2k^{2} - 2k + 1) + (k-1)\cdot k^{2} = k^{3} + k^{2} - 2k + 1$, $\Delta(2k^{2} - 2k + 1) = 1$.

By Rule~2 applied to $(2k^{2} - 2k + 1) + (k-1)\cdot 2 = 2k^{2} - 1$ and $(3k - 2) + (k-1)\cdot (2k - 1) = 2k^{2} - 1$, $\Delta(2k^{2} - 1) = 3$.

By Rule~2 applied to $2 + (k-1)\cdot 1 = k + 1$ and $k^{2} + (k-1)\cdot (k + 1) = 2k^{2} - 1$, $\Delta(k + 1) = 2$.

By Rule~2 applied to $k\cdot (k + 1) = k^{2} + k$ and $k^{2} + (k-1)\cdot (k^{2} + k) = k^{3} + k^{2} - k$, $\Delta(k^{2} + k) = 1$.

By Rule~2 applied to $(k^{2} + k - 1) + (k-1)\cdot (k^{2} + k) = k^{3} + k^{2} - 1$ and $(2k^{2} - 1) + (k-1)\cdot k^{2} = k^{3} + k^{2} - 1$, $\Delta(k^{3} + k^{2} - 1) = 2$.

By Rule~2 applied to $(k^{2} + k) + (k-1)\cdot 1 = k^{2} + 2k - 1$ and $(3k - 2) + (k-1)\cdot (k^{2} + 2k - 1) = k^{3} + k^{2} - 1$, $\Delta(k^{2} + 2k - 1) = 3$.

By Rule~2 applied to $(k^{3} + k^{2} - 1) + (k-1)\cdot k = k^{3} + 2k^{2} - k - 1$ and $(k^{2} + 2k - 2) + (k-1)\cdot (k^{2} + 2k - 1) = k^{3} + 2k^{2} - k - 1$, $\Delta(k^{3} + 2k^{2} - k - 1) = 1$.

By Rule~2 applied to $(k^{3} + 2k^{2} - k - 1) + (k-1)\cdot 1 = k^{3} + 2k^{2} - 2$ and $(k^{3} + k^{2} - 1) + (k-1)\cdot (k + 1) = k^{3} + 2k^{2} - 2$, $\Delta(k^{3} + 2k^{2} - 2) = 3$.

By Rule~2 applied to $1 + (k-1)\cdot (k + 2) = k^{2} + k - 1$ and $(k^{3} + k^{2} - k) + (k-1)\cdot (k + 2) = k^{3} + 2k^{2} - 2$, $\Delta(k + 2) = 2$.

\emph{Contradiction (Rule~3) at $k^{2} + k + 1$:} if $\Delta(k^{2} + k + 1)=1$ then $1 + (k-1)\cdot (k^{2} + k + 1) = k^{3}$ is monochromatic; if $\Delta(k^{2} + k + 1)=2$ then $(k + 2) + (k-1)\cdot (k + 1) = k^{2} + k + 1$ is monochromatic; if $\Delta(k^{2} + k + 1)=3$ then $(2k^{2} - 1) + (k-1)\cdot (k^{2} + k + 1) = k^{3} + 2k^{2} - 2$ is monochromatic.

\case{2.1.2.2}{$\Delta(k^{2} + k - 1) = 2$}\quad {\small($A_1\supseteq \{1, k^{3}\},\ A_2\supseteq \{k, k^{2} + k - 1\},\ A_3\supseteq \{k^{2}, k^{3} + k - 1\}$)}

By Rule~2 applied to $(k^{3} - k + 1) + (k-1)\cdot 1 = k^{3}$ and $k + (k-1)\cdot (k^{2} + k - 1) = k^{3} - k + 1$, $\Delta(k^{3} - k + 1) = 3$.

Since $k + (k-1)\cdot (k + 1) = k^{2} + k - 1$ must not be monochromatic in $A_2$, either $\Delta(k + 1) = 1$ or $\Delta(k + 1) = 3$.

\case{2.1.2.2.1}{$\Delta(k + 1) = 1$}\quad {\small($A_1\supseteq \{1, k + 1, k^{3}\},\ A_2\supseteq \{k, k^{2} + k - 1\},\ A_3\supseteq \{k^{2}, k^{3} - k + 1, k^{3} + k - 1\}$)}

By Rule~2 applied to $2 + (k-1)\cdot 1 = k + 1$ and $(k^{3} - k + 1) + (k-1)\cdot 2 = k^{3} + k - 1$, $\Delta(2) = 2$.

By Rule~2 applied to $(k + 1) + (k-1)\cdot 1 = 2k$ and $k\cdot 2 = 2k$, $\Delta(2k) = 3$.

Since $(2k - 1) + (k-1)\cdot k = k^{2} + k - 1$ must not be monochromatic in $A_2$, either $\Delta(2k - 1) = 1$ or $\Delta(2k - 1) = 3$.

\case{2.1.2.2.1.1}{$\Delta(2k - 1) = 1$}\quad {\small($A_1\supseteq \{1, k + 1, 2k - 1, k^{3}\},\ A_2\supseteq \{2, k, k^{2} + k - 1\},\ A_3\supseteq \{2k, k^{2}, k^{3} - k + 1, k^{3} + k - 1\}$)}

By Rule~2 applied to $(2k - 1) + (k-1)\cdot 1 = 3k - 2$ and $k + (k-1)\cdot 2 = 3k - 2$, $\Delta(3k - 2) = 3$.

By Rule~2 applied to $k^{3} + (k-1)\cdot (2k - 1) = k^{3} + 2k^{2} - 3k + 1$ and $(k^{3} - k + 1) + (k-1)\cdot 2k = k^{3} + 2k^{2} - 3k + 1$, $\Delta(k^{3} + 2k^{2} - 3k + 1) = 2$.

By Rule~2 applied to $(2k - 1) + (k-1)\cdot (k + 1) = k^{2} + 2k - 2$ and $(k^{2} + k - 1) + (k-1)\cdot (k^{2} + 2k - 2) = k^{3} + 2k^{2} - 3k + 1$, $\Delta(k^{2} + 2k - 2) = 3$.

By Rule~2 applied to $(k^{3} + 2k^{2} - 3k + 1) + (k-1)\cdot 2 = k^{3} + 2k^{2} - k - 1$ and $(k^{3} + k - 1) + (k-1)\cdot 2k = k^{3} + 2k^{2} - k - 1$, $\Delta(k^{3} + 2k^{2} - k - 1) = 1$.

\emph{Contradiction (Rule~3) at $k^{3} + k^{2} - k$:} if $\Delta(k^{3} + k^{2} - k)=1$ then $(k^{3} + k^{2} - k) + (k-1)\cdot (k + 1) = k^{3} + 2k^{2} - k - 1$ is monochromatic; if $\Delta(k^{3} + k^{2} - k)=2$ then $k\cdot (k^{2} + k - 1) = k^{3} + k^{2} - k$ is monochromatic; if $\Delta(k^{3} + k^{2} - k)=3$ then $(3k - 2) + (k-1)\cdot (k^{2} + 2k - 2) = k^{3} + k^{2} - k$ is monochromatic.

\case{2.1.2.2.1.2}{$\Delta(2k - 1) = 3$}\quad {\small($A_1\supseteq \{1, k + 1, k^{3}\},\ A_2\supseteq \{2, k, k^{2} + k - 1\},\ A_3\supseteq \{2k - 1, 2k, k^{2}, k^{3} - k + 1, k^{3} + k - 1\}$)}

By Rule~2 applied to $k\cdot (k + 1) = k^{2} + k$ and $(2k - 1) + (k-1)\cdot (k^{2} + k) = k^{3} + k - 1$, $\Delta(k^{2} + k) = 2$.

By Rule~2 applied to $k^{3} + (k-1)\cdot (k + 1) = k^{3} + k^{2} - 1$ and $(k^{2} + k - 1) + (k-1)\cdot (k^{2} + k) = k^{3} + k^{2} - 1$, $\Delta(k^{3} + k^{2} - 1) = 3$.

By Rule~2 applied to $2 + (k-1)\cdot (k + 2) = k^{2} + k$ and $(k^{3} - k + 1) + (k-1)\cdot (k + 2) = k^{3} + k^{2} - 1$, $\Delta(k + 2) = 1$.

By Rule~2 applied to $k\cdot (k + 2) = k^{2} + 2k$ and $(2k - 1) + (k-1)\cdot (k^{2} + 2k) = k^{3} + k^{2} - 1$, $\Delta(k^{2} + 2k) = 2$.

By Rule~2 applied to $(k^{2} + k - 1) + (k-1)\cdot (k^{2} + 2k) = k^{3} + 2k^{2} - k - 1$ and $(k^{3} + k - 1) + (k-1)\cdot 2k = k^{3} + 2k^{2} - k - 1$, $\Delta(k^{3} + 2k^{2} - k - 1) = 1$.

By Rule~2 applied to $(k^{3} + 2k^{2} - 2k) + (k-1)\cdot 1 = k^{3} + 2k^{2} - k - 1$ and $(k^{3} + k - 1) + (k-1)\cdot (2k - 1) = k^{3} + 2k^{2} - 2k$, $\Delta(k^{3} + 2k^{2} - 2k) = 2$.

By Rule~2 applied to $(k^{3} + k^{2} - k) + (k-1)\cdot (k + 1) = k^{3} + 2k^{2} - k - 1$ and $k\cdot (k^{2} + k - 1) = k^{3} + k^{2} - k$, $\Delta(k^{3} + k^{2} - k) = 3$.

\emph{Contradiction (Rule~3) at $k^{2} + 2k - 1$:} if $\Delta(k^{2} + 2k - 1)=1$ then $(k + 1) + (k-1)\cdot (k + 2) = k^{2} + 2k - 1$ is monochromatic; if $\Delta(k^{2} + 2k - 1)=2$ then $(k^{2} + k - 1) + (k-1)\cdot (k^{2} + 2k - 1) = k^{3} + 2k^{2} - 2k$ is monochromatic; if $\Delta(k^{2} + 2k - 1)=3$ then $(2k - 1) + (k-1)\cdot (k^{2} + 2k - 1) = k^{3} + k^{2} - k$ is monochromatic.

\case{2.1.2.2.2}{$\Delta(k + 1) = 3$}\quad {\small($A_1\supseteq \{1, k^{3}\},\ A_2\supseteq \{k, k^{2} + k - 1\},\ A_3\supseteq \{k + 1, k^{2}, k^{3} - k + 1, k^{3} + k - 1\}$)}

By Rule~2 applied to $(k^{2} + k - 1) + (k-1)\cdot k = 2k^{2} - 1$ and $k^{2} + (k-1)\cdot (k + 1) = 2k^{2} - 1$, $\Delta(2k^{2} - 1) = 1$.

By Rule~2 applied to $k\cdot (k^{2} + k - 1) = k^{3} + k^{2} - k$ and $(k^{3} - k + 1) + (k-1)\cdot (k + 1) = k^{3} + k^{2} - k$, $\Delta(k^{3} + k^{2} - k) = 1$.

By Rule~2 applied to $(2k^{2} - 1) + (k-1)\cdot (k^{2} - 1) = k^{3} + k^{2} - k$ and $k^{2} + (k-1)\cdot (k^{2} - 1) = k^{3} - k + 1$, $\Delta(k^{2} - 1) = 2$.

Since $(2k - 1) + (k-1)\cdot k = k^{2} + k - 1$ must not be monochromatic in $A_2$, either $\Delta(2k - 1) = 1$ or $\Delta(2k - 1) = 3$.

\case{2.1.2.2.2.1}{$\Delta(2k - 1) = 1$}\quad {\small($A_1\supseteq \{1, 2k - 1, 2k^{2} - 1, k^{3}, k^{3} + k^{2} - k\},\ A_2\supseteq \{k, k^{2} - 1, k^{2} + k - 1\},\ A_3\supseteq \{k + 1, k^{2}, k^{3} - k + 1, k^{3} + k - 1\}$)}

By Rule~2 applied to $1 + (k-1)\cdot 2 = 2k - 1$ and $(k^{3} - k + 1) + (k-1)\cdot 2 = k^{3} + k - 1$, $\Delta(2) = 2$.

By Rule~2 applied to $(2k - 1) + (k-1)\cdot 2k = 2k^{2} - 1$ and $k\cdot 2 = 2k$, $\Delta(2k) = 3$.

By Rule~2 applied to $k^{3} + (k-1)\cdot (2k - 1) = k^{3} + 2k^{2} - 3k + 1$ and $(k^{3} - k + 1) + (k-1)\cdot 2k = k^{3} + 2k^{2} - 3k + 1$, $\Delta(k^{3} + 2k^{2} - 3k + 1) = 2$.

By Rule~2 applied to $(2k - 1) + (k-1)\cdot (k^{2} + 2k - 1) = k^{3} + k^{2} - k$ and $2k + (k-1)\cdot (k + 1) = k^{2} + 2k - 1$, $\Delta(k^{2} + 2k - 1) = 2$.

By Rule~2 applied to $(k^{3} + k^{2} - 2k + 1) + (k-1)\cdot 1 = k^{3} + k^{2} - k$ and $(k^{3} + k^{2} - 2k + 1) + (k-1)\cdot k = k^{3} + 2k^{2} - 3k + 1$, $\Delta(k^{3} + k^{2} - 2k + 1) = 3$.

By Rule~2 applied to $(k^{3} + 2k^{2} - 3k + 1) + (k-1)\cdot 2 = k^{3} + 2k^{2} - k - 1$ and $(k^{3} + k - 1) + (k-1)\cdot 2k = k^{3} + 2k^{2} - k - 1$, $\Delta(k^{3} + 2k^{2} - k - 1) = 1$.

\emph{Contradiction (Rule~3) at $k^{3} + 2k^{2} - 2k$:} if $\Delta(k^{3} + 2k^{2} - 2k)=1$ then $(k^{3} + 2k^{2} - 2k) + (k-1)\cdot 1 = k^{3} + 2k^{2} - k - 1$ is monochromatic; if $\Delta(k^{3} + 2k^{2} - 2k)=2$ then $(k^{2} + k - 1) + (k-1)\cdot (k^{2} + 2k - 1) = k^{3} + 2k^{2} - 2k$ is monochromatic; if $\Delta(k^{3} + 2k^{2} - 2k)=3$ then $(k^{3} + k^{2} - 2k + 1) + (k-1)\cdot (k + 1) = k^{3} + 2k^{2} - 2k$ is monochromatic.

\case{2.1.2.2.2.2}{$\Delta(2k - 1) = 3$}\quad {\small($A_1\supseteq \{1, 2k^{2} - 1, k^{3}, k^{3} + k^{2} - k\},\ A_2\supseteq \{k, k^{2} - 1, k^{2} + k - 1\},\ A_3\supseteq \{k + 1, 2k - 1, k^{2}, k^{3} - k + 1, k^{3} + k - 1\}$)}

By Rule~2 applied to $(k - 1) + (k-1)\cdot k = k^{2} - 1$ and $(2k - 1) + (k-1)\cdot (k - 1) = k^{2}$, $\Delta(k - 1) = 1$.

By Rule~2 applied to $(k^{2} + 2k - 2) + (k-1)\cdot (k - 1) = 2k^{2} - 1$ and $(2k - 1) + (k-1)\cdot (k + 1) = k^{2} + 2k - 2$, $\Delta(k^{2} + 2k - 2) = 2$.

By Rule~2 applied to $k\cdot (k^{2} + 2k - 2) = k^{3} + 2k^{2} - 2k$ and $(k^{3} + k - 1) + (k-1)\cdot (2k - 1) = k^{3} + 2k^{2} - 2k$, $\Delta(k^{3} + 2k^{2} - 2k) = 1$.

By Rule~2 applied to $(k^{3} + k^{2} - k) + (k-1)\cdot 1 = k^{3} + k^{2} - 1$ and $(k^{2} + 2k - 2) + (k-1)\cdot (k^{2} + k - 1) = k^{3} + k^{2} - 1$, $\Delta(k^{3} + k^{2} - 1) = 3$.

By Rule~2 applied to $(k^{3} + k^{2} - k) + (k-1)\cdot (k - 1) = k^{3} + 2k^{2} - 3k + 1$ and $(k^{2} + k - 1) + (k-1)\cdot (k^{2} + 2k - 2) = k^{3} + 2k^{2} - 3k + 1$, $\Delta(k^{3} + 2k^{2} - 3k + 1) = 3$.

By Rule~2 applied to $k + (k-1)\cdot (k + 2) = k^{2} + 2k - 2$ and $(k^{3} - k + 1) + (k-1)\cdot (k + 2) = k^{3} + k^{2} - 1$, $\Delta(k + 2) = 1$.

By Rule~2 applied to $k^{3} + (k-1)\cdot 2k = k^{3} + 2k^{2} - 2k$ and $(k^{3} - k + 1) + (k-1)\cdot 2k = k^{3} + 2k^{2} - 3k + 1$, $\Delta(2k) = 2$.

By Rule~2 applied to $2k + (k-1)\cdot k = k^{2} + k$ and $(2k - 1) + (k-1)\cdot (k^{2} + k) = k^{3} + k - 1$, $\Delta(k^{2} + k) = 1$.

By Rule~2 applied to $2 + (k-1)\cdot (k + 2) = k^{2} + k$ and $(k^{3} - k + 1) + (k-1)\cdot 2 = k^{3} + k - 1$, $\Delta(2) = 2$.

\emph{Contradiction:} $k\cdot 2 = 2k$ is now monochromatic.

\case{2.2}{$\Delta(k^{3}) = 2$}\quad {\small($A_1\supseteq \{1\},\ A_2\supseteq \{k, k^{3}\},\ A_3\supseteq \{k^{2}\}$)}

Since $k^{3} + (k-1)\cdot k = k^{3} + k^{2} - k$ must not be monochromatic in $A_2$, either $\Delta(k^{3} + k^{2} - k) = 1$ or $\Delta(k^{3} + k^{2} - k) = 3$.

\case{2.2.1}{$\Delta(k^{3} + k^{2} - k) = 1$}\quad {\small($A_1\supseteq \{1, k^{3} + k^{2} - k\},\ A_2\supseteq \{k, k^{3}\},\ A_3\supseteq \{k^{2}\}$)}

Since $(k^{3} + k^{2} - 2k + 1) + (k-1)\cdot 1 = k^{3} + k^{2} - k$ must not be monochromatic in $A_1$, either $\Delta(k^{3} + k^{2} - 2k + 1) = 2$ or $\Delta(k^{3} + k^{2} - 2k + 1) = 3$.

\case{2.2.1.1}{$\Delta(k^{3} + k^{2} - 2k + 1) = 2$}\quad {\small($A_1\supseteq \{1, k^{3} + k^{2} - k\},\ A_2\supseteq \{k, k^{3}, k^{3} + k^{2} - 2k + 1\},\ A_3\supseteq \{k^{2}\}$)}

Since $(k^{3} + k^{2} - k) + (k-1)\cdot 1 = k^{3} + k^{2} - 1$ must not be monochromatic in $A_1$, either $\Delta(k^{3} + k^{2} - 1) = 2$ or $\Delta(k^{3} + k^{2} - 1) = 3$.

\case{2.2.1.1.1}{$\Delta(k^{3} + k^{2} - 1) = 2$}\quad {\small($A_1\supseteq \{1, k^{3} + k^{2} - k\},\ A_2\supseteq \{k, k^{3}, k^{3} + k^{2} - 2k + 1, k^{3} + k^{2} - 1\},\ A_3\supseteq \{k^{2}\}$)}

Since $k^{3} + (k-1)\cdot (k + 1) = k^{3} + k^{2} - 1$ must not be monochromatic in $A_2$, either $\Delta(k + 1) = 1$ or $\Delta(k + 1) = 3$.

\case{2.2.1.1.1.1}{$\Delta(k + 1) = 1$}\quad {\small($A_1\supseteq \{1, k + 1, k^{3} + k^{2} - k\},\ A_2\supseteq \{k, k^{3}, k^{3} + k^{2} - 2k + 1, k^{3} + k^{2} - 1\},\ A_3\supseteq \{k^{2}\}$)}

By Rule~2 applied to $2 + (k-1)\cdot 1 = k + 1$ and $(k^{3} + k^{2} - 2k + 1) + (k-1)\cdot 2 = k^{3} + k^{2} - 1$, $\Delta(2) = 3$.

By Rule~2 applied to $(k + 1) + (k-1)\cdot 1 = 2k$ and $k\cdot 2 = 2k$, $\Delta(2k) = 2$.

By Rule~2 applied to $k\cdot (k + 1) = k^{2} + k$ and $k + (k-1)\cdot (k^{2} + k) = k^{3}$, $\Delta(k^{2} + k) = 3$.

By Rule~2 applied to $(k^{3} - k + 1) + (k-1)\cdot (k + 1) = k^{3} + k^{2} - k$ and $(k^{3} - k + 1) + (k-1)\cdot k = k^{3} + k^{2} - 2k + 1$, $\Delta(k^{3} - k + 1) = 3$.

By Rule~2 applied to $(k^{3} + k^{2} - k) + (k-1)\cdot (k + 1) = k^{3} + 2k^{2} - k - 1$ and $(k^{3} + k^{2} - 1) + (k-1)\cdot k = k^{3} + 2k^{2} - k - 1$, $\Delta(k^{3} + 2k^{2} - k - 1) = 3$.

By Rule~2 applied to $(k^{3} + k - 1) + (k-1)\cdot k = k^{3} + k^{2} - 1$ and $(k^{3} - k + 1) + (k-1)\cdot 2 = k^{3} + k - 1$, $\Delta(k^{3} + k - 1) = 1$.

By Rule~2 applied to $(k^{3} + k - 1) + (k-1)\cdot (k - 1) = k^{3} + k^{2} - k$ and $k^{3} + (k-1)\cdot (k - 1) = k^{3} + k^{2} - 2k + 1$, $\Delta(k - 1) = 3$.

By Rule~2 applied to $(k^{3} + k^{2} - 2k + 1) + (k-1)\cdot k = k^{3} + 2k^{2} - 3k + 1$ and $(k^{3} + 2k^{2} - 3k + 1) + (k-1)\cdot 2 = k^{3} + 2k^{2} - k - 1$, $\Delta(k^{3} + 2k^{2} - 3k + 1) = 1$.

By Rule~2 applied to $(k^{3} + 2k^{2} - 3k + 1) + (k-1)\cdot 1 = k^{3} + 2k^{2} - 2k$ and $k^{3} + (k-1)\cdot 2k = k^{3} + 2k^{2} - 2k$, $\Delta(k^{3} + 2k^{2} - 2k) = 3$.

Since $2 + (k-1)\cdot (k + 2) = k^{2} + k$ must not be monochromatic in $A_3$, either $\Delta(k + 2) = 1$ or $\Delta(k + 2) = 2$.

\case{2.2.1.1.1.1.1}{$\Delta(k + 2) = 1$}\quad {\small($A_1\supseteq \{1, k + 1, k + 2, k^{3} + k - 1, k^{3} + k^{2} - k, k^{3} + 2k^{2} - 3k + 1\},\ A_2\supseteq \{k, 2k, k^{3}, k^{3} + k^{2} - 2k + 1, k^{3} + k^{2} - 1\},\ A_3\supseteq \{2, k - 1, k^{2}, k^{2} + k, k^{3} - k + 1, k^{3} + 2k^{2} - 2k, k^{3} + 2k^{2} - k - 1\}$)}

By Rule~2 applied to $(k + 1) + (k-1)\cdot (k + 2) = k^{2} + 2k - 1$ and $k + (k-1)\cdot (k^{2} + 2k - 1) = k^{3} + k^{2} - 2k + 1$, $\Delta(k^{2} + 2k - 1) = 3$.

By Rule~2 applied to $1 + (k-1)\cdot (k + 2) = k^{2} + k - 1$ and $(k^{2} + k - 1) + (k-1)\cdot (k^{2} + 2k - 1) = k^{3} + 2k^{2} - 2k$, $\Delta(k^{2} + k - 1) = 2$.

By Rule~2 applied to $(k^{2} + k - 1) + (k-1)\cdot k = 2k^{2} - 1$ and $(2k^{2} - 1) + (k-1)\cdot (k^{2} + k) = k^{3} + 2k^{2} - k - 1$, $\Delta(2k^{2} - 1) = 1$.

\emph{Contradiction (Rule~3) at $k^{2} - 1$:} if $\Delta(k^{2} - 1)=1$ then $(2k^{2} - 1) + (k-1)\cdot (k^{2} - 1) = k^{3} + k^{2} - k$ is monochromatic; if $\Delta(k^{2} - 1)=2$ then $(k^{2} + k - 1) + (k-1)\cdot (k^{2} - 1) = k^{3}$ is monochromatic; if $\Delta(k^{2} - 1)=3$ then $k^{2} + (k-1)\cdot (k^{2} - 1) = k^{3} - k + 1$ is monochromatic.

\case{2.2.1.1.1.1.2}{$\Delta(k + 2) = 2$}\quad {\small($A_1\supseteq \{1, k + 1, k^{3} + k - 1, k^{3} + k^{2} - k, k^{3} + 2k^{2} - 3k + 1\},\ A_2\supseteq \{k, k + 2, 2k, k^{3}, k^{3} + k^{2} - 2k + 1, k^{3} + k^{2} - 1\},\ A_3\supseteq \{2, k - 1, k^{2}, k^{2} + k, k^{3} - k + 1, k^{3} + 2k^{2} - 2k, k^{3} + 2k^{2} - k - 1\}$)}

By Rule~2 applied to $k + (k-1)\cdot (k + 2) = k^{2} + 2k - 2$ and $k^{2} + (k-1)\cdot 2 = k^{2} + 2k - 2$, $\Delta(k^{2} + 2k - 2) = 1$.

By Rule~2 applied to $(2k - 1) + (k-1)\cdot (k + 1) = k^{2} + 2k - 2$ and $(2k - 1) + (k-1)\cdot (k - 1) = k^{2}$, $\Delta(2k - 1) = 2$.

By Rule~2 applied to $(2k - 1) + (k-1)\cdot 2k = 2k^{2} - 1$ and $(2k^{2} - 1) + (k-1)\cdot (k^{2} + k) = k^{3} + 2k^{2} - k - 1$, $\Delta(2k^{2} - 1) = 1$.

By Rule~2 applied to $(k^{2} + 2k - 2) + (k-1)\cdot (k^{2} - 1) = k^{3} + k - 1$ and $k^{2} + (k-1)\cdot (k^{2} - 1) = k^{3} - k + 1$, $\Delta(k^{2} - 1) = 2$.

By Rule~2 applied to $(k^{2} + k - 1) + (k-1)\cdot 1 = k^{2} + 2k - 2$ and $(2k - 1) + (k-1)\cdot k = k^{2} + k - 1$, $\Delta(k^{2} + k - 1) = 3$.

\emph{Contradiction (Rule~3) at $2k^{2} + k - 2$:} if $\Delta(2k^{2} + k - 2)=1$ then $(2k^{2} - 1) + (k-1)\cdot 1 = 2k^{2} + k - 2$ is monochromatic; if $\Delta(2k^{2} + k - 2)=2$ then $(2k^{2} + k - 2) + (k-1)\cdot (k^{2} - 1) = k^{3} + k^{2} - 1$ is monochromatic; if $\Delta(2k^{2} + k - 2)=3$ then $(2k^{2} + k - 2) + (k-1)\cdot (k^{2} + k - 1) = k^{3} + 2k^{2} - k - 1$ is monochromatic.

\case{2.2.1.1.1.2}{$\Delta(k + 1) = 3$}\quad {\small($A_1\supseteq \{1, k^{3} + k^{2} - k\},\ A_2\supseteq \{k, k^{3}, k^{3} + k^{2} - 2k + 1, k^{3} + k^{2} - 1\},\ A_3\supseteq \{k + 1, k^{2}\}$)}

By Rule~2 applied to $k + (k-1)\cdot (k^{2} + k) = k^{3}$ and $k\cdot (k + 1) = k^{2} + k$, $\Delta(k^{2} + k) = 1$.

By Rule~2 applied to $1 + (k-1)\cdot (k^{2} + k) = k^{3} - k + 1$ and $(k^{3} - k + 1) + (k-1)\cdot k = k^{3} + k^{2} - 2k + 1$, $\Delta(k^{3} - k + 1) = 3$.

By Rule~2 applied to $(k^{2} + k) + (k-1)\cdot 1 = k^{2} + 2k - 1$ and $k + (k-1)\cdot (k^{2} + 2k - 1) = k^{3} + k^{2} - 2k + 1$, $\Delta(k^{2} + 2k - 1) = 3$.

By Rule~2 applied to $(k^{3} + k^{2} - 2k + 1) + (k-1)\cdot k = k^{3} + 2k^{2} - 3k + 1$ and $k^{2} + (k-1)\cdot (k^{2} + 2k - 1) = k^{3} + 2k^{2} - 3k + 1$, $\Delta(k^{3} + 2k^{2} - 3k + 1) = 1$.

By Rule~2 applied to $(k^{3} + k^{2} - k) + (k-1)\cdot (k - 1) = k^{3} + 2k^{2} - 3k + 1$ and $k^{3} + (k-1)\cdot (k - 1) = k^{3} + k^{2} - 2k + 1$, $\Delta(k - 1) = 3$.

By Rule~2 applied to $(2k^{2} - 2k + 1) + (k-1)\cdot (k^{2} + k) = k^{3} + 2k^{2} - 3k + 1$ and $k^{2} + (k-1)\cdot (k - 1) = 2k^{2} - 2k + 1$, $\Delta(2k^{2} - 2k + 1) = 2$.

By Rule~2 applied to $k + (k-1)\cdot (2k - 1) = 2k^{2} - 2k + 1$ and $(2k - 1) + (k-1)\cdot (k - 1) = k^{2}$, $\Delta(2k - 1) = 1$.

By Rule~2 applied to $1 + (k-1)\cdot 2 = 2k - 1$ and $(k^{3} + k^{2} - 2k + 1) + (k-1)\cdot 2 = k^{3} + k^{2} - 1$, $\Delta(2) = 3$.

\emph{Contradiction (Rule~3) at $k^{3} + k - 1$:} if $\Delta(k^{3} + k - 1)=1$ then $(2k - 1) + (k-1)\cdot (k^{2} + k) = k^{3} + k - 1$ is monochromatic; if $\Delta(k^{3} + k - 1)=2$ then $(k^{3} + k - 1) + (k-1)\cdot k = k^{3} + k^{2} - 1$ is monochromatic; if $\Delta(k^{3} + k - 1)=3$ then $(k^{3} - k + 1) + (k-1)\cdot 2 = k^{3} + k - 1$ is monochromatic.

\case{2.2.1.1.2}{$\Delta(k^{3} + k^{2} - 1) = 3$}\quad {\small($A_1\supseteq \{1, k^{3} + k^{2} - k\},\ A_2\supseteq \{k, k^{3}, k^{3} + k^{2} - 2k + 1\},\ A_3\supseteq \{k^{2}, k^{3} + k^{2} - 1\}$)}

Since $k + (k-1)\cdot (k^{2} + k) = k^{3}$ must not be monochromatic in $A_2$, either $\Delta(k^{2} + k) = 1$ or $\Delta(k^{2} + k) = 3$.

\case{2.2.1.1.2.1}{$\Delta(k^{2} + k) = 1$}\quad {\small($A_1\supseteq \{1, k^{2} + k, k^{3} + k^{2} - k\},\ A_2\supseteq \{k, k^{3}, k^{3} + k^{2} - 2k + 1\},\ A_3\supseteq \{k^{2}, k^{3} + k^{2} - 1\}$)}

By Rule~2 applied to $(k^{2} + k) + (k-1)\cdot 1 = k^{2} + 2k - 1$ and $k + (k-1)\cdot (k^{2} + 2k - 1) = k^{3} + k^{2} - 2k + 1$, $\Delta(k^{2} + 2k - 1) = 3$.

By Rule~2 applied to $(k^{3} + k^{2} - 2k + 1) + (k-1)\cdot k = k^{3} + 2k^{2} - 3k + 1$ and $k^{2} + (k-1)\cdot (k^{2} + 2k - 1) = k^{3} + 2k^{2} - 3k + 1$, $\Delta(k^{3} + 2k^{2} - 3k + 1) = 1$.

By Rule~2 applied to $(k^{3} + k^{2} - k) + (k-1)\cdot (k - 1) = k^{3} + 2k^{2} - 3k + 1$ and $k^{3} + (k-1)\cdot (k - 1) = k^{3} + k^{2} - 2k + 1$, $\Delta(k - 1) = 3$.

By Rule~2 applied to $(k^{3} + 2k^{2} - 3k + 1) + (k-1)\cdot 1 = k^{3} + 2k^{2} - 2k$ and $(k^{3} + k^{2} - 1) + (k-1)\cdot (k - 1) = k^{3} + 2k^{2} - 2k$, $\Delta(k^{3} + 2k^{2} - 2k) = 2$.

By Rule~2 applied to $(2k^{2} - 2k + 1) + (k-1)\cdot (k^{2} + k) = k^{3} + 2k^{2} - 3k + 1$ and $k^{2} + (k-1)\cdot (k - 1) = 2k^{2} - 2k + 1$, $\Delta(2k^{2} - 2k + 1) = 2$.

By Rule~2 applied to $k + (k-1)\cdot (2k - 1) = 2k^{2} - 2k + 1$ and $(2k - 1) + (k-1)\cdot (k - 1) = k^{2}$, $\Delta(2k - 1) = 1$.

By Rule~2 applied to $(2k - 1) + (k-1)\cdot 1 = 3k - 2$ and $(3k - 2) + (k-1)\cdot (k^{2} + 2k - 1) = k^{3} + k^{2} - 1$, $\Delta(3k - 2) = 2$.

By Rule~2 applied to $1 + (k-1)\cdot 2 = 2k - 1$ and $k + (k-1)\cdot 2 = 3k - 2$, $\Delta(2) = 3$.

By Rule~2 applied to $k^{3} + (k-1)\cdot 2k = k^{3} + 2k^{2} - 2k$ and $k\cdot 2 = 2k$, $\Delta(2k) = 1$.

By Rule~2 applied to $(k + 1) + (k-1)\cdot 1 = 2k$ and $(k^{3} + k^{2} - 2k + 1) + (k-1)\cdot (k + 1) = k^{3} + 2k^{2} - 2k$, $\Delta(k + 1) = 3$.

By Rule~2 applied to $(3k - 2) + (k-1)\cdot k = k^{2} + 2k - 2$ and $k^{2} + (k-1)\cdot 2 = k^{2} + 2k - 2$, $\Delta(k^{2} + 2k - 2) = 1$.

\emph{Contradiction (Rule~3) at $k^{2} + k - 2$:} if $\Delta(k^{2} + k - 2)=1$ then $(k^{2} + 2k - 2) + (k-1)\cdot (k^{2} + k - 2) = k^{3} + k^{2} - k$ is monochromatic; if $\Delta(k^{2} + k - 2)=2$ then $(3k - 2) + (k-1)\cdot (k^{2} + k - 2) = k^{3}$ is monochromatic; if $\Delta(k^{2} + k - 2)=3$ then $(k - 1) + (k-1)\cdot (k + 1) = k^{2} + k - 2$ is monochromatic.

\case{2.2.1.1.2.2}{$\Delta(k^{2} + k) = 3$}\quad {\small($A_1\supseteq \{1, k^{3} + k^{2} - k\},\ A_2\supseteq \{k, k^{3}, k^{3} + k^{2} - 2k + 1\},\ A_3\supseteq \{k^{2}, k^{2} + k, k^{3} + k^{2} - 1\}$)}

Since $(k^{2} + k - 1) + (k-1)\cdot (k^{2} + k) = k^{3} + k^{2} - 1$ must not be monochromatic in $A_3$, either $\Delta(k^{2} + k - 1) = 1$ or $\Delta(k^{2} + k - 1) = 2$.

\case{2.2.1.1.2.2.1}{$\Delta(k^{2} + k - 1) = 1$}\quad {\small($A_1\supseteq \{1, k^{2} + k - 1, k^{3} + k^{2} - k\},\ A_2\supseteq \{k, k^{3}, k^{3} + k^{2} - 2k + 1\},\ A_3\supseteq \{k^{2}, k^{2} + k, k^{3} + k^{2} - 1\}$)}

\emph{Contradiction:} $k\cdot (k^{2} + k - 1) = k^{3} + k^{2} - k$ is now monochromatic.

\case{2.2.1.1.2.2.2}{$\Delta(k^{2} + k - 1) = 2$}\quad {\small($A_1\supseteq \{1, k^{3} + k^{2} - k\},\ A_2\supseteq \{k, k^{2} + k - 1, k^{3}, k^{3} + k^{2} - 2k + 1\},\ A_3\supseteq \{k^{2}, k^{2} + k, k^{3} + k^{2} - 1\}$)}

By Rule~2 applied to $(k^{2} + k - 1) + (k-1)\cdot k = 2k^{2} - 1$ and $(2k^{2} - 1) + (k-1)\cdot k^{2} = k^{3} + k^{2} - 1$, $\Delta(2k^{2} - 1) = 1$.

By Rule~2 applied to $(2k^{2} - 1) + (k-1)\cdot (k^{2} - 1) = k^{3} + k^{2} - k$ and $(k^{2} + k - 1) + (k-1)\cdot (k^{2} - 1) = k^{3}$, $\Delta(k^{2} - 1) = 3$.

By Rule~2 applied to $k + (k-1)\cdot (k^{2} + k - 1) = k^{3} - k + 1$ and $k^{2} + (k-1)\cdot (k^{2} - 1) = k^{3} - k + 1$, $\Delta(k^{3} - k + 1) = 1$.

By Rule~2 applied to $(k^{3} - k + 1) + (k-1)\cdot (k + 1) = k^{3} + k^{2} - k$ and $k + (k-1)\cdot (k + 1) = k^{2} + k - 1$, $\Delta(k + 1) = 3$.

\emph{Contradiction:} $k\cdot (k + 1) = k^{2} + k$ is now monochromatic.

\case{2.2.1.2}{$\Delta(k^{3} + k^{2} - 2k + 1) = 3$}\quad {\small($A_1\supseteq \{1, k^{3} + k^{2} - k\},\ A_2\supseteq \{k, k^{3}\},\ A_3\supseteq \{k^{2}, k^{3} + k^{2} - 2k + 1\}$)}

Since $k^{2} + (k-1)\cdot (k^{2} + k - 1) = k^{3} + k^{2} - 2k + 1$ must not be monochromatic in $A_3$, either $\Delta(k^{2} + k - 1) = 1$ or $\Delta(k^{2} + k - 1) = 2$.

\case{2.2.1.2.1}{$\Delta(k^{2} + k - 1) = 1$}\quad {\small($A_1\supseteq \{1, k^{2} + k - 1, k^{3} + k^{2} - k\},\ A_2\supseteq \{k, k^{3}\},\ A_3\supseteq \{k^{2}, k^{3} + k^{2} - 2k + 1\}$)}

\emph{Contradiction:} $k\cdot (k^{2} + k - 1) = k^{3} + k^{2} - k$ is now monochromatic.

\case{2.2.1.2.2}{$\Delta(k^{2} + k - 1) = 2$}\quad {\small($A_1\supseteq \{1, k^{3} + k^{2} - k\},\ A_2\supseteq \{k, k^{2} + k - 1, k^{3}\},\ A_3\supseteq \{k^{2}, k^{3} + k^{2} - 2k + 1\}$)}

Since $k + (k-1)\cdot (k + 1) = k^{2} + k - 1$ must not be monochromatic in $A_2$, either $\Delta(k + 1) = 1$ or $\Delta(k + 1) = 3$.

\case{2.2.1.2.2.1}{$\Delta(k + 1) = 1$}\quad {\small($A_1\supseteq \{1, k + 1, k^{3} + k^{2} - k\},\ A_2\supseteq \{k, k^{2} + k - 1, k^{3}\},\ A_3\supseteq \{k^{2}, k^{3} + k^{2} - 2k + 1\}$)}

By Rule~2 applied to $k\cdot (k + 1) = k^{2} + k$ and $k + (k-1)\cdot (k^{2} + k) = k^{3}$, $\Delta(k^{2} + k) = 3$.

By Rule~2 applied to $(k^{3} - k + 1) + (k-1)\cdot (k + 1) = k^{3} + k^{2} - k$ and $k + (k-1)\cdot (k^{2} + k - 1) = k^{3} - k + 1$, $\Delta(k^{3} - k + 1) = 3$.

By Rule~2 applied to $(k^{2} + k - 1) + (k-1)\cdot (k^{2} - 1) = k^{3}$ and $k^{2} + (k-1)\cdot (k^{2} - 1) = k^{3} - k + 1$, $\Delta(k^{2} - 1) = 1$.

By Rule~2 applied to $(2k^{2} - 1) + (k-1)\cdot (k^{2} - 1) = k^{3} + k^{2} - k$ and $(k^{2} + k - 1) + (k-1)\cdot k = 2k^{2} - 1$, $\Delta(2k^{2} - 1) = 3$.

By Rule~2 applied to $(k^{3} + k^{2} - k) + (k-1)\cdot (k + 1) = k^{3} + 2k^{2} - k - 1$ and $(2k^{2} - 1) + (k-1)\cdot (k^{2} + k) = k^{3} + 2k^{2} - k - 1$, $\Delta(k^{3} + 2k^{2} - k - 1) = 2$.

\emph{Contradiction (Rule~3) at $k^{3} + k^{2} - 1$:} if $\Delta(k^{3} + k^{2} - 1)=1$ then $(k^{3} + k^{2} - k) + (k-1)\cdot 1 = k^{3} + k^{2} - 1$ is monochromatic; if $\Delta(k^{3} + k^{2} - 1)=2$ then $(k^{3} + k^{2} - 1) + (k-1)\cdot k = k^{3} + 2k^{2} - k - 1$ is monochromatic; if $\Delta(k^{3} + k^{2} - 1)=3$ then $(2k^{2} - 1) + (k-1)\cdot k^{2} = k^{3} + k^{2} - 1$ is monochromatic.

\case{2.2.1.2.2.2}{$\Delta(k + 1) = 3$}\quad {\small($A_1\supseteq \{1, k^{3} + k^{2} - k\},\ A_2\supseteq \{k, k^{2} + k - 1, k^{3}\},\ A_3\supseteq \{k + 1, k^{2}, k^{3} + k^{2} - 2k + 1\}$)}

By Rule~2 applied to $k + (k-1)\cdot (k^{2} + k) = k^{3}$ and $k\cdot (k + 1) = k^{2} + k$, $\Delta(k^{2} + k) = 1$.

By Rule~2 applied to $(k^{2} + k - 1) + (k-1)\cdot k = 2k^{2} - 1$ and $k^{2} + (k-1)\cdot (k + 1) = 2k^{2} - 1$, $\Delta(2k^{2} - 1) = 1$.

By Rule~2 applied to $(2k^{2} - 1) + (k-1)\cdot (k^{2} - 1) = k^{3} + k^{2} - k$ and $(k^{2} + k - 1) + (k-1)\cdot (k^{2} - 1) = k^{3}$, $\Delta(k^{2} - 1) = 3$.

\emph{Contradiction (Rule~3) at $k^{3} - k + 1$:} if $\Delta(k^{3} - k + 1)=1$ then $1 + (k-1)\cdot (k^{2} + k) = k^{3} - k + 1$ is monochromatic; if $\Delta(k^{3} - k + 1)=2$ then $k + (k-1)\cdot (k^{2} + k - 1) = k^{3} - k + 1$ is monochromatic; if $\Delta(k^{3} - k + 1)=3$ then $k^{2} + (k-1)\cdot (k^{2} - 1) = k^{3} - k + 1$ is monochromatic.

\case{2.2.2}{$\Delta(k^{3} + k^{2} - k) = 3$}\quad {\small($A_1\supseteq \{1\},\ A_2\supseteq \{k, k^{3}\},\ A_3\supseteq \{k^{2}, k^{3} + k^{2} - k\}$)}

By Rule~2 applied to $k + (k-1)\cdot (k^{2} + k) = k^{3}$ and $k^{2} + (k-1)\cdot (k^{2} + k) = k^{3} + k^{2} - k$, $\Delta(k^{2} + k) = 1$.

Since $1 + (k-1)\cdot (k^{2} + k) = k^{3} - k + 1$ must not be monochromatic in $A_1$, either $\Delta(k^{3} - k + 1) = 2$ or $\Delta(k^{3} - k + 1) = 3$.

\case{2.2.2.1}{$\Delta(k^{3} - k + 1) = 2$}\quad {\small($A_1\supseteq \{1, k^{2} + k\},\ A_2\supseteq \{k, k^{3} - k + 1, k^{3}\},\ A_3\supseteq \{k^{2}, k^{3} + k^{2} - k\}$)}

Since $k + (k-1)\cdot (k^{2} + k - 1) = k^{3} - k + 1$ must not be monochromatic in $A_2$, either $\Delta(k^{2} + k - 1) = 1$ or $\Delta(k^{2} + k - 1) = 3$.

\case{2.2.2.1.1}{$\Delta(k^{2} + k - 1) = 1$}\quad {\small($A_1\supseteq \{1, k^{2} + k - 1, k^{2} + k\},\ A_2\supseteq \{k, k^{3} - k + 1, k^{3}\},\ A_3\supseteq \{k^{2}, k^{3} + k^{2} - k\}$)}

Since $(k^{2} + k - 1) + (k-1)\cdot 1 = k^{2} + 2k - 2$ must not be monochromatic in $A_1$, either $\Delta(k^{2} + 2k - 2) = 2$ or $\Delta(k^{2} + 2k - 2) = 3$.

\case{2.2.2.1.1.1}{$\Delta(k^{2} + 2k - 2) = 2$}\quad {\small($A_1\supseteq \{1, k^{2} + k - 1, k^{2} + k\},\ A_2\supseteq \{k, k^{2} + 2k - 2, k^{3} - k + 1, k^{3}\},\ A_3\supseteq \{k^{2}, k^{3} + k^{2} - k\}$)}

By Rule~2 applied to $1 + (k-1)\cdot (k + 2) = k^{2} + k - 1$ and $k + (k-1)\cdot (k + 2) = k^{2} + 2k - 2$, $\Delta(k + 2) = 3$.

By Rule~2 applied to $(k^{2} + 2k - 2) + (k-1)\cdot k = 2k^{2} + k - 2$ and $k^{2} + (k-1)\cdot (k + 2) = 2k^{2} + k - 2$, $\Delta(2k^{2} + k - 2) = 1$.

By Rule~2 applied to $(2k^{2} + k - 2) + (k-1)\cdot (k^{2} + k) = k^{3} + 2k^{2} - 2$ and $(k^{3} + k^{2} - k) + (k-1)\cdot (k + 2) = k^{3} + 2k^{2} - 2$, $\Delta(k^{3} + 2k^{2} - 2) = 2$.

By Rule~2 applied to $(k^{2} + 2k - 2) + (k-1)\cdot (k^{2} + 2k) = k^{3} + 2k^{2} - 2$ and $k\cdot (k + 2) = k^{2} + 2k$, $\Delta(k^{2} + 2k) = 1$.

By Rule~2 applied to $1 + (k-1)\cdot (k^{2} + 2k) = k^{3} + k^{2} - 2k + 1$ and $(k^{3} - k + 1) + (k-1)\cdot k = k^{3} + k^{2} - 2k + 1$, $\Delta(k^{3} + k^{2} - 2k + 1) = 3$.

By Rule~2 applied to $(2k^{2} + k - 2) + (k-1)\cdot (k^{2} + k - 1) = k^{3} + 2k^{2} - k - 1$ and $(k^{3} + k^{2} - 2k + 1) + (k-1)\cdot (k + 2) = k^{3} + 2k^{2} - k - 1$, $\Delta(k^{3} + 2k^{2} - k - 1) = 2$.

By Rule~2 applied to $(k^{2} + k - 1) + (k-1)\cdot (k^{2} + k) = k^{3} + k^{2} - 1$ and $(k^{3} + k^{2} - 1) + (k-1)\cdot k = k^{3} + 2k^{2} - k - 1$, $\Delta(k^{3} + k^{2} - 1) = 3$.

By Rule~2 applied to $(2k^{2} - 1) + (k-1)\cdot 1 = 2k^{2} + k - 2$ and $(2k^{2} - 1) + (k-1)\cdot k^{2} = k^{3} + k^{2} - 1$, $\Delta(2k^{2} - 1) = 2$.

\emph{Contradiction (Rule~3) at $k^{2} + k + 1$:} if $\Delta(k^{2} + k + 1)=1$ then $(k^{2} + k + 1) + (k-1)\cdot 1 = k^{2} + 2k$ is monochromatic; if $\Delta(k^{2} + k + 1)=2$ then $(2k^{2} - 1) + (k-1)\cdot (k^{2} + k + 1) = k^{3} + 2k^{2} - 2$ is monochromatic; if $\Delta(k^{2} + k + 1)=3$ then $k^{2} + (k-1)\cdot (k^{2} + k + 1) = k^{3} + k^{2} - 1$ is monochromatic.

\case{2.2.2.1.1.2}{$\Delta(k^{2} + 2k - 2) = 3$}\quad {\small($A_1\supseteq \{1, k^{2} + k - 1, k^{2} + k\},\ A_2\supseteq \{k, k^{3} - k + 1, k^{3}\},\ A_3\supseteq \{k^{2}, k^{2} + 2k - 2, k^{3} + k^{2} - k\}$)}

Since $(k^{2} + 2k - 2) + (k-1)\cdot (k^{2} + k - 2) = k^{3} + k^{2} - k$ must not be monochromatic in $A_3$, either $\Delta(k^{2} + k - 2) = 1$ or $\Delta(k^{2} + k - 2) = 2$.

\case{2.2.2.1.1.2.1}{$\Delta(k^{2} + k - 2) = 1$}\quad {\small($A_1\supseteq \{1, k^{2} + k - 2, k^{2} + k - 1, k^{2} + k\},\ A_2\supseteq \{k, k^{3} - k + 1, k^{3}\},\ A_3\supseteq \{k^{2}, k^{2} + 2k - 2, k^{3} + k^{2} - k\}$)}

By Rule~2 applied to $(k^{2} + k - 1) + (k-1)\cdot (k^{2} + k - 2) = k^{3} + k^{2} - 2k + 1$ and $(k^{3} - k + 1) + (k-1)\cdot k = k^{3} + k^{2} - 2k + 1$, $\Delta(k^{3} + k^{2} - 2k + 1) = 3$.

Since $(3k - 2) + (k-1)\cdot (k^{2} + 2k - 2) = k^{3} + k^{2} - k$ must not be monochromatic in $A_3$, either $\Delta(3k - 2) = 1$ or $\Delta(3k - 2) = 2$.

\case{2.2.2.1.1.2.1.1}{$\Delta(3k - 2) = 1$}\quad {\small($A_1\supseteq \{1, 3k - 2, k^{2} + k - 2, k^{2} + k - 1, k^{2} + k\},\ A_2\supseteq \{k, k^{3} - k + 1, k^{3}\},\ A_3\supseteq \{k^{2}, k^{2} + 2k - 2, k^{3} + k^{2} - 2k + 1, k^{3} + k^{2} - k\}$)}

By Rule~2 applied to $(2k - 1) + (k-1)\cdot 1 = 3k - 2$ and $(2k - 1) + (k-1)\cdot (k^{2} + 2k - 2) = k^{3} + k^{2} - 2k + 1$, $\Delta(2k - 1) = 2$.

By Rule~2 applied to $k + (k-1)\cdot (2k - 1) = 2k^{2} - 2k + 1$ and $(2k^{2} - 2k + 1) + (k-1)\cdot k^{2} = k^{3} + k^{2} - 2k + 1$, $\Delta(2k^{2} - 2k + 1) = 1$.

By Rule~2 applied to $(2k^{2} - 2k + 1) + (k-1)\cdot (k^{2} + k) = k^{3} + 2k^{2} - 3k + 1$ and $k^{3} + (k-1)\cdot (2k - 1) = k^{3} + 2k^{2} - 3k + 1$, $\Delta(k^{3} + 2k^{2} - 3k + 1) = 3$.

By Rule~2 applied to $(3k - 2) + (k-1)\cdot (k - 1) = k^{2} + k - 1$ and $(k^{3} + k^{2} - k) + (k-1)\cdot (k - 1) = k^{3} + 2k^{2} - 3k + 1$, $\Delta(k - 1) = 2$.

By Rule~2 applied to $(k^{2} - 1) + (k-1)\cdot 1 = k^{2} + k - 2$ and $(k - 1) + (k-1)\cdot k = k^{2} - 1$, $\Delta(k^{2} - 1) = 3$.

By Rule~2 applied to $(3k - 2) + (k-1)\cdot (k^{2} + k - 1) = k^{3} + k - 1$ and $(k^{2} + 2k - 2) + (k-1)\cdot (k^{2} - 1) = k^{3} + k - 1$, $\Delta(k^{3} + k - 1) = 2$.

By Rule~2 applied to $(k^{3} + k - 1) + (k-1)\cdot (2k - 1) = k^{3} + 2k^{2} - 2k$ and $k\cdot (k^{2} + 2k - 2) = k^{3} + 2k^{2} - 2k$, $\Delta(k^{3} + 2k^{2} - 2k) = 1$.

By Rule~2 applied to $(k^{2} + k - 1) + (k-1)\cdot (k^{2} + k) = k^{3} + k^{2} - 1$ and $(k^{3} + k - 1) + (k-1)\cdot k = k^{3} + k^{2} - 1$, $\Delta(k^{3} + k^{2} - 1) = 3$.

\emph{Contradiction (Rule~3) at $k^{2} + k + 1$:} if $\Delta(k^{2} + k + 1)=1$ then $(2k^{2} - 2k + 1) + (k-1)\cdot (k^{2} + k + 1) = k^{3} + 2k^{2} - 2k$ is monochromatic; if $\Delta(k^{2} + k + 1)=2$ then $k + (k-1)\cdot (k^{2} + k + 1) = k^{3} + k - 1$ is monochromatic; if $\Delta(k^{2} + k + 1)=3$ then $k^{2} + (k-1)\cdot (k^{2} + k + 1) = k^{3} + k^{2} - 1$ is monochromatic.

\case{2.2.2.1.1.2.1.2}{$\Delta(3k - 2) = 2$}\quad {\small($A_1\supseteq \{1, k^{2} + k - 2, k^{2} + k - 1, k^{2} + k\},\ A_2\supseteq \{k, 3k - 2, k^{3} - k + 1, k^{3}\},\ A_3\supseteq \{k^{2}, k^{2} + 2k - 2, k^{3} + k^{2} - 2k + 1, k^{3} + k^{2} - k\}$)}

By Rule~2 applied to $k + (k-1)\cdot 2 = 3k - 2$ and $k^{2} + (k-1)\cdot 2 = k^{2} + 2k - 2$, $\Delta(2) = 1$.

By Rule~2 applied to $1 + (k-1)\cdot 2 = 2k - 1$ and $(2k - 1) + (k-1)\cdot (k^{2} + 2k - 2) = k^{3} + k^{2} - 2k + 1$, $\Delta(2k - 1) = 2$.

By Rule~2 applied to $k + (k-1)\cdot (2k - 1) = 2k^{2} - 2k + 1$ and $(2k^{2} - 2k + 1) + (k-1)\cdot k^{2} = k^{3} + k^{2} - 2k + 1$, $\Delta(2k^{2} - 2k + 1) = 1$.

By Rule~2 applied to $(2k^{2} - 2k + 1) + (k-1)\cdot 2 = 2k^{2} - 1$ and $(3k - 2) + (k-1)\cdot (2k - 1) = 2k^{2} - 1$, $\Delta(2k^{2} - 1) = 3$.

By Rule~2 applied to $(k^{2} + k - 1) + (k-1)\cdot (k^{2} + k) = k^{3} + k^{2} - 1$ and $(2k^{2} - 1) + (k-1)\cdot k^{2} = k^{3} + k^{2} - 1$, $\Delta(k^{3} + k^{2} - 1) = 2$.

\emph{Contradiction (Rule~3) at $k + 1$:} if $\Delta(k + 1)=1$ then $2 + (k-1)\cdot 1 = k + 1$ is monochromatic; if $\Delta(k + 1)=2$ then $k^{3} + (k-1)\cdot (k + 1) = k^{3} + k^{2} - 1$ is monochromatic; if $\Delta(k + 1)=3$ then $k^{2} + (k-1)\cdot (k + 1) = 2k^{2} - 1$ is monochromatic.

\case{2.2.2.1.1.2.2}{$\Delta(k^{2} + k - 2) = 2$}\quad {\small($A_1\supseteq \{1, k^{2} + k - 1, k^{2} + k\},\ A_2\supseteq \{k, k^{2} + k - 2, k^{3} - k + 1, k^{3}\},\ A_3\supseteq \{k^{2}, k^{2} + 2k - 2, k^{3} + k^{2} - k\}$)}

By Rule~2 applied to $(3k - 2) + (k-1)\cdot (k^{2} + k - 2) = k^{3}$ and $(3k - 2) + (k-1)\cdot (k^{2} + 2k - 2) = k^{3} + k^{2} - k$, $\Delta(3k - 2) = 1$.

By Rule~2 applied to $(2k - 1) + (k-1)\cdot 1 = 3k - 2$ and $(2k - 1) + (k-1)\cdot (k^{2} + k - 2) = k^{3} - k + 1$, $\Delta(2k - 1) = 3$.

By Rule~2 applied to $(3k - 2) + (k-1)\cdot (k - 1) = k^{2} + k - 1$ and $(2k - 1) + (k-1)\cdot (k - 1) = k^{2}$, $\Delta(k - 1) = 2$.

By Rule~2 applied to $(k - 1) + (k-1)\cdot (k + 1) = k^{2} + k - 2$ and $(2k - 1) + (k-1)\cdot (k + 1) = k^{2} + 2k - 2$, $\Delta(k + 1) = 1$.

\emph{Contradiction:} $k\cdot (k + 1) = k^{2} + k$ is now monochromatic.

\case{2.2.2.1.2}{$\Delta(k^{2} + k - 1) = 3$}\quad {\small($A_1\supseteq \{1, k^{2} + k\},\ A_2\supseteq \{k, k^{3} - k + 1, k^{3}\},\ A_3\supseteq \{k^{2}, k^{2} + k - 1, k^{3} + k^{2} - k\}$)}

\emph{Contradiction:} $k\cdot (k^{2} + k - 1) = k^{3} + k^{2} - k$ is now monochromatic.

\case{2.2.2.2}{$\Delta(k^{3} - k + 1) = 3$}\quad {\small($A_1\supseteq \{1, k^{2} + k\},\ A_2\supseteq \{k, k^{3}\},\ A_3\supseteq \{k^{2}, k^{3} - k + 1, k^{3} + k^{2} - k\}$)}

Since $(k^{3} - k + 1) + (k-1)\cdot (k + 1) = k^{3} + k^{2} - k$ must not be monochromatic in $A_3$, either $\Delta(k + 1) = 1$ or $\Delta(k + 1) = 2$.

\case{2.2.2.2.1}{$\Delta(k + 1) = 1$}\quad {\small($A_1\supseteq \{1, k + 1, k^{2} + k\},\ A_2\supseteq \{k, k^{3}\},\ A_3\supseteq \{k^{2}, k^{3} - k + 1, k^{3} + k^{2} - k\}$)}

\emph{Contradiction:} $k\cdot (k + 1) = k^{2} + k$ is now monochromatic.

\case{2.2.2.2.2}{$\Delta(k + 1) = 2$}\quad {\small($A_1\supseteq \{1, k^{2} + k\},\ A_2\supseteq \{k, k + 1, k^{3}\},\ A_3\supseteq \{k^{2}, k^{3} - k + 1, k^{3} + k^{2} - k\}$)}

Since $k + (k-1)\cdot (k + 1) = k^{2} + k - 1$ must not be monochromatic in $A_2$, either $\Delta(k^{2} + k - 1) = 1$ or $\Delta(k^{2} + k - 1) = 3$.

\case{2.2.2.2.2.1}{$\Delta(k^{2} + k - 1) = 1$}\quad {\small($A_1\supseteq \{1, k^{2} + k - 1, k^{2} + k\},\ A_2\supseteq \{k, k + 1, k^{3}\},\ A_3\supseteq \{k^{2}, k^{3} - k + 1, k^{3} + k^{2} - k\}$)}

By Rule~2 applied to $(k^{2} + k - 1) + (k-1)\cdot (k^{2} + k) = k^{3} + k^{2} - 1$ and $k^{3} + (k-1)\cdot (k + 1) = k^{3} + k^{2} - 1$, $\Delta(k^{3} + k^{2} - 1) = 3$.

By Rule~2 applied to $1 + (k-1)\cdot (k + 2) = k^{2} + k - 1$ and $(k^{3} - k + 1) + (k-1)\cdot (k + 2) = k^{3} + k^{2} - 1$, $\Delta(k + 2) = 2$.

By Rule~2 applied to $(k^{2} + k - 1) + (k-1)\cdot 1 = k^{2} + 2k - 2$ and $k + (k-1)\cdot (k + 2) = k^{2} + 2k - 2$, $\Delta(k^{2} + 2k - 2) = 3$.

By Rule~2 applied to $(k^{2} + k) + (k-1)\cdot 1 = k^{2} + 2k - 1$ and $(k + 1) + (k-1)\cdot (k + 2) = k^{2} + 2k - 1$, $\Delta(k^{2} + 2k - 1) = 3$.

By Rule~2 applied to $(k + 2) + (k-1)\cdot (k + 1) = k^{2} + k + 1$ and $k^{2} + (k-1)\cdot (k^{2} + k + 1) = k^{3} + k^{2} - 1$, $\Delta(k^{2} + k + 1) = 1$.

By Rule~2 applied to $(k^{2} + k + 1) + (k-1)\cdot 1 = k^{2} + 2k$ and $k\cdot (k + 2) = k^{2} + 2k$, $\Delta(k^{2} + 2k) = 3$.

Since $k^{2} + (k-1)\cdot 2 = k^{2} + 2k - 2$ must not be monochromatic in $A_3$, either $\Delta(2) = 1$ or $\Delta(2) = 2$.

\case{2.2.2.2.2.1.1}{$\Delta(2) = 1$}\quad {\small($A_1\supseteq \{1, 2, k^{2} + k - 1, k^{2} + k, k^{2} + k + 1\},\ A_2\supseteq \{k, k + 1, k + 2, k^{3}\},\ A_3\supseteq \{k^{2}, k^{2} + 2k - 2, k^{2} + 2k - 1, k^{2} + 2k, k^{3} - k + 1, k^{3} + k^{2} - k, k^{3} + k^{2} - 1\}$)}

By Rule~2 applied to $1 + (k-1)\cdot 2 = 2k - 1$ and $(2k - 1) + (k-1)\cdot (k^{2} + 2k) = k^{3} + k^{2} - 1$, $\Delta(2k - 1) = 2$.

By Rule~2 applied to $k^{3} + (k-1)\cdot (2k - 1) = k^{3} + 2k^{2} - 3k + 1$ and $k^{2} + (k-1)\cdot (k^{2} + 2k - 1) = k^{3} + 2k^{2} - 3k + 1$, $\Delta(k^{3} + 2k^{2} - 3k + 1) = 1$.

By Rule~2 applied to $(k^{3} + 2k^{2} - 3k + 1) + (k-1)\cdot 1 = k^{3} + 2k^{2} - 2k$ and $k^{2} + (k-1)\cdot (k^{2} + 2k) = k^{3} + 2k^{2} - 2k$, $\Delta(k^{3} + 2k^{2} - 2k) = 2$.

By Rule~2 applied to $(2k^{2} - 2k + 1) + (k-1)\cdot (k^{2} + k) = k^{3} + 2k^{2} - 3k + 1$ and $k + (k-1)\cdot (2k - 1) = 2k^{2} - 2k + 1$, $\Delta(2k^{2} - 2k + 1) = 3$.

By Rule~2 applied to $(k^{3} + k^{2} - 2k + 1) + (k-1)\cdot (k + 1) = k^{3} + 2k^{2} - 2k$ and $(2k^{2} - 2k + 1) + (k-1)\cdot k^{2} = k^{3} + k^{2} - 2k + 1$, $\Delta(k^{3} + k^{2} - 2k + 1) = 1$.

By Rule~2 applied to $(k^{2} + k - 1) + (k-1)\cdot (k^{2} + k - 2) = k^{3} + k^{2} - 2k + 1$ and $(k^{2} + 2k - 2) + (k-1)\cdot (k^{2} + k - 2) = k^{3} + k^{2} - k$, $\Delta(k^{2} + k - 2) = 2$.

By Rule~2 applied to $(3k - 2) + (k-1)\cdot (k^{2} + k - 2) = k^{3}$ and $(3k - 2) + (k-1)\cdot (k^{2} + 2k - 1) = k^{3} + k^{2} - 1$, $\Delta(3k - 2) = 1$.

\emph{Contradiction (Rule~3) at $k - 1$:} if $\Delta(k - 1)=1$ then $(3k - 2) + (k-1)\cdot (k - 1) = k^{2} + k - 1$ is monochromatic; if $\Delta(k - 1)=2$ then $(k - 1) + (k-1)\cdot (k + 1) = k^{2} + k - 2$ is monochromatic; if $\Delta(k - 1)=3$ then $k^{2} + (k-1)\cdot (k - 1) = 2k^{2} - 2k + 1$ is monochromatic.

\case{2.2.2.2.2.1.2}{$\Delta(2) = 2$}\quad {\small($A_1\supseteq \{1, k^{2} + k - 1, k^{2} + k, k^{2} + k + 1\},\ A_2\supseteq \{2, k, k + 1, k + 2, k^{3}\},\ A_3\supseteq \{k^{2}, k^{2} + 2k - 2, k^{2} + 2k - 1, k^{2} + 2k, k^{3} - k + 1, k^{3} + k^{2} - k, k^{3} + k^{2} - 1\}$)}

By Rule~2 applied to $k + (k-1)\cdot 2 = 3k - 2$ and $(3k - 2) + (k-1)\cdot (k^{2} + 2k - 1) = k^{3} + k^{2} - 1$, $\Delta(3k - 2) = 1$.

By Rule~2 applied to $(2k - 1) + (k-1)\cdot 1 = 3k - 2$ and $(2k - 1) + (k-1)\cdot (k^{2} + 2k) = k^{3} + k^{2} - 1$, $\Delta(2k - 1) = 2$.

By Rule~2 applied to $k^{3} + (k-1)\cdot (2k - 1) = k^{3} + 2k^{2} - 3k + 1$ and $k^{2} + (k-1)\cdot (k^{2} + 2k - 1) = k^{3} + 2k^{2} - 3k + 1$, $\Delta(k^{3} + 2k^{2} - 3k + 1) = 1$.

By Rule~2 applied to $(k^{3} + 2k^{2} - 3k + 1) + (k-1)\cdot 1 = k^{3} + 2k^{2} - 2k$ and $k^{2} + (k-1)\cdot (k^{2} + 2k) = k^{3} + 2k^{2} - 2k$, $\Delta(k^{3} + 2k^{2} - 2k) = 2$.

By Rule~2 applied to $(3k - 2) + (k-1)\cdot (k^{2} + k - 1) = k^{3} + k - 1$ and $(k^{3} + k - 1) + (k-1)\cdot (2k - 1) = k^{3} + 2k^{2} - 2k$, $\Delta(k^{3} + k - 1) = 3$.

By Rule~2 applied to $(3k - 2) + (k-1)\cdot (k - 1) = k^{2} + k - 1$ and $(k^{3} + k - 1) + (k-1)\cdot (k - 1) = k^{3} + k^{2} - k$, $\Delta(k - 1) = 2$.

By Rule~2 applied to $(k - 1) + (k-1)\cdot (k + 1) = k^{2} + k - 2$ and $(k^{2} + 2k - 2) + (k-1)\cdot (k^{2} + k - 2) = k^{3} + k^{2} - k$, $\Delta(k^{2} + k - 2) = 1$.

\emph{Contradiction (Rule~3) at $k^{2} - 1$:} if $\Delta(k^{2} - 1)=1$ then $(k^{2} - 1) + (k-1)\cdot 1 = k^{2} + k - 2$ is monochromatic; if $\Delta(k^{2} - 1)=2$ then $(k - 1) + (k-1)\cdot k = k^{2} - 1$ is monochromatic; if $\Delta(k^{2} - 1)=3$ then $k^{2} + (k-1)\cdot (k^{2} - 1) = k^{3} - k + 1$ is monochromatic.

\case{2.2.2.2.2.2}{$\Delta(k^{2} + k - 1) = 3$}\quad {\small($A_1\supseteq \{1, k^{2} + k\},\ A_2\supseteq \{k, k + 1, k^{3}\},\ A_3\supseteq \{k^{2}, k^{2} + k - 1, k^{3} - k + 1, k^{3} + k^{2} - k\}$)}

\emph{Contradiction:} $k\cdot (k^{2} + k - 1) = k^{3} + k^{2} - k$ is now monochromatic.

\par}

\section{The $144$ identities among the elements of $\Xmin$}\label{app:identities}
The identities are ordered by the value of $z$ at large $k$. Identities with $x=y$ are written $k\cdot x=kx$.

\begin{multicols}{2}\footnotesize\raggedright\noindent
1.~$k\cdot 1 = k$\\
2.~$2 + (k-1)\cdot 1 = k + 1$\\
3.~$k + (k-1)\cdot 1 = 2k - 1$\\
4.~$1 + (k-1)\cdot 2 = 2k - 1$\\
5.~$(k + 1) + (k-1)\cdot 1 = 2k$\\
6.~$k\cdot 2 = 2k$\\
7.~$(2k - 1) + (k-1)\cdot 1 = 3k - 2$\\
8.~$k + (k-1)\cdot 2 = 3k - 2$\\
9.~$(k - 1) + (k-1)\cdot k = k^{2} - 1$\\
10.~$(2k - 1) + (k-1)\cdot (k - 1) = k^{2}$\\
11.~$k\cdot k = k^{2}$\\
12.~$1 + (k-1)\cdot (k + 1) = k^{2}$\\
13.~$(k^{2} - 1) + (k-1)\cdot 1 = k^{2} + k - 2$\\
14.~$(k - 1) + (k-1)\cdot (k + 1) = k^{2} + k - 2$\\
15.~$k^{2} + (k-1)\cdot 1 = k^{2} + k - 1$\\
16.~$(3k - 2) + (k-1)\cdot (k - 1) = k^{2} + k - 1$\\
17.~$(2k - 1) + (k-1)\cdot k = k^{2} + k - 1$\\
18.~$k + (k-1)\cdot (k + 1) = k^{2} + k - 1$\\
19.~$1 + (k-1)\cdot (k + 2) = k^{2} + k - 1$\\
20.~$2k + (k-1)\cdot k = k^{2} + k$\\
21.~$k\cdot (k + 1) = k^{2} + k$\\
22.~$2 + (k-1)\cdot (k + 2) = k^{2} + k$\\
23.~$(k + 2) + (k-1)\cdot (k + 1) = k^{2} + k + 1$\\
24.~$(k^{2} + k - 1) + (k-1)\cdot 1 = k^{2} + 2k - 2$\\
25.~$k^{2} + (k-1)\cdot 2 = k^{2} + 2k - 2$\\
26.~$(3k - 2) + (k-1)\cdot k = k^{2} + 2k - 2$\\
27.~$(2k - 1) + (k-1)\cdot (k + 1) = k^{2} + 2k - 2$\\
28.~$k + (k-1)\cdot (k + 2) = k^{2} + 2k - 2$\\
29.~$(k^{2} + k) + (k-1)\cdot 1 = k^{2} + 2k - 1$\\
30.~$2k + (k-1)\cdot (k + 1) = k^{2} + 2k - 1$\\
31.~$(k + 1) + (k-1)\cdot (k + 2) = k^{2} + 2k - 1$\\
32.~$(k^{2} + k + 1) + (k-1)\cdot 1 = k^{2} + 2k$\\
33.~$k\cdot (k + 2) = k^{2} + 2k$\\
34.~$k^{2} + (k-1)\cdot (k - 1) = 2k^{2} - 2k + 1$\\
35.~$k + (k-1)\cdot (2k - 1) = 2k^{2} - 2k + 1$\\
36.~$1 + (k-1)\cdot 2k = 2k^{2} - 2k + 1$\\
37.~$(2k^{2} - 2k + 1) + (k-1)\cdot 2 = 2k^{2} - 1$\\
38.~$(k^{2} + 2k - 2) + (k-1)\cdot (k - 1) = 2k^{2} - 1$\\
39.~$(k^{2} + k - 1) + (k-1)\cdot k = 2k^{2} - 1$\\
40.~$k^{2} + (k-1)\cdot (k + 1) = 2k^{2} - 1$\\
41.~$(3k - 2) + (k-1)\cdot (2k - 1) = 2k^{2} - 1$\\
42.~$(2k - 1) + (k-1)\cdot 2k = 2k^{2} - 1$\\
43.~$(2k^{2} - 1) + (k-1)\cdot 1 = 2k^{2} + k - 2$\\
44.~$(k^{2} + 2k - 2) + (k-1)\cdot k = 2k^{2} + k - 2$\\
45.~$(k^{2} + k - 1) + (k-1)\cdot (k + 1) = 2k^{2} + k - 2$\\
46.~$k^{2} + (k-1)\cdot (k + 2) = 2k^{2} + k - 2$\\
47.~$(3k - 2) + (k-1)\cdot 2k = 2k^{2} + k - 2$\\
48.~$(2k - 1) + (k-1)\cdot (k^{2} - 2) = k^{3} - k^{2} + 1$\\
49.~$k + (k-1)\cdot (k^{2} - 1) = k^{3} - k^{2} + 1$\\
50.~$1 + (k-1)\cdot k^{2} = k^{3} - k^{2} + 1$\\
51.~$(k^{3} - k^{2} + 1) + (k-1)\cdot k = k^{3} - k + 1$\\
52.~$(k^{2} + k - 1) + (k-1)\cdot (k^{2} - 2) = k^{3} - k + 1$\\
53.~$k^{2} + (k-1)\cdot (k^{2} - 1) = k^{3} - k + 1$\\
54.~$(2k - 1) + (k-1)\cdot (k^{2} + k - 2) = k^{3} - k + 1$\\
55.~$k + (k-1)\cdot (k^{2} + k - 1) = k^{3} - k + 1$\\
56.~$1 + (k-1)\cdot (k^{2} + k) = k^{3} - k + 1$\\
57.~$(k^{3} - k + 1) + (k-1)\cdot 1 = k^{3}$\\
58.~$(k^{3} - k^{2} + 1) + (k-1)\cdot (k + 1) = k^{3}$\\
59.~$(k^{2} + 2k - 2) + (k-1)\cdot (k^{2} - 2) = k^{3}$\\
60.~$(k^{2} + k - 1) + (k-1)\cdot (k^{2} - 1) = k^{3}$\\
61.~$k\cdot k^{2} = k^{3}$\\
62.~$(3k - 2) + (k-1)\cdot (k^{2} + k - 2) = k^{3}$\\
63.~$(2k - 1) + (k-1)\cdot (k^{2} + k - 1) = k^{3}$\\
64.~$k + (k-1)\cdot (k^{2} + k) = k^{3}$\\
65.~$1 + (k-1)\cdot (k^{2} + k + 1) = k^{3}$\\
66.~$k^{3} + (k-1)\cdot 1 = k^{3} + k - 1$\\
67.~$(k^{3} - k + 1) + (k-1)\cdot 2 = k^{3} + k - 1$\\
68.~$(k^{3} - k^{2} + 1) + (k-1)\cdot (k + 2) = k^{3} + k - 1$\\
69.~$(k^{2} + 2k - 2) + (k-1)\cdot (k^{2} - 1) = k^{3} + k - 1$\\
70.~$(k^{2} + k - 1) + (k-1)\cdot k^{2} = k^{3} + k - 1$\\
71.~$(3k - 2) + (k-1)\cdot (k^{2} + k - 1) = k^{3} + k - 1$\\
72.~$(2k - 1) + (k-1)\cdot (k^{2} + k) = k^{3} + k - 1$\\
73.~$k + (k-1)\cdot (k^{2} + k + 1) = k^{3} + k - 1$\\
74.~$k^{3} + (k-1)\cdot (k - 1) = k^{3} + k^{2} - 2k + 1$\\
75.~$(k^{3} - k + 1) + (k-1)\cdot k = k^{3} + k^{2} - 2k + 1$\\
76.~$(k^{3} - k^{2} + 1) + (k-1)\cdot 2k = k^{3} + k^{2} - 2k + 1$\\
77.~$(2k^{2} - 1) + (k-1)\cdot (k^{2} - 2) = k^{3} + k^{2} - 2k + 1$\\
78.~$(2k^{2} - 2k + 1) + (k-1)\cdot k^{2} = k^{3} + k^{2} - 2k + 1$\\
79.~$(k^{2} + k - 1) + (k-1)\cdot (k^{2} + k - 2) = k^{3} + k^{2} - 2k + 1$\\
80.~$k^{2} + (k-1)\cdot (k^{2} + k - 1) = k^{3} + k^{2} - 2k + 1$\\
81.~$(2k - 1) + (k-1)\cdot (k^{2} + 2k - 2) = k^{3} + k^{2} - 2k + 1$\\
82.~$k + (k-1)\cdot (k^{2} + 2k - 1) = k^{3} + k^{2} - 2k + 1$\\
83.~$1 + (k-1)\cdot (k^{2} + 2k) = k^{3} + k^{2} - 2k + 1$\\
84.~$(k^{3} + k^{2} - 2k + 1) + (k-1)\cdot 1 = k^{3} + k^{2} - k$\\
85.~$(k^{3} + k - 1) + (k-1)\cdot (k - 1) = k^{3} + k^{2} - k$\\
86.~$k^{3} + (k-1)\cdot k = k^{3} + k^{2} - k$\\
87.~$(k^{3} - k + 1) + (k-1)\cdot (k + 1) = k^{3} + k^{2} - k$\\
88.~$(2k^{2} + k - 2) + (k-1)\cdot (k^{2} - 2) = k^{3} + k^{2} - k$\\
89.~$(2k^{2} - 1) + (k-1)\cdot (k^{2} - 1) = k^{3} + k^{2} - k$\\
90.~$(k^{2} + 2k - 2) + (k-1)\cdot (k^{2} + k - 2) = k^{3} + k^{2} - k$\\
91.~$k\cdot (k^{2} + k - 1) = k^{3} + k^{2} - k$\\
92.~$k^{2} + (k-1)\cdot (k^{2} + k) = k^{3} + k^{2} - k$\\
93.~$(3k - 2) + (k-1)\cdot (k^{2} + 2k - 2) = k^{3} + k^{2} - k$\\
94.~$(2k - 1) + (k-1)\cdot (k^{2} + 2k - 1) = k^{3} + k^{2} - k$\\
95.~$k + (k-1)\cdot (k^{2} + 2k) = k^{3} + k^{2} - k$\\
96.~$(k^{3} + k^{2} - k) + (k-1)\cdot 1 = k^{3} + k^{2} - 1$\\
97.~$(k^{3} + k^{2} - 2k + 1) + (k-1)\cdot 2 = k^{3} + k^{2} - 1$\\
98.~$(k^{3} + k - 1) + (k-1)\cdot k = k^{3} + k^{2} - 1$\\
99.~$k^{3} + (k-1)\cdot (k + 1) = k^{3} + k^{2} - 1$\\
100.~$(k^{3} - k + 1) + (k-1)\cdot (k + 2) = k^{3} + k^{2} - 1$\\
101.~$(2k^{2} + k - 2) + (k-1)\cdot (k^{2} - 1) = k^{3} + k^{2} - 1$\\
102.~$(2k^{2} - 1) + (k-1)\cdot k^{2} = k^{3} + k^{2} - 1$\\
103.~$(k^{2} + 2k - 2) + (k-1)\cdot (k^{2} + k - 1) = k^{3} + k^{2} - 1$\\
104.~$(k^{2} + k - 1) + (k-1)\cdot (k^{2} + k) = k^{3} + k^{2} - 1$\\
105.~$k^{2} + (k-1)\cdot (k^{2} + k + 1) = k^{3} + k^{2} - 1$\\
106.~$(3k - 2) + (k-1)\cdot (k^{2} + 2k - 1) = k^{3} + k^{2} - 1$\\
107.~$(2k - 1) + (k-1)\cdot (k^{2} + 2k) = k^{3} + k^{2} - 1$\\
108.~$(k^{3} + k^{2} - k) + (k-1)\cdot (k - 1) = k^{3} + 2k^{2} - 3k + 1$\\
109.~$(k^{3} + k^{2} - 2k + 1) + (k-1)\cdot k = k^{3} + 2k^{2} - 3k + 1$\\
110.~$k^{3} + (k-1)\cdot (2k - 1) = k^{3} + 2k^{2} - 3k + 1$\\
111.~$(k^{3} - k + 1) + (k-1)\cdot 2k = k^{3} + 2k^{2} - 3k + 1$\\
112.~$(2k^{2} - 1) + (k-1)\cdot (k^{2} + k - 2) = k^{3} + 2k^{2} - 3k + 1$\\
113.~$(2k^{2} - 2k + 1) + (k-1)\cdot (k^{2} + k) = k^{3} + 2k^{2} - 3k + 1$\\
114.~$(k^{2} + k - 1) + (k-1)\cdot (k^{2} + 2k - 2) = k^{3} + 2k^{2} - 3k + 1$\\
115.~$k^{2} + (k-1)\cdot (k^{2} + 2k - 1) = k^{3} + 2k^{2} - 3k + 1$\\
116.~$(k^{3} + 2k^{2} - 3k + 1) + (k-1)\cdot 1 = k^{3} + 2k^{2} - 2k$\\
117.~$(k^{3} + k^{2} - 1) + (k-1)\cdot (k - 1) = k^{3} + 2k^{2} - 2k$\\
118.~$(k^{3} + k^{2} - k) + (k-1)\cdot k = k^{3} + 2k^{2} - 2k$\\
119.~$(k^{3} + k^{2} - 2k + 1) + (k-1)\cdot (k + 1) = k^{3} + 2k^{2} - 2k$\\
120.~$(k^{3} + k - 1) + (k-1)\cdot (2k - 1) = k^{3} + 2k^{2} - 2k$\\
121.~$k^{3} + (k-1)\cdot 2k = k^{3} + 2k^{2} - 2k$\\
122.~$(2k^{2} + k - 2) + (k-1)\cdot (k^{2} + k - 2) = k^{3} + 2k^{2} - 2k$\\
123.~$(2k^{2} - 1) + (k-1)\cdot (k^{2} + k - 1) = k^{3} + 2k^{2} - 2k$\\
124.~$(2k^{2} - 2k + 1) + (k-1)\cdot (k^{2} + k + 1) = k^{3} + 2k^{2} - 2k$\\
125.~$k\cdot (k^{2} + 2k - 2) = k^{3} + 2k^{2} - 2k$\\
126.~$(k^{2} + k - 1) + (k-1)\cdot (k^{2} + 2k - 1) = k^{3} + 2k^{2} - 2k$\\
127.~$k^{2} + (k-1)\cdot (k^{2} + 2k) = k^{3} + 2k^{2} - 2k$\\
128.~$(k^{3} + 2k^{2} - 2k) + (k-1)\cdot 1 = k^{3} + 2k^{2} - k - 1$\\
129.~$(k^{3} + 2k^{2} - 3k + 1) + (k-1)\cdot 2 = k^{3} + 2k^{2} - k - 1$\\
130.~$(k^{3} + k^{2} - 1) + (k-1)\cdot k = k^{3} + 2k^{2} - k - 1$\\
131.~$(k^{3} + k^{2} - k) + (k-1)\cdot (k + 1) = k^{3} + 2k^{2} - k - 1$\\
132.~$(k^{3} + k^{2} - 2k + 1) + (k-1)\cdot (k + 2) = k^{3} + 2k^{2} - k - 1$\\
133.~$(k^{3} + k - 1) + (k-1)\cdot 2k = k^{3} + 2k^{2} - k - 1$\\
134.~$(2k^{2} + k - 2) + (k-1)\cdot (k^{2} + k - 1) = k^{3} + 2k^{2} - k - 1$\\
135.~$(2k^{2} - 1) + (k-1)\cdot (k^{2} + k) = k^{3} + 2k^{2} - k - 1$\\
136.~$(k^{2} + 2k - 2) + (k-1)\cdot (k^{2} + 2k - 1) = k^{3} + 2k^{2} - k - 1$\\
137.~$(k^{2} + k - 1) + (k-1)\cdot (k^{2} + 2k) = k^{3} + 2k^{2} - k - 1$\\
138.~$(k^{3} + 2k^{2} - k - 1) + (k-1)\cdot 1 = k^{3} + 2k^{2} - 2$\\
139.~$(k^{3} + 2k^{2} - 2k) + (k-1)\cdot 2 = k^{3} + 2k^{2} - 2$\\
140.~$(k^{3} + k^{2} - 1) + (k-1)\cdot (k + 1) = k^{3} + 2k^{2} - 2$\\
141.~$(k^{3} + k^{2} - k) + (k-1)\cdot (k + 2) = k^{3} + 2k^{2} - 2$\\
142.~$(2k^{2} + k - 2) + (k-1)\cdot (k^{2} + k) = k^{3} + 2k^{2} - 2$\\
143.~$(2k^{2} - 1) + (k-1)\cdot (k^{2} + k + 1) = k^{3} + 2k^{2} - 2$\\
144.~$(k^{2} + 2k - 2) + (k-1)\cdot (k^{2} + 2k) = k^{3} + 2k^{2} - 2$
\end{multicols}

\section{Witnesses of minimality}\label{app:witness}
For each $p\in\Xmin$ the row below gives a good map of $\Xmin\setminus\{p\}$ as a string of $33$ colors. The elements are in the order of Table~\ref{tab:X}, grouped as $9$ linear, $13$ quadratic and $11$ cubic, and the removed element is marked ``\texttt{-}''. Each map avoids all identities among the remaining $32$ elements, which proves Theorem~\ref{thm:critical}(4).

\begin{longtable}{ll}
\toprule removed $p$ & colors of the other $32$ elements\\ \midrule\endhead
\bottomrule\endfoot
$1$ & \texttt{-22211111\,2112222233333\,33333333333} \\
$2$ & \texttt{1-2222221\,1113333111112\,33332223333} \\
$k - 1$ & \texttt{13-232113\,2112323321221\,23313121233} \\
$k$ & \texttt{121-22212\,2212211133333\,33333333333} \\
$k + 1$ & \texttt{1222-1113\,2112222223333\,33333333333} \\
$k + 2$ & \texttt{12122-313\,1232112222213\,33213333333} \\
$2k - 1$ & \texttt{111223-21\,1113333211221\,33312223333} \\
$2k$ & \texttt{1212111-3\,2232222222213\,33133333333} \\
$3k - 2$ & \texttt{12123321-\,1112322111322\,23332213333} \\
$k^{2} - 2$ & \texttt{111222323\,-132333111222\,13122213333} \\
$k^{2} - 1$ & \texttt{111223323\,1-31333211121\,13122213333} \\
$k^{2}$ & \texttt{122211333\,13-2232122112\,21213333333} \\
$k^{2} + k - 2$ & \texttt{121222211\,111-333311322\,23333122233} \\
$k^{2} + k - 1$ & \texttt{121221211\,2212-12133333\,33333333333} \\
$k^{2} + k$ & \texttt{122211331\,11322-2222213\,21213333333} \\
$k^{2} + k + 1$ & \texttt{121211331\,113222-222113\,23123331333} \\
$k^{2} + 2k - 2$ & \texttt{111223323\,1132333-11222\,33122213333} \\
$k^{2} + 2k - 1$ & \texttt{121221213\,11333111-2322\,13222133333} \\
$k^{2} + 2k$ & \texttt{121231211\,111232211-322\,23331223333} \\
$2k^{2} - 2k + 1$ & \texttt{121222213\,1112333111-22\,33323213333} \\
$2k^{2} - 1$ & \texttt{121222313\,22313331112-1\,13122213333} \\
$2k^{2} + k - 2$ & \texttt{131222211\,121233331131-\,31333122233} \\
$k^{3} - k^{2} + 1$ & \texttt{111223221\,2112332331331\,-3333122213} \\
$k^{3} - k + 1$ & \texttt{122211333\,3333232122212\,2-213133333} \\
$k^{3}$ & \texttt{122211331\,2132222222213\,21-13333333} \\
$k^{3} + k - 1$ & \texttt{121233213\,1112322111322\,233-2213333} \\
$k^{3} + k^{2} - 2k + 1$ & \texttt{111223333\,3333311222321\,1222-133333} \\
$k^{3} + k^{2} - k$ & \texttt{121221211\,2212312122322\,33333-33333} \\
$k^{3} + k^{2} - 1$ & \texttt{122231113\,3311222232233\,233121-3333} \\
$k^{3} + 2k^{2} - 3k + 1$ & \texttt{131221311\,1132311232223\,1322213-123} \\
$k^{3} + 2k^{2} - 2k$ & \texttt{123222211\,1213313131312\,33322233-33} \\
$k^{3} + 2k^{2} - k - 1$ & \texttt{121231211\,1112321113322\,233312233-3} \\
$k^{3} + 2k^{2} - 2$ & \texttt{121223333\,3333112222121\,1221133333-} \\
\end{longtable}

\end{document}